\documentclass[a4paper,psfig,10pt]{amsart}
\usepackage{epic, eepic,amssymb}

\input xy
\xyoption{all}

\theoremstyle{definition}
\newtheorem{thm}{Theorem}[section]
\newtheorem{lem}{Lemma}[section]
\newtheorem{defn}{Definition}[section]

\newtheorem{prop}{Proposition}[section]

\DeclareMathOperator{\ohom}{hom}

\DeclareMathOperator{\oAut}{Aut}

\DeclareMathOperator{\odiff}{Diffeo^+}
\DeclareMathOperator{\oid}{id}

\def\bR{\mathbb{R}}

\def\bN{\mathbb{N}}

\def\bf1{\mathbf{1}}

\def\cA{\mathcal{A}}
\def\cB{\mathcal{B}}

\def\cS{\mathcal{S}}
\def\cD{\mathcal{D}}
\def\cI{\mathcal{I}}
\def\cJ{\mathcal{J}}
\def\cE{\mathcal{E}}

\def\la{\langle}
\def\ra{\rangle}

\def\va{\varphi(\alpha)}
\def\vg{\varphi(\gamma)}
\def\vai{\varphi(\alpha^{-1})}
\def\vgi{\varphi(\gamma^{-1})}
\def\rg{\rho(g)}
\def\rgi{\rho(g^{-1})}

\newcommand{\gap}{\\ [1.5mm]}

\newcommand{\nattrans}{\stackrel{\bullet}{\longrightarrow}}

\newcommand{\additive}{\widehat{{{k\!\!-\!\!\mbox{Add}}}}}
\newcommand{\monunit}{\bf1}
\newcommand{\aub}{{\underline{\alpha}}}
\newcommand{\aubi}{{\underline{\alpha}^{-1}}}

\newcommand{\rotate}[1] {\widetilde{#1}}
\newcommand{\xyreflect}[1] {\widehat{#1}}

\title[Homotopy Quantum Field Theories and Tortile Structures]{Homotopy Quantum Field Theories and \\ Tortile Structures}
\author{Mark Brightwell}
\address{University of Aarhus\\
8000 \AA rhus\\
Denmark}
\author{Paul Turner}
\address{Heriot-Watt University\\ Edinburgh EH14 4AS\\Scotland}

\begin{document}
\begin{abstract}
We study a variation of Turaev's homotopy quantum field theories using 2-categories
of surfaces. We define the homotopy surface 2-category of a space $X$ and
define an $\cS_X$-structure to be a monoidal 2-functor from this to the
2-category of idempotent-complete additive $k$-linear categories. We
initiate the study of the algebraic structure arising from these functors. In particular we show that, under certain conditions, an
$\cS_X$-structure gives rise to a lax tortile $\pi$-category when the
background space is an Eilenberg-Maclane space $X=K(\pi,1)$, and to a
tortile category with lax $\pi_2X$-action when the background space is
simply-connected.
\end{abstract}
\maketitle
\section*{Introduction}
The motivation for this paper was to construct approximations to a
conformal version of homotopy quantum field theory using 2-categories.
A homotopy quantum field theory, as defined by Turaev in
\cite{Turaev:HomotopyFieldTheory2D}, is a variant of a topological
quantum field theory in which manifolds come equipped with a map to
some auxiliary space $X$. From a geometrical point of view a 1+1
dimensional homotopy quantum field theory is something like a vector
bundle on the free loop space of $X$ with a generalised flat
connection, giving ``parallel transport'' across surfaces. The
definition can be formulated in terms of representations of categories
of cobordisms in $X$ \cite{BrightwellTurner:Representations,
  Rodrigues:HomotopyQFT} and in 1+1 dimensions classification theorems
in terms of generalised Frobenius algebras are possible
\cite{Turaev:HomotopyFieldTheory2D, BrightwellTurner:Representations}.
From a string
theory point of view a conformal version of this set up is required
and Segal \cite{Segal:Elliptic} has defined a category whose
representations would provide this. In the background free case the
category is that of Riemann surfaces appearing in the definition of
conformal field theory. Tillmann \cite{Tillmann:Discrete,
  Tillmann:S-Structures} has pioneered the use of 2-categories to
approximate this category, by replacing the morphism spaces (of
Riemann surfaces) with categories whose classifying spaces have the
same rational homotopy type. The relevant 2-category is one whose
objects are circles, whose morphisms are surfaces and whose
2-morphisms are (path components of) diffeomorphisms of surfaces, and its 2-representations are closely
related to conformal field theory.  In this paper, we generalise
Tillmann's work to study structure arising from representations of
2-categories where a background space is incorporated which provides
an approximation to a conformal variant of 1+1-dimensional homotopy
quantum field theory. For a space $X$, we call such a representation an $\cS_X$-structure. We restrict ourselves to examining the genus zero part and show how the additive categories arising have
a rich structure inherited from the underlying geometry.  In order to
have rigid duality in these categories we require further assumptions
and we discuss how certain self-dual representations give rise to
this.

In 2+1 dimensions Turaev has produced examples of homotopy quantum field theories where the background
space is an Eilenberg-Maclane space, by generalising the definition of
a modular category \cite{Turaev:HomotopyFieldTheory3D}. He similarly introduced the notion of homotopy modular functor. As in the background free case, $\cS_X$-structures are closely related to both these notions (more details are given in section \ref{sec:structures}). 
And indeed the categorical structures we obtain are very close to those used by Turaev \cite{Turaev:HomotopyFieldTheory3D} (or to be more precise, close to the genus zero part, which is all we consider in this paper). We get lax versions of his structure, and
more importantly rigid duality is not present from the outset.

In more detail the paper contains the following. In section \ref{sec:2-category} we start
by constructing a 2-category $\cS_X$ whose objects are circles mapped
into $X$, whose morphisms are surfaces (with boundary) mapped to $X$
(with deformation up to homotopy), and whose 2-morphisms are path components of diffeomorphisms between these. This model is based on the construction
of Tillmann \cite{Tillmann:Discrete} in the background free case. We also
introduce two operators on the 2-category, reflection and rotation,
which play a central role in the geometric arguments later in the
paper. In section \ref{sec:structures} we define an $\cS_X$-structure as a 2-functor from $\cS_X$ to the 2-category of
idempotent complete additive categories over an algebraically closed
field $k$ and discuss dual $\cS_X$-structures. In section \ref{section:balanced} we
associate a category to an $\cS_X$-structure and prove the following theorem

\setcounter{section}{3}
\begin{thm}
(a) Let $\pi$ be a discrete group and let $X$ be an
  Eilenberg-Maclane space $K(\pi,1)$. The $k$-additive category $\cA$
  associated to an $\cS_X$-structure is a balanced $\pi$-category.  \\
(b) For any space $X$, the subcategory $\cA_1$ is a balanced
  category with $\pi_2X$-action.\\
(c) The categories above are semi-simple Artinian categories.
\end{thm}

Definitions of the categorical structures appearing in the above theorem can be found in Appendix \ref{app:defns}. Finally, in section \ref{section:tortile} we consider self dual theories and prove that
the associated category in this case has rigid duals and thus a tortile structure:
\setcounter{section}{4}
\setcounter{thm}{0}
\begin{thm}
(a) Let $\pi$ be a discrete group and let $X$ be an
  Eilenberg-Maclane space $K(\pi,1)$. The $k$-additive category $\cA$
  associated to an $\cS_X$-structure which is lax self dual with respect
  to hom,  is a semi-simple Artinian lax tortile $\pi$-category.  \\
(b) For any space $X$, the subcategory $\cA_1$ is a 
semi-simple Artinian tortile category with lax $\pi_2X$-action.
\end{thm}

As already noted, we gather together the definitions of
appropriate variants of balanced and tortile categories in an appendix,
and a further appendix recalls the basics of additive categories over
$k$ (or $k$-additive categories) and Tillmann's involution.

\setcounter{section}{0}
\section{The homotopy surface 2-category of a space}
\label{sec:2-category}

We begin with a few recollections about 2-categories mainly to
establish terminology.  Recall that
a \emph{2-category} $\cB$ is essentially a category in which the
morphism sets are categories and composition $\cB(A,B)\times
\cB(B,C)\rightarrow \cB(A,C)$ is functorial. The morphisms of the
morphism categories are called 2-morphisms. We shall denote objects by
$A,B,C,\dots$ 1-morphisms by $f,g,h,\dots$ and 2-morphisms by
$\alpha,\beta,\gamma,\dots$. By $\cB_{A,B}(f,g)$ we mean the set of
2-morphisms between 1-morphisms $f,g\in \cB(A,B)$. 2-morphisms have two kinds of compositions $\circ _1$
and $\circ_2$ called vertical and horizontal composition:
\begin{eqnarray*}
&&\circ_1:\cB_{A,B}(f,g)\times \cB_{A,B}(g,h) \rightarrow \cB_{A,B}(f,h)\\
&&\circ_2:\cB_{A,B}(f_1,g_1)\times \cB_{B,C}(f_2.g_2) \rightarrow \cB_{A,C}(f_1f_2,g_1g_2)
\end{eqnarray*}
$\cB$ comes equipped with associativity and identity 2-isomorphisms
$a_{fgh}$, $l_f$ and $r_f$ for 1-morphisms $f,g,h$:
\begin{eqnarray*}
&&a_{fgh}:(fg)h\rightarrow f(gh)\\
&&l_f:1_Af\rightarrow f\\
&&r_f:f1_B\rightarrow f
\end{eqnarray*}
where $f$ is in $\cB(A,B)$, and satisfying the associativity pentagons
and identity triangles. A
\emph{strict} 2-category is one in which all associativity and
identity 2-isomorphisms are identities. 
If  $G,H:\cD \rightarrow \cE$ are 2-functors between strict 2-categories $\cD$ and $\cE$ then a  {\em pseudo
2-natural transformation} $G\rightarrow H$ is a collection of
1-morphisms $N_U: GU\rightarrow HU$ and 2-isomorphisms
\[
\xymatrix{GU\ar[d]_{Gf}\ar[rr]^{N_U} &
{\raisebox{-2ex}{$~$}}\ar@{=>}[d]^{Nf} & HU \ar[d]^{Hf}\\ 
GV\ar[rr]_{N_V} & {\raisebox{2ex}{$~$}}& HV
}
\]
for objects $U,V\in \cD$ and morphism $f:U\rightarrow V$. These must satisfy
\[
N(id_U)=id_U,\quad N(fg)=N(f)N(g)
\]
and 
\[
\xymatrix{
& GU\ar@/_1.5pc/[dd]_{Gh} \ar@/^1.5pc/[dd]^{Gf} \ar[rr]^{N_U} & &
HU\ar[dd]_{Hf}\ar@{} [ddr] |{=} & GU \ar[dd]^{Gh}\ar[rr]^{N_U} & & HU
\ar@/_1.5pc/[dd]_{Hh} \ar@/^1.5pc/[dd]^{Hf} &\\
 & &  \ar@{}[ll]|{\Longleftarrow}^{G\gamma} &
 \ar@{}[l]|{\Longleftarrow \;\;\;\;\;\;\;\;\;\;}^<>(.8){Nf} & &
  \ar@{}[l]|{\;\;\;\;\;\;\;\;\;\;\Longleftarrow}^<>(.2){Nh}
 & &
 \ar@{}[ll]|{\Longleftarrow}^{H\gamma} \\
& GV \ar[rr]_{N_V} & & HV & GV \ar[rr]_{N_V}& & HV &
}
\]
for morphisms $f,h\colon U \rightarrow V$ and 2-morphism $\gamma\colon
f \rightarrow h$. 

A \emph{strict monoidal} 2-category is a 2-category $\cB$ with a
strict 2-functor $\otimes:\cB \times \cB \rightarrow \cB$ which is
strictly associative and has a strict left and right identity element
$\bf1$. The monoidal structure is \emph{semi-strict} if $\otimes$ is allowed to have non-trivial 2-isomorphisms
$\otimes_{f_i,g_i}:(f_1f_2)\otimes (g_1g_2)\rightarrow (f_1\otimes
g_1)(f_2\otimes g_2)$ (with $\otimes$ being strict elsewhere; $\otimes$ is no longer a 2-functor). 

A \emph{monoidal} 2-functor $G:\cD \rightarrow \cE$ between semi-strict monoidal strict 2-categories comes equipped with isomorphisms $M^G_{U,V}: G(U\otimes V)\rightarrow GU\otimes GV$, for objects $U,V\in \cD$; these must satisfy
\[
\xymatrix{G(U\otimes V) \ar[dd]_{G(f\otimes g)} \ar[rr]^{M^G_{U,V}} && GU\otimes GV\ar[dd]^{Gf\otimes Gg}\\ \\
G(U'\otimes V') \ar[rr]^{M^G_{U',V'}} &&GU'\otimes GV'
}
\]
for morphisms $f:U\rightarrow U'$ and $g:V\rightarrow V'$ in $\cD$. We also require that via the isomorphisms $M^G$, $G(\gamma \otimes \delta)=G(\gamma)\otimes G(\delta)$ for 2-morphisms $\gamma$, $\delta$ (to be more precise the equality ressembles that of 2-morphisms given in the definition of 2-natural transformation above, but with the squares now the identity 2-morphisms of the previous diagram)

A pseudo 2-natural transformation $N$ between monoidal 2-functors $G$ and $H$ is said to be \emph{monoidal} if the following diagram commutes
\[
\xymatrix{G(U\otimes V) \ar[rr]^{N_{U\otimes V}} \ar[dd]_{M^G_{U,V}} && H(U\otimes V) \ar[dd]^{M^H_{U,V}}\\ \\
GU\otimes GV \ar[rr]^{N_U\otimes N_V} && HU\otimes HV
}
\] 
We also require that for morphisms $f$ and $g$, $N_{f\otimes g}=N_f\otimes N_g$ under the isomorphisms $M^G$ and $M^F$ (again we omit the obvious pasting diagram). Notice the definitions just given clearly admit laxer and stricter versions; for further details on 2-categories we refer to \cite{KellyStreet}\cite{Maclane}.

The 2-categories central to this paper are ones, roughly speaking,
whose objects are  collections of loops in a space $X$, whose
morphisms are 2-manifolds (with boundary) in $X$, and whose
2-morphisms are diffeomorphisms of 2-manifolds in $X$.
In practice great care is needed with the definition and all models
have their own advantages and deficiencies. As in topological quantum field theory it is
extremely useful to be able to make geometric arguments using surfaces
and we will extensively use surface diagram manipulation. For this
we need a fairly explicit hold on surfaces and we generalise the
2-category of Tillmann \cite{Tillmann:Discrete} in which all surfaces are embedded in
$\bR^3$.

Let $S_m$ denote $m$ circles of
radius $1/4$ centred at $(1,0),(2,0),\dots ,(m,0)$ in $\bR^2$ and let
$S_0$ denote the empty set. By {\em surface} we shall mean a smooth
cobordism in $\bR^3$ with boundary circles lying on the planes $z=0$
and $z=t$ for some $t>0$.  The circles are oriented counterclockwise when viewed from $z>>0$. We require that in some neighbourhood of
each boundary component, the surfaces be straight cylinders of radius
$1/4$ and further for simplicity that projection onto the
$z$-coordinate is a Morse function. Notice that such surfaces are
canonically oriented by choosing inward pointing normals. The boundary
components in the plane $z=0$ are called {\em inputs} and those on the
plane $z=t$ are called {\em outputs}. Two surfaces $\Sigma_1$ and
$\Sigma_2$ can be glued together by shifting $\Sigma_2$ vertically by
$t_1$ (the height of $\Sigma_1$) and gluing along the boundary
circles; the result is again smooth since the collars are straight
cylinders.

\begin{defn} \label{defn:SDX}
Let $X$ be a based topological space. Define a 2-category $\cS_{D,X}$ as follows.
\begin{itemize}
\item Objects: based continuous functions $s\colon S_m \rightarrow X$,
  for $m\in \bN$. 
\item 1-Morphisms:  continuous functions $g\colon \Sigma \rightarrow
  X$ where $\Sigma$ is a surface (as above). The source and target are
  $g|_{z=0}$ and $g|_{z=t}$ respectively. On the straight boundary
  collars $g$ must factor through the projection to the boundary.
\item 2-Morphisms: orientation preserving diffeomorphisms $T\colon
  \Sigma_1 \rightarrow \Sigma_2$ that fix boundary collars pointwise
  and such that the following diagram commutes up to
  basepoint-preserving homotopy relative to the boundary.
\[
\xymatrix{\Sigma_1 \ar[rr]^{T} \ar[dr]_{g_1} & & \Sigma_2
  \ar[dl]^{g_2}\\ & X}
\]
\end{itemize}
Composition of 1-morphisms is defined by gluing surfaces and taking
the induced map to $X$. Similarly vertical composition $\circ_1$ of
2-morphisms is composition of diffeomorphisms and horizontal
composition $\circ_2$ is union of diffeomorphisms induced by the
gluing of surfaces. Where there is no ambiguity we will write $s$
for the object $s\colon S_m \rightarrow X$ and $g$ for the
morphism $g\colon \Sigma \rightarrow X$.
\end{defn}

Strictly speaking this is not a 2-category as there are no identity
morphisms. This will be remedied in what follows where we introduce
limited isotopy of surfaces needed  in order to have a well
defined monoidal product given by disjoint union. Again we appeal to
\cite{Tillmann:Discrete} and use Tillmann's RS-moves which we now recall. Let
$\Sigma$ be a surface. By choosing
$0=t_0<t_1<\dots t_{k-1}<t_k=t$ cut up $\Sigma$ into slices $\Sigma_i
\subset \bR^2\times [t_{i-1},t_i]$ for $i=1,\dots ,k$, with components
$\Sigma_i^j$ for $j=1,\dots j_i$. Define a {\em rescaling} to be a
collection  of functions $R_{i,j}(x,y,z)=(x,y,r_{i,j}(z))$
where $r_{i,j}:[t_{i-1},t_i]\rightarrow [t_{i-1},\tilde{t_i}]$ is a
smooth function with derivative 1 in a neighbourhood of the boundary
of $[t_{i-1},t_i]$. Here $\tilde{t_i}$ is independent of $j$, and the
surfaces $\Sigma_{i+1},\dots,\Sigma_k$ are shifted vertically by the
appropriate amount. Define a {\em shift} to be a collection of functions
$
S_{i,j}(x,y,z)=(x,y,z)+(f_{i,j}^1(z),f_{i,j}^2(z),0)
$
where $f_{i,j}^l:[t_{i-1},t_i]\rightarrow \bR$ are smooth functions
vanishing on some neighbourhood of the boundary of $[t_{i-1},t_i]$. 
A 2-morphism $T$ in $\cS_{D,X}$ which is defined by a finite number of
rescalings and shifts is called an {\em RS 2-morphism}. 

We now
define the \emph{homotopy surface 2-category} $\cS_X$ by taking
successive quotients of $\cS_{D,X}$ on the 1 and 2-morphisms. 
Define an equivalence relation $\equiv_2$ on the set of 2-morphisms
$\cS_{D,X}(g_1,g_2)$ by setting $T_1 \equiv_2 T_2$ if and only if they
are in the same connected component of
$\odiff(\Sigma_1,\Sigma_2;\partial)$. Next define an equivalence
relation $\equiv_1$ on the 1-morphisms by setting $g\equiv_1 g^\prime$
if and only if  $g$ and $g^\prime$ are
isomorphic via an RS 2-morphism.

\begin{defn}\label{sx}
The {\em homotopy surface 2-category of} $X$ is defined to be the
2-category obtained by quotienting 2-morphisms in $\cS_{D,X}$ by
$\equiv_2$ followed by quotienting 1-morphisms by
$\equiv_1$. 
\end{defn}

The identity 1-morphism in $\cS_X (s, s)$ can be taken to be a
collection of straight cylinders mapping to $X$ via $s \circ p$ where
$p$ is projection onto $\bR^2\times 0$.  Standard category theory
implies that on taking the quotients above the result is a 2-category.
Also, notice that by taking $X$ to be a one point space the resultant
category is Tillmann's. We can also form a category from the above
2-category by identifying all 2-isomorphic 1-morphisms and this is a
model for the homotopy surface category of a space defined in
\cite{BrightwellTurner:Representations}.

We can now define a monoidal structure $\cS_X$ using disjoint union.
Define a 2-functor  $\sqcup
:\cS_X\times \cS_X\rightarrow \cS_X$ on objects $s\colon S_m \rightarrow X$ and
  $s^\prime\colon S_n \rightarrow X$ to be the obvious map $s\sqcup
s^\prime\colon S_{m+n} \rightarrow X$. On 1-morphisms $g\colon \Sigma \rightarrow X$ and
  $g^\prime\colon \Sigma^\prime \rightarrow X$ let $\Sigma \cup
  \Sigma^\prime$ be the surface obtained by rescaling the height of
  $\Sigma^\prime$ to that of $\Sigma$ and shifting the result by a
  diffeomorphism $(x,y,z) \mapsto (x+f(z), y, z)$ where $f\colon
  [0,t] \rightarrow \bR $ is constantly $k$ in a neighbourhood of
  $0$ and constantly $l$ in a neighbourhood of $t$ and such that
  result is disjoint from $\Sigma$. Then $g\sqcup g^\prime$ is
  defined by the induced maps. Finally, on 2-morphisms we take the 
disjoint union of diffeomorphisms. 
The analysis in \cite{Tillmann:Discrete} proves that 
$\cS_X$ is a semi-strict monoidal strict 2-category.

The arguments in this paper rest crucially on the fact that surfaces
in $\bR^3$ can be manipulated by a number of geometric operations
defined as follows.  Given $s\colon S_1 \rightarrow X$ let $s^{-1}$ be
the map $s$ precomposed with reflection in the line $x=1$.  For
$s=s_1\sqcup \cdots \sqcup s_n \colon S_n \rightarrow X$ let $s^{-1} =
s_n^{-1}\sqcup \cdots \sqcup s_1^{-1}$.  \gap

\begin{bfseries}Reflection.\end{bfseries}
There is a contravariant monoidal equivalence of
2-categories  
\[
\xyreflect {\makebox[0.6cm]{}}\colon  \cS_X \rightarrow \cS_X
\]
defined on objects by $s \mapsto s$ and on a morphism $g\colon \Sigma
\rightarrow X$ by reflecting $\Sigma$ in the plane
$z=0$ and then translating the result. The map to $X$ is the induced
one and we denote the result by $\xyreflect g$. Note that reflecting
in the plane $z=0$ leaves the object circles and their maps to $X$
unchanged. A 2-morphism will induce a 2-morphism on reflected
1-morphisms.

\vspace{0.3cm}
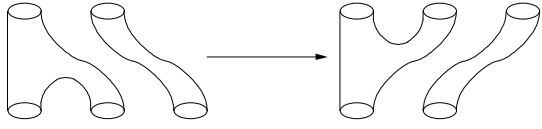
\begin{center}
\setlength{\unitlength}{0.00029167in}
\begingroup\makeatletter\ifx\SetFigFont\undefined%
\gdef\SetFigFont#1#2#3#4#5{%
  \reset@font\fontsize{#1}{#2pt}%
  \fontfamily{#3}\fontseries{#4}\fontshape{#5}%
  \selectfont}%
\fi\endgroup%
{\renewcommand{\dashlinestretch}{30}
\begin{picture}(9624,2129)(0,-10)
\put(7812,1957){\ellipse{600}{300}}
\put(6312,157){\ellipse{600}{300}}
\put(6312,1957){\ellipse{600}{300}}
\path(6612,1957)(6612,1807)
\path(7512,1957)(7512,1807)
\path(6612,157)(6612,158)(6612,161)
        (6612,166)(6613,174)(6613,185)
        (6614,199)(6616,216)(6618,235)
        (6621,257)(6625,281)(6631,307)
        (6638,334)(6646,363)(6657,393)
        (6669,425)(6685,458)(6704,493)
        (6726,531)(6753,571)(6785,614)
        (6822,659)(6864,707)(6912,757)
        (6950,794)(6990,830)(7029,864)
        (7067,896)(7103,924)(7137,949)
        (7169,971)(7198,989)(7224,1005)
        (7249,1017)(7271,1028)(7291,1036)
        (7310,1043)(7328,1048)(7345,1053)
        (7362,1057)(7379,1061)(7396,1066)
        (7414,1071)(7433,1078)(7453,1086)
        (7475,1097)(7500,1109)(7526,1125)
        (7555,1143)(7587,1165)(7621,1190)
        (7657,1218)(7695,1250)(7734,1284)
        (7774,1320)(7812,1357)(7860,1407)
        (7902,1455)(7939,1500)(7971,1543)
        (7998,1583)(8020,1621)(8039,1656)
        (8055,1689)(8067,1721)(8078,1751)
        (8086,1780)(8093,1807)(8099,1833)
        (8103,1857)(8106,1879)(8108,1898)
        (8110,1915)(8111,1929)(8111,1940)
        (8112,1948)(8112,1953)(8112,1956)(8112,1957)
\path(6012,1957)(6012,157)
\path(6612,1807)(6612,1805)(6613,1801)
        (6615,1794)(6617,1784)(6620,1771)
        (6624,1755)(6629,1737)(6635,1718)
        (6642,1697)(6650,1676)(6660,1654)
        (6671,1631)(6684,1608)(6699,1584)
        (6717,1558)(6738,1533)(6762,1507)
        (6785,1485)(6807,1466)(6827,1450)
        (6845,1436)(6860,1425)(6871,1416)
        (6881,1409)(6888,1403)(6894,1398)
        (6900,1394)(6905,1391)(6912,1387)
        (6920,1384)(6930,1380)(6943,1375)
        (6960,1371)(6980,1366)(7005,1362)
        (7032,1358)(7062,1357)(7092,1358)
        (7119,1362)(7144,1366)(7164,1371)
        (7181,1375)(7194,1380)(7204,1384)
        (7212,1387)(7219,1391)(7225,1395)
        (7230,1398)(7236,1403)(7243,1409)
        (7253,1416)(7264,1425)(7279,1436)
        (7297,1450)(7317,1466)(7339,1485)
        (7362,1507)(7386,1533)(7407,1558)
        (7425,1584)(7440,1608)(7453,1631)
        (7464,1654)(7474,1676)(7482,1697)
        (7489,1718)(7495,1737)(7500,1755)
        (7504,1771)(7507,1784)(7509,1794)
        (7511,1801)(7512,1805)(7512,1807)
\put(9312,1957){\ellipse{600}{300}}
\put(7812,157){\ellipse{600}{300}}
\path(8112,157)(8112,158)(8112,161)
        (8112,166)(8113,174)(8113,185)
        (8114,199)(8116,216)(8118,235)
        (8121,257)(8125,281)(8131,307)
        (8138,334)(8146,363)(8157,393)
        (8169,425)(8185,458)(8204,493)
        (8226,531)(8253,571)(8285,614)
        (8322,659)(8364,707)(8412,757)
        (8450,794)(8490,830)(8529,864)
        (8567,896)(8603,924)(8637,949)
        (8669,971)(8698,989)(8724,1005)
        (8749,1017)(8771,1028)(8791,1036)
        (8810,1043)(8828,1048)(8845,1053)
        (8862,1057)(8879,1061)(8896,1066)
        (8914,1071)(8933,1078)(8953,1086)
        (8975,1097)(9000,1109)(9026,1125)
        (9055,1143)(9087,1165)(9121,1190)
        (9157,1218)(9195,1250)(9234,1284)
        (9274,1320)(9312,1357)(9360,1407)
        (9402,1455)(9439,1500)(9471,1543)
        (9498,1583)(9520,1621)(9539,1656)
        (9555,1689)(9567,1721)(9578,1751)
        (9586,1780)(9593,1807)(9599,1833)
        (9603,1857)(9606,1879)(9608,1898)
        (9610,1915)(9611,1929)(9611,1940)
        (9612,1948)(9612,1953)(9612,1956)(9612,1957)
\path(7512,157)(7512,158)(7512,161)
        (7512,166)(7513,174)(7513,185)
        (7514,199)(7516,216)(7518,235)
        (7521,257)(7525,281)(7531,307)
        (7538,334)(7546,363)(7557,393)
        (7569,425)(7585,458)(7604,493)
        (7626,531)(7653,571)(7685,614)
        (7722,659)(7764,707)(7812,757)
        (7850,794)(7890,830)(7929,864)
        (7967,896)(8003,924)(8037,949)
        (8069,971)(8098,989)(8124,1005)
        (8149,1017)(8171,1028)(8191,1036)
        (8210,1043)(8228,1048)(8245,1053)
        (8262,1057)(8279,1061)(8296,1066)
        (8314,1071)(8333,1078)(8353,1086)
        (8375,1097)(8400,1109)(8426,1125)
        (8455,1143)(8487,1165)(8521,1190)
        (8557,1218)(8595,1250)(8634,1284)
        (8674,1320)(8712,1357)(8760,1407)
        (8802,1455)(8839,1500)(8871,1543)
        (8898,1583)(8920,1621)(8939,1656)
        (8955,1689)(8967,1721)(8978,1751)
        (8986,1780)(8993,1807)(8999,1833)
        (9003,1857)(9006,1879)(9008,1898)
        (9010,1915)(9011,1929)(9011,1940)
        (9012,1948)(9012,1953)(9012,1956)(9012,1957)
\put(1812,157){\ellipse{600}{300}}
\put(312,1957){\ellipse{600}{300}}
\put(312,157){\ellipse{600}{300}}
\put(3312,157){\ellipse{600}{300}}
\put(1812,1957){\ellipse{600}{300}}
\path(612,157)(612,307)
\path(1512,157)(1512,307)
\path(3612,1132)(5712,1132)
\blacken\thicklines
\path(5592.000,1102.000)(5712.000,1132.000)(5592.000,1162.000)(5592.000,1102.000)
\thinlines
\path(612,1957)(612,1956)(612,1953)
        (612,1948)(613,1940)(613,1929)
        (614,1915)(616,1898)(618,1879)
        (621,1857)(625,1833)(631,1807)
        (638,1780)(646,1751)(657,1721)
        (669,1689)(685,1656)(704,1621)
        (726,1583)(753,1543)(785,1500)
        (822,1455)(864,1407)(912,1357)
        (950,1320)(990,1284)(1029,1250)
        (1067,1218)(1103,1190)(1137,1165)
        (1169,1143)(1198,1125)(1224,1109)
        (1249,1097)(1271,1086)(1291,1078)
        (1310,1071)(1328,1066)(1345,1061)
        (1362,1057)(1379,1053)(1396,1048)
        (1414,1043)(1433,1036)(1453,1028)
        (1475,1017)(1500,1005)(1526,989)
        (1555,971)(1587,949)(1621,924)
        (1657,896)(1695,864)(1734,830)
        (1774,794)(1812,757)(1860,707)
        (1902,659)(1939,614)(1971,571)
        (1998,531)(2020,493)(2039,458)
        (2055,425)(2067,393)(2078,363)
        (2086,334)(2093,307)(2099,281)
        (2103,257)(2106,235)(2108,216)
        (2110,199)(2111,185)(2111,174)
        (2112,166)(2112,161)(2112,158)(2112,157)
\path(12,157)(12,1957)
\path(612,307)(612,309)(613,313)
        (615,320)(617,330)(620,343)
        (624,359)(629,377)(635,396)
        (642,417)(650,438)(660,460)
        (671,483)(684,506)(699,530)
        (717,556)(738,581)(762,607)
        (785,629)(807,648)(827,664)
        (845,678)(860,689)(871,698)
        (881,705)(888,711)(894,716)
        (900,720)(905,723)(912,727)
        (920,730)(930,734)(943,739)
        (960,743)(980,748)(1005,752)
        (1032,756)(1062,757)(1092,756)
        (1119,752)(1144,748)(1164,743)
        (1181,739)(1194,734)(1204,730)
        (1212,727)(1219,723)(1225,719)
        (1230,716)(1236,711)(1243,705)
        (1253,698)(1264,689)(1279,678)
        (1297,664)(1317,648)(1339,629)
        (1362,607)(1386,581)(1407,556)
        (1425,530)(1440,506)(1453,483)
        (1464,460)(1474,438)(1482,417)
        (1489,396)(1495,377)(1500,359)
        (1504,343)(1507,330)(1509,320)
        (1511,313)(1512,309)(1512,307)
\path(2112,1957)(2112,1956)(2112,1953)
        (2112,1948)(2113,1940)(2113,1929)
        (2114,1915)(2116,1898)(2118,1879)
        (2121,1857)(2125,1833)(2131,1807)
        (2138,1780)(2146,1751)(2157,1721)
        (2169,1689)(2185,1656)(2204,1621)
        (2226,1583)(2253,1543)(2285,1500)
        (2322,1455)(2364,1407)(2412,1357)
        (2450,1320)(2490,1284)(2529,1250)
        (2567,1218)(2603,1190)(2637,1165)
        (2669,1143)(2698,1125)(2724,1109)
        (2749,1097)(2771,1086)(2791,1078)
        (2810,1071)(2828,1066)(2845,1061)
        (2862,1057)(2879,1053)(2896,1048)
        (2914,1043)(2933,1036)(2953,1028)
        (2975,1017)(3000,1005)(3026,989)
        (3055,971)(3087,949)(3121,924)
        (3157,896)(3195,864)(3234,830)
        (3274,794)(3312,757)(3360,707)
        (3402,659)(3439,614)(3471,571)
        (3498,531)(3520,493)(3539,458)
        (3555,425)(3567,393)(3578,363)
        (3586,334)(3593,307)(3599,281)
        (3603,257)(3606,235)(3608,216)
        (3610,199)(3611,185)(3611,174)
        (3612,166)(3612,161)(3612,158)(3612,157)
\path(1512,1957)(1512,1956)(1512,1953)
        (1512,1948)(1513,1940)(1513,1929)
        (1514,1915)(1516,1898)(1518,1879)
        (1521,1857)(1525,1833)(1531,1807)
        (1538,1780)(1546,1751)(1557,1721)
        (1569,1689)(1585,1656)(1604,1621)
        (1626,1583)(1653,1543)(1685,1500)
        (1722,1455)(1764,1407)(1812,1357)
        (1850,1320)(1890,1284)(1929,1250)
        (1967,1218)(2003,1190)(2037,1165)
        (2069,1143)(2098,1125)(2124,1109)
        (2149,1097)(2171,1086)(2191,1078)
        (2210,1071)(2228,1066)(2245,1061)
        (2262,1057)(2279,1053)(2296,1048)
        (2314,1043)(2333,1036)(2353,1028)
        (2375,1017)(2400,1005)(2426,989)
        (2455,971)(2487,949)(2521,924)
        (2557,896)(2595,864)(2634,830)
        (2674,794)(2712,757)(2760,707)
        (2802,659)(2839,614)(2871,571)
        (2898,531)(2920,493)(2939,458)
        (2955,425)(2967,393)(2978,363)
        (2986,334)(2993,307)(2999,281)
        (3003,257)(3006,235)(3008,216)
        (3010,199)(3011,185)(3011,174)
        (3012,166)(3012,161)(3012,158)(3012,157)
\end{picture}
}
\end{center}
\vspace{0.3cm}

\begin{bfseries}Rotation.\end{bfseries}
$\cS_X$ has another identification with its opposite by rotating
surfaces by 180 degrees. Define a contravariant 
2-functor 
\[
\rotate {\makebox[0.6cm]{}}\colon  \cS_X \rightarrow \cS_X
\]
on objects by $s \mapsto s^{-1}$ and on a morphism $g\colon \Sigma
\rightarrow X$ by rotating $\Sigma$ by 180 degrees
around the $y$-axis and then adjusting by using a diffeomorphism
$(x,y,z) \mapsto (x+f(z),y,z)$ where $f\colon [0,t_1]\rightarrow \bR$
is as in the definition of monoidal product.  The map to $X$ is the
induced one. Denote the result by $\rotate g$ and again take
induced 2-morphisms. Note that this is ``anti-monoidal'' rather than
monoidal. 

\vspace{0.3cm}
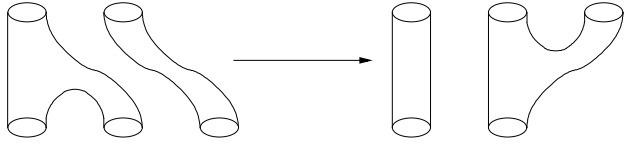
\begin{center}
\setlength{\unitlength}{0.00033333in}
\begingroup\makeatletter\ifx\SetFigFont\undefined%
\gdef\SetFigFont#1#2#3#4#5{%
  \reset@font\fontsize{#1}{#2pt}%
  \fontfamily{#3}\fontseries{#4}\fontshape{#5}%
  \selectfont}%
\fi\endgroup%
{\renewcommand{\dashlinestretch}{30}
\begin{picture}(9624,2129)(0,-10)
\put(9312,1957){\ellipse{600}{300}}
\put(7812,157){\ellipse{600}{300}}
\put(7812,1957){\ellipse{600}{300}}
\path(8112,1957)(8112,1807)
\path(9012,1957)(9012,1807)
\path(8112,157)(8112,158)(8112,161)
        (8112,166)(8113,174)(8113,185)
        (8114,199)(8116,216)(8118,235)
        (8121,257)(8125,281)(8131,307)
        (8138,334)(8146,363)(8157,393)
        (8169,425)(8185,458)(8204,493)
        (8226,531)(8253,571)(8285,614)
        (8322,659)(8364,707)(8412,757)
        (8450,794)(8490,830)(8529,864)
        (8567,896)(8603,924)(8637,949)
        (8669,971)(8698,989)(8724,1005)
        (8749,1017)(8771,1028)(8791,1036)
        (8810,1043)(8828,1048)(8845,1053)
        (8862,1057)(8879,1061)(8896,1066)
        (8914,1071)(8933,1078)(8953,1086)
        (8975,1097)(9000,1109)(9026,1125)
        (9055,1143)(9087,1165)(9121,1190)
        (9157,1218)(9195,1250)(9234,1284)
        (9274,1320)(9312,1357)(9360,1407)
        (9402,1455)(9439,1500)(9471,1543)
        (9498,1583)(9520,1621)(9539,1656)
        (9555,1689)(9567,1721)(9578,1751)
        (9586,1780)(9593,1807)(9599,1833)
        (9603,1857)(9606,1879)(9608,1898)
        (9610,1915)(9611,1929)(9611,1940)
        (9612,1948)(9612,1953)(9612,1956)(9612,1957)
\path(7512,1957)(7512,157)
\path(8112,1807)(8112,1805)(8113,1801)
        (8115,1794)(8117,1784)(8120,1771)
        (8124,1755)(8129,1737)(8135,1718)
        (8142,1697)(8150,1676)(8160,1654)
        (8171,1631)(8184,1608)(8199,1584)
        (8217,1558)(8238,1533)(8262,1507)
        (8285,1485)(8307,1466)(8327,1450)
        (8345,1436)(8360,1425)(8371,1416)
        (8381,1409)(8388,1403)(8394,1398)
        (8400,1394)(8405,1391)(8412,1387)
        (8420,1384)(8430,1380)(8443,1375)
        (8460,1371)(8480,1366)(8505,1362)
        (8532,1358)(8562,1357)(8592,1358)
        (8619,1362)(8644,1366)(8664,1371)
        (8681,1375)(8694,1380)(8704,1384)
        (8712,1387)(8719,1391)(8725,1395)
        (8730,1398)(8736,1403)(8743,1409)
        (8753,1416)(8764,1425)(8779,1436)
        (8797,1450)(8817,1466)(8839,1485)
        (8862,1507)(8886,1533)(8907,1558)
        (8925,1584)(8940,1608)(8953,1631)
        (8964,1654)(8974,1676)(8982,1697)
        (8989,1718)(8995,1737)(9000,1755)
        (9004,1771)(9007,1784)(9009,1794)
        (9011,1801)(9012,1805)(9012,1807)
\put(1812,157){\ellipse{600}{300}}
\put(312,1957){\ellipse{600}{300}}
\put(312,157){\ellipse{600}{300}}
\put(3312,157){\ellipse{600}{300}}
\put(1812,1957){\ellipse{600}{300}}
\put(6312,1957){\ellipse{600}{300}}
\put(6312,157){\ellipse{600}{300}}
\path(612,157)(612,307)
\path(1512,157)(1512,307)
\path(6012,157)(6012,1957)
\path(6612,1957)(6612,157)
\path(3537,1207)(5637,1207)
\blacken\thicklines
\path(5517.000,1177.000)(5637.000,1207.000)(5517.000,1237.000)(5517.000,1177.000)
\thinlines
\path(612,1957)(612,1956)(612,1953)
        (612,1948)(613,1940)(613,1929)
        (614,1915)(616,1898)(618,1879)
        (621,1857)(625,1833)(631,1807)
        (638,1780)(646,1751)(657,1721)
        (669,1689)(685,1656)(704,1621)
        (726,1583)(753,1543)(785,1500)
        (822,1455)(864,1407)(912,1357)
        (950,1320)(990,1284)(1029,1250)
        (1067,1218)(1103,1190)(1137,1165)
        (1169,1143)(1198,1125)(1224,1109)
        (1249,1097)(1271,1086)(1291,1078)
        (1310,1071)(1328,1066)(1345,1061)
        (1362,1057)(1379,1053)(1396,1048)
        (1414,1043)(1433,1036)(1453,1028)
        (1475,1017)(1500,1005)(1526,989)
        (1555,971)(1587,949)(1621,924)
        (1657,896)(1695,864)(1734,830)
        (1774,794)(1812,757)(1860,707)
        (1902,659)(1939,614)(1971,571)
        (1998,531)(2020,493)(2039,458)
        (2055,425)(2067,393)(2078,363)
        (2086,334)(2093,307)(2099,281)
        (2103,257)(2106,235)(2108,216)
        (2110,199)(2111,185)(2111,174)
        (2112,166)(2112,161)(2112,158)(2112,157)
\path(12,157)(12,1957)
\path(612,307)(612,309)(613,313)
        (615,320)(617,330)(620,343)
        (624,359)(629,377)(635,396)
        (642,417)(650,438)(660,460)
        (671,483)(684,506)(699,530)
        (717,556)(738,581)(762,607)
        (785,629)(807,648)(827,664)
        (845,678)(860,689)(871,698)
        (881,705)(888,711)(894,716)
        (900,720)(905,723)(912,727)
        (920,730)(930,734)(943,739)
        (960,743)(980,748)(1005,752)
        (1032,756)(1062,757)(1092,756)
        (1119,752)(1144,748)(1164,743)
        (1181,739)(1194,734)(1204,730)
        (1212,727)(1219,723)(1225,719)
        (1230,716)(1236,711)(1243,705)
        (1253,698)(1264,689)(1279,678)
        (1297,664)(1317,648)(1339,629)
        (1362,607)(1386,581)(1407,556)
        (1425,530)(1440,506)(1453,483)
        (1464,460)(1474,438)(1482,417)
        (1489,396)(1495,377)(1500,359)
        (1504,343)(1507,330)(1509,320)
        (1511,313)(1512,309)(1512,307)
\path(2112,1957)(2112,1956)(2112,1953)
        (2112,1948)(2113,1940)(2113,1929)
        (2114,1915)(2116,1898)(2118,1879)
        (2121,1857)(2125,1833)(2131,1807)
        (2138,1780)(2146,1751)(2157,1721)
        (2169,1689)(2185,1656)(2204,1621)
        (2226,1583)(2253,1543)(2285,1500)
        (2322,1455)(2364,1407)(2412,1357)
        (2450,1320)(2490,1284)(2529,1250)
        (2567,1218)(2603,1190)(2637,1165)
        (2669,1143)(2698,1125)(2724,1109)
        (2749,1097)(2771,1086)(2791,1078)
        (2810,1071)(2828,1066)(2845,1061)
        (2862,1057)(2879,1053)(2896,1048)
        (2914,1043)(2933,1036)(2953,1028)
        (2975,1017)(3000,1005)(3026,989)
        (3055,971)(3087,949)(3121,924)
        (3157,896)(3195,864)(3234,830)
        (3274,794)(3312,757)(3360,707)
        (3402,659)(3439,614)(3471,571)
        (3498,531)(3520,493)(3539,458)
        (3555,425)(3567,393)(3578,363)
        (3586,334)(3593,307)(3599,281)
        (3603,257)(3606,235)(3608,216)
        (3610,199)(3611,185)(3611,174)
        (3612,166)(3612,161)(3612,158)(3612,157)
\path(1512,1957)(1512,1956)(1512,1953)
        (1512,1948)(1513,1940)(1513,1929)
        (1514,1915)(1516,1898)(1518,1879)
        (1521,1857)(1525,1833)(1531,1807)
        (1538,1780)(1546,1751)(1557,1721)
        (1569,1689)(1585,1656)(1604,1621)
        (1626,1583)(1653,1543)(1685,1500)
        (1722,1455)(1764,1407)(1812,1357)
        (1850,1320)(1890,1284)(1929,1250)
        (1967,1218)(2003,1190)(2037,1165)
        (2069,1143)(2098,1125)(2124,1109)
        (2149,1097)(2171,1086)(2191,1078)
        (2210,1071)(2228,1066)(2245,1061)
        (2262,1057)(2279,1053)(2296,1048)
        (2314,1043)(2333,1036)(2353,1028)
        (2375,1017)(2400,1005)(2426,989)
        (2455,971)(2487,949)(2521,924)
        (2557,896)(2595,864)(2634,830)
        (2674,794)(2712,757)(2760,707)
        (2802,659)(2839,614)(2871,571)
        (2898,531)(2920,493)(2939,458)
        (2955,425)(2967,393)(2978,363)
        (2986,334)(2993,307)(2999,281)
        (3003,257)(3006,235)(3008,216)
        (3010,199)(3011,185)(3011,174)
        (3012,166)(3012,161)(3012,158)(3012,157)
\end{picture}
}

\end{center}
\vspace{0.3cm}

The composite
\[
\xymatrix{\cS_X \ar[r]^{\rotate {\makebox[0.6cm]{}}} & \cS_X
  \ar[r]^{\xyreflect {\makebox[0.6cm]{}}} & \cS_X}
\]
is 
given on objects by 
$s \mapsto s^{-1}$
and on a morphism $g\colon \Sigma
\rightarrow X$ by reflecting in the plane
$x=0$ and then adjusting similarly to above.

\section{$\cS_X$-structures}\label{sec:structures}

Now for the main definition of the paper. Let $k$ be an algebraically
closed field and let $\additive$ be the
2-category of idempotent complete $k$-additive categories with monoidal structure given by
the tensor product (see Appendix \ref{app:add} for details).

\begin{defn}\label{defn:hqft}
  Let $X$ be a based space. An \emph{$\cS_X$-structure} is a monoidal 2-functor of strict 2-categories
  $F\colon \cS_X \rightarrow \additive $ such that the crossed
  cylinders are mapped to the functor changing components.
\end{defn}

Thus to each collection of loops $s\colon S_m \rightarrow X$ we assign
an additive category; to each surface $g\colon \Sigma
\rightarrow X$ a functor of additive categories and to each
diffeomorphism $T$ of surfaces we assign a natural transformation of
functors. This assignment is monoidal taking disjoint union to tensor
product.

We make a few remarks on how the above definition relates to other notions in homotopy quantum field theory. Let us first consider the background free case, or topological quantum field theory. As was the point of view in \cite{Tillmann:S-Structures}, $\cS$-structures have the flavour of a modular functor in dimension 2, since they essentially consist of some kind of representations of the mapping class groups. To be more precise, the standard notion of a modular functor appears in an $\cS$-structure from those surfaces with an empty target 1-manifold (see \cite{Tillmann:S-Structures} for details). On the other hand $\cS$-structures fit nicely into the picture of extended topological quantum field theories in dimension 3. This becomes clear if one considers the definition given in \cite{KerlerLyubashenko}, where an extended TQFT is defined as a representation of the double category of circles, surfaces between these, and relative 3-cobordisms (a 2-category is a special case of a double category). From this point of view $\cS$ should be viewed as a sub-2-category of the latter (we avoid being precise here). One could repeat the above remarks for the more general $\cS_X$-structures, comparing them to homotopy modular functors and extended homotopy quantum field theories in dimension 3. 

By recalling that there is a contravariant functor
$(-)^\vee\colon \additive \rightarrow \additive$ taking an additive
category $\cA$ to its dual $\cA^\vee = \additive (\cA , \hat{k})$,
the dual of an $\cS_X$-structure is defined as follows.

\begin{defn}
  The {\em dual} of $F\colon \cS_X \rightarrow \additive$, denoted
  $F^\vee$ is defined as the composite
\[
\cS_X \stackrel{\xyreflect{\makebox[0.6cm]{}}}{\rightarrow} \cS_X
\stackrel{F}{\rightarrow} \additive \stackrel{(-)^\vee}{\rightarrow}
\additive
\]
\end{defn}

Later we will restrict ourselves to theories which are self-dual in a
way we now make precise.  Given a finitely generated $k$-additive
category $\cA$ i.e. one whose morphism vector spaces are finitely
generated, one can define a functor $hom\colon \cA \rightarrow
\cA^\vee$ given by $Y\mapsto \cA (Y, -)$. Suppose $F$ is an $\cS_X$-structure taking values among $k$-additive categories that are
finitely generated (in fact this is always the case) then there is a
family of maps
\[
\{hom_s \colon F(s) \rightarrow F^\vee (s)
\}_{s\in Ob(\cS_X)}.
\]

\begin{defn}\label{defn:selfdual}
  An $\cS_X$-structure $F$ is {\em lax self dual with respect to hom} if
  the above family provide a monoidal pseudo 2-natural transformation $N$ between
  $F$ and $F^\vee$, and cylinders $g\colon I \rightarrow X$ satisfy
  $N(g)=F(g)^{-1}$.
\end{defn}

Two important consequences of this definition are the following.

\begin{list}{S-\Roman{enumi}}{\usecounter{enumi}}
\item \label{S1} If $g\colon \Sigma \rightarrow X $ is  a morphism then
  $F(\xyreflect g)$ is right adjoint to $F(g)$
  i.e. there are natural isomorphisms
\[
\cB(F(g)(U), V) \cong \cA(U, F(\xyreflect g)(V))
\]
where $F(g)\colon \cA \rightarrow \cB$. For the collapsed cylinders this isomorphism is the identity.
\item \label{S2} If 
  $T\colon g_1 \rightarrow g_2$ is a 2-morphism
  then the following diagram commutes
\[
\xymatrix{\cB(F({g_1})U,V) \ar[r]^\simeq &
  \cA(U,F(\xyreflect {g_1})V)
  \ar[d]^{F(\xyreflect{T})_V}\\ 
  \cB(F({g_2})U,V) \ar[u]^{F({T})_U} \ar[r]^\simeq &
  \cA(U,F(\xyreflect {g_2})V)}
\]
where 
  $\xyreflect{T}:\xyreflect{g_1}\rightarrow
  \xyreflect{g_2}$ is  obtained from $T$ under
  reflection and the vertical arrows are given by pre and post-composition.
\end{list}

\section{Balanced categories from $\cS_X$-structures}\label{section:balanced}
Let $X$ be a based space. An  $\cS_X$-structure $F\colon \cS_X \rightarrow
\additive $ determines a collection of categories $\{
\cA_{\alpha}\}_{\alpha\in \pi}$ indexed by the group $\pi = \pi_1X$ as
follows.  A loop $\alpha$ in $X$ determines a map $s_\alpha\colon
S_1\rightarrow X$, taking $(1,-1/4)$ as basepoint of $S_1$. For each
element of $\pi=\pi_1(X)$ choose a representative loop $\alpha$ and
set

\[
\cA_{\alpha} = F({s_\alpha}) \;\;\;\;\;\;\; \cA =
\bigsqcup_{\alpha\in \pi} \cA_\alpha 
\]

The following proposition refers to the definitions of balanced
categories in Appendix \ref{app:defns}.

\begin{thm}\label{prop:balanced}
(a) Let $\pi$ be a discrete group and let $X$ be a based
  Eilenberg-Maclane space $K(\pi,1)$. The $k$-additive category $\cA$
  associated to an $\cS_X$-structure is a balanced $\pi$-category.  \\
(b) For any space $X$, the subcategory $\cA_1$ is a balanced
  category with $\pi_2X$-action.\\
(c) The categories above are semi-simple Artinian categories.
\end{thm}

The proof of this proposition will take up the rest of this
section. 

Firstly consider part (a) in which $X=K(\pi,1)$. To define a monoidal structure
on $\cA$ let $\alpha, \beta \in \pi$ and pick a pair of pants surface
$P$ with two inputs and one output and let $p_{\alpha,\beta}\colon P
\rightarrow X$ be a map inducing the maps indicated below.

\vspace{0.3cm}
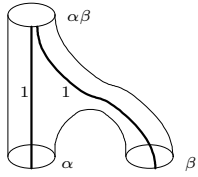
\begin{center}
\setlength{\unitlength}{0.00041667in}
\begingroup\makeatletter\ifx\SetFigFont\undefined%
\gdef\SetFigFont#1#2#3#4#5{%
  \reset@font\fontsize{#1}{#2pt}%
  \fontfamily{#3}\fontseries{#4}\fontshape{#5}%
  \selectfont}%
\fi\endgroup%
{\renewcommand{\dashlinestretch}{30}
\begin{picture}(2554,2176)(0,-10)
\put(1812,204){\ellipse{600}{300}}
\put(312,2004){\ellipse{600}{300}}
\put(312,204){\ellipse{600}{300}}
\path(612,204)(612,354)
\path(1512,204)(1512,354)
\thicklines
\path(312,54)(312,1854)
\thinlines
\path(612,2004)(612,2003)(612,2000)
        (612,1995)(613,1987)(613,1976)
        (614,1962)(616,1945)(618,1926)
        (621,1904)(625,1880)(631,1854)
        (638,1827)(646,1798)(657,1768)
        (669,1736)(685,1703)(704,1668)
        (726,1630)(753,1590)(785,1547)
        (822,1502)(864,1454)(912,1404)
        (950,1367)(990,1331)(1029,1297)
        (1067,1265)(1103,1237)(1137,1212)
        (1169,1190)(1198,1172)(1224,1156)
        (1249,1144)(1271,1133)(1291,1125)
        (1310,1118)(1328,1113)(1345,1108)
        (1362,1104)(1379,1100)(1396,1095)
        (1414,1090)(1433,1083)(1453,1075)
        (1475,1064)(1500,1052)(1526,1036)
        (1555,1018)(1587,996)(1621,971)
        (1657,943)(1695,911)(1734,877)
        (1774,841)(1812,804)(1860,754)
        (1902,706)(1939,661)(1971,618)
        (1998,578)(2020,540)(2039,505)
        (2055,472)(2067,440)(2078,410)
        (2086,381)(2093,354)(2099,328)
        (2103,304)(2106,282)(2108,263)
        (2110,246)(2111,232)(2111,221)
        (2112,213)(2112,208)(2112,205)(2112,204)
\path(12,204)(12,2004)
\path(612,354)(612,356)(613,360)
        (615,367)(617,377)(620,390)
        (624,406)(629,424)(635,443)
        (642,464)(650,485)(660,507)
        (671,530)(684,553)(699,577)
        (717,603)(738,628)(762,654)
        (785,676)(807,695)(827,711)
        (845,725)(860,736)(871,745)
        (881,752)(888,758)(894,763)
        (900,767)(905,770)(912,774)
        (920,777)(930,781)(943,786)
        (960,790)(980,795)(1005,799)
        (1032,803)(1062,804)(1092,803)
        (1119,799)(1144,795)(1164,790)
        (1181,786)(1194,781)(1204,777)
        (1212,774)(1219,770)(1225,766)
        (1230,763)(1236,758)(1243,752)
        (1253,745)(1264,736)(1279,725)
        (1297,711)(1317,695)(1339,676)
        (1362,654)(1386,628)(1407,603)
        (1425,577)(1440,553)(1453,530)
        (1464,507)(1474,485)(1482,464)
        (1489,443)(1495,424)(1500,406)
        (1504,390)(1507,377)(1509,367)
        (1511,360)(1512,356)(1512,354)
\thicklines
\path(387,1854)(387,1853)(387,1850)
        (387,1845)(388,1837)(388,1826)
        (389,1812)(391,1795)(393,1776)
        (396,1754)(400,1730)(406,1704)
        (413,1677)(421,1648)(432,1618)
        (444,1586)(460,1553)(479,1518)
        (501,1480)(528,1440)(560,1397)
        (597,1352)(639,1304)(687,1254)
        (725,1217)(765,1181)(804,1147)
        (842,1115)(878,1087)(912,1062)
        (944,1040)(973,1022)(999,1006)
        (1024,994)(1046,983)(1066,975)
        (1085,968)(1103,963)(1120,958)
        (1137,954)(1154,950)(1171,945)
        (1189,940)(1208,933)(1228,925)
        (1250,914)(1275,902)(1301,886)
        (1330,868)(1362,846)(1396,821)
        (1432,793)(1470,761)(1509,727)
        (1549,691)(1587,654)(1635,604)
        (1677,556)(1714,511)(1746,468)
        (1773,428)(1795,390)(1814,355)
        (1830,322)(1842,290)(1853,260)
        (1861,231)(1868,204)(1874,178)
        (1878,154)(1881,132)(1883,113)
        (1885,96)(1886,82)(1886,71)
        (1887,63)(1887,58)(1887,55)(1887,54)
\put(162,954){\makebox(0,0)[lb]{\smash{{{\SetFigFont{6}{7.2}{\rmdefault}{\mddefault}{\updefault}$1$}}}}}
\put(687,954){\makebox(0,0)[lb]{\smash{{{\SetFigFont{6}{7.2}{\rmdefault}{\mddefault}{\updefault}$1$}}}}}
\put(762,1929){\makebox(0,0)[lb]{\smash{{{\SetFigFont{6}{7.2}{\rmdefault}{\mddefault}{\updefault}$\alpha\beta$}}}}}
\put(2262,54){\makebox(0,0)[lb]{\smash{{{\SetFigFont{6}{7.2}{\rmdefault}{\mddefault}{\updefault}$\beta$}}}}}
\put(687,54){\makebox(0,0)[lb]{\smash{{{\SetFigFont{6}{7.2}{\rmdefault}{\mddefault}{\updefault}$\alpha$}}}}}
\end{picture}
}

\end{center}
\vspace{0.3cm}

Here the label indicates the map on the given line or boundary
component and any two choices of $p_{\alpha,\beta}$ give the same 1-morphism in $\cS_X$.
Define functors
\begin{eqnarray*}
*_{\alpha,\beta}&=&F({p_{\alpha,\beta}}):\cA_\alpha \otimes
\cA_\beta \rightarrow \cA_{\alpha\beta}\\   
*&=&\bigsqcup_{\alpha,\beta\in \pi}*_{\alpha,\beta}:\cA \otimes \cA
\rightarrow \cA  
\end{eqnarray*}

Now choose a disc $D$ with one input only and let
$d\colon D\rightarrow X$ be the collapse map to a basepoint of $X$. 
Noting that $F(d)\colon \hat k \rightarrow \cA_1$ define
\[
\monunit = F(D^d)(k) \in \cA_1
\]

It can now be seen that 
$(\cA,*,\monunit)$ is a $\pi$-graded monoidal category by 
choosing a diffeomorphism

\vspace{0.3cm}
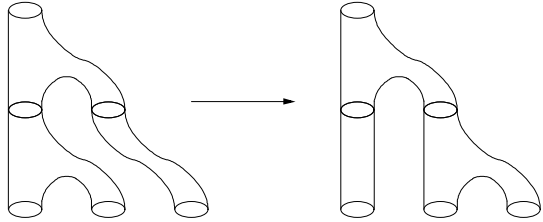
\begin{center}
\setlength{\unitlength}{0.00029167in}
\begingroup\makeatletter\ifx\SetFigFont\undefined%
\gdef\SetFigFont#1#2#3#4#5{%
  \reset@font\fontsize{#1}{#2pt}%
  \fontfamily{#3}\fontseries{#4}\fontshape{#5}%
  \selectfont}%
\fi\endgroup%
{\renewcommand{\dashlinestretch}{30}
\begin{picture}(9624,3929)(0,-10)
\put(9312,157){\ellipse{600}{300}}
\put(7812,1957){\ellipse{600}{300}}
\put(7812,157){\ellipse{600}{300}}
\path(8112,157)(8112,307)
\path(9012,157)(9012,307)
\path(8112,1957)(8112,1956)(8112,1953)
        (8112,1948)(8113,1940)(8113,1929)
        (8114,1915)(8116,1898)(8118,1879)
        (8121,1857)(8125,1833)(8131,1807)
        (8138,1780)(8146,1751)(8157,1721)
        (8169,1689)(8185,1656)(8204,1621)
        (8226,1583)(8253,1543)(8285,1500)
        (8322,1455)(8364,1407)(8412,1357)
        (8450,1320)(8490,1284)(8529,1250)
        (8567,1218)(8603,1190)(8637,1165)
        (8669,1143)(8698,1125)(8724,1109)
        (8749,1097)(8771,1086)(8791,1078)
        (8810,1071)(8828,1066)(8845,1061)
        (8862,1057)(8879,1053)(8896,1048)
        (8914,1043)(8933,1036)(8953,1028)
        (8975,1017)(9000,1005)(9026,989)
        (9055,971)(9087,949)(9121,924)
        (9157,896)(9195,864)(9234,830)
        (9274,794)(9312,757)(9360,707)
        (9402,659)(9439,614)(9471,571)
        (9498,531)(9520,493)(9539,458)
        (9555,425)(9567,393)(9578,363)
        (9586,334)(9593,307)(9599,281)
        (9603,257)(9606,235)(9608,216)
        (9610,199)(9611,185)(9611,174)
        (9612,166)(9612,161)(9612,158)(9612,157)
\path(7512,157)(7512,1957)
\path(8112,307)(8112,309)(8113,313)
        (8115,320)(8117,330)(8120,343)
        (8124,359)(8129,377)(8135,396)
        (8142,417)(8150,438)(8160,460)
        (8171,483)(8184,506)(8199,530)
        (8217,556)(8238,581)(8262,607)
        (8285,629)(8307,648)(8327,664)
        (8345,678)(8360,689)(8371,698)
        (8381,705)(8388,711)(8394,716)
        (8400,720)(8405,723)(8412,727)
        (8420,730)(8430,734)(8443,739)
        (8460,743)(8480,748)(8505,752)
        (8532,756)(8562,757)(8592,756)
        (8619,752)(8644,748)(8664,743)
        (8681,739)(8694,734)(8704,730)
        (8712,727)(8719,723)(8725,719)
        (8730,716)(8736,711)(8743,705)
        (8753,698)(8764,689)(8779,678)
        (8797,664)(8817,648)(8839,629)
        (8862,607)(8886,581)(8907,556)
        (8925,530)(8940,506)(8953,483)
        (8964,460)(8974,438)(8982,417)
        (8989,396)(8995,377)(9000,359)
        (9004,343)(9007,330)(9009,320)
        (9011,313)(9012,309)(9012,307)
\put(1812,1957){\ellipse{600}{300}}
\put(312,3757){\ellipse{600}{300}}
\put(312,1957){\ellipse{600}{300}}
\path(612,1957)(612,2107)
\path(1512,1957)(1512,2107)
\path(612,3757)(612,3756)(612,3753)
        (612,3748)(613,3740)(613,3729)
        (614,3715)(616,3698)(618,3679)
        (621,3657)(625,3633)(631,3607)
        (638,3580)(646,3551)(657,3521)
        (669,3489)(685,3456)(704,3421)
        (726,3383)(753,3343)(785,3300)
        (822,3255)(864,3207)(912,3157)
        (950,3120)(990,3084)(1029,3050)
        (1067,3018)(1103,2990)(1137,2965)
        (1169,2943)(1198,2925)(1224,2909)
        (1249,2897)(1271,2886)(1291,2878)
        (1310,2871)(1328,2866)(1345,2861)
        (1362,2857)(1379,2853)(1396,2848)
        (1414,2843)(1433,2836)(1453,2828)
        (1475,2817)(1500,2805)(1526,2789)
        (1555,2771)(1587,2749)(1621,2724)
        (1657,2696)(1695,2664)(1734,2630)
        (1774,2594)(1812,2557)(1860,2507)
        (1902,2459)(1939,2414)(1971,2371)
        (1998,2331)(2020,2293)(2039,2258)
        (2055,2225)(2067,2193)(2078,2163)
        (2086,2134)(2093,2107)(2099,2081)
        (2103,2057)(2106,2035)(2108,2016)
        (2110,1999)(2111,1985)(2111,1974)
        (2112,1966)(2112,1961)(2112,1958)(2112,1957)
\path(12,1957)(12,3757)
\path(612,2107)(612,2109)(613,2113)
        (615,2120)(617,2130)(620,2143)
        (624,2159)(629,2177)(635,2196)
        (642,2217)(650,2238)(660,2260)
        (671,2283)(684,2306)(699,2330)
        (717,2356)(738,2381)(762,2407)
        (785,2429)(807,2448)(827,2464)
        (845,2478)(860,2489)(871,2498)
        (881,2505)(888,2511)(894,2516)
        (900,2520)(905,2523)(912,2527)
        (920,2530)(930,2534)(943,2539)
        (960,2543)(980,2548)(1005,2552)
        (1032,2556)(1062,2557)(1092,2556)
        (1119,2552)(1144,2548)(1164,2543)
        (1181,2539)(1194,2534)(1204,2530)
        (1212,2527)(1219,2523)(1225,2519)
        (1230,2516)(1236,2511)(1243,2505)
        (1253,2498)(1264,2489)(1279,2478)
        (1297,2464)(1317,2448)(1339,2429)
        (1362,2407)(1386,2381)(1407,2356)
        (1425,2330)(1440,2306)(1453,2283)
        (1464,2260)(1474,2238)(1482,2217)
        (1489,2196)(1495,2177)(1500,2159)
        (1504,2143)(1507,2130)(1509,2120)
        (1511,2113)(1512,2109)(1512,2107)
\put(7812,1957){\ellipse{600}{300}}
\put(6312,3757){\ellipse{600}{300}}
\put(6312,1957){\ellipse{600}{300}}
\path(6612,1957)(6612,2107)
\path(7512,1957)(7512,2107)
\path(6612,3757)(6612,3756)(6612,3753)
        (6612,3748)(6613,3740)(6613,3729)
        (6614,3715)(6616,3698)(6618,3679)
        (6621,3657)(6625,3633)(6631,3607)
        (6638,3580)(6646,3551)(6657,3521)
        (6669,3489)(6685,3456)(6704,3421)
        (6726,3383)(6753,3343)(6785,3300)
        (6822,3255)(6864,3207)(6912,3157)
        (6950,3120)(6990,3084)(7029,3050)
        (7067,3018)(7103,2990)(7137,2965)
        (7169,2943)(7198,2925)(7224,2909)
        (7249,2897)(7271,2886)(7291,2878)
        (7310,2871)(7328,2866)(7345,2861)
        (7362,2857)(7379,2853)(7396,2848)
        (7414,2843)(7433,2836)(7453,2828)
        (7475,2817)(7500,2805)(7526,2789)
        (7555,2771)(7587,2749)(7621,2724)
        (7657,2696)(7695,2664)(7734,2630)
        (7774,2594)(7812,2557)(7860,2507)
        (7902,2459)(7939,2414)(7971,2371)
        (7998,2331)(8020,2293)(8039,2258)
        (8055,2225)(8067,2193)(8078,2163)
        (8086,2134)(8093,2107)(8099,2081)
        (8103,2057)(8106,2035)(8108,2016)
        (8110,1999)(8111,1985)(8111,1974)
        (8112,1966)(8112,1961)(8112,1958)(8112,1957)
\path(6012,1957)(6012,3757)
\path(6612,2107)(6612,2109)(6613,2113)
        (6615,2120)(6617,2130)(6620,2143)
        (6624,2159)(6629,2177)(6635,2196)
        (6642,2217)(6650,2238)(6660,2260)
        (6671,2283)(6684,2306)(6699,2330)
        (6717,2356)(6738,2381)(6762,2407)
        (6785,2429)(6807,2448)(6827,2464)
        (6845,2478)(6860,2489)(6871,2498)
        (6881,2505)(6888,2511)(6894,2516)
        (6900,2520)(6905,2523)(6912,2527)
        (6920,2530)(6930,2534)(6943,2539)
        (6960,2543)(6980,2548)(7005,2552)
        (7032,2556)(7062,2557)(7092,2556)
        (7119,2552)(7144,2548)(7164,2543)
        (7181,2539)(7194,2534)(7204,2530)
        (7212,2527)(7219,2523)(7225,2519)
        (7230,2516)(7236,2511)(7243,2505)
        (7253,2498)(7264,2489)(7279,2478)
        (7297,2464)(7317,2448)(7339,2429)
        (7362,2407)(7386,2381)(7407,2356)
        (7425,2330)(7440,2306)(7453,2283)
        (7464,2260)(7474,2238)(7482,2217)
        (7489,2196)(7495,2177)(7500,2159)
        (7504,2143)(7507,2130)(7509,2120)
        (7511,2113)(7512,2109)(7512,2107)
\put(1812,157){\ellipse{600}{300}}
\put(312,1957){\ellipse{600}{300}}
\put(312,157){\ellipse{600}{300}}
\put(3312,157){\ellipse{600}{300}}
\put(1812,1957){\ellipse{600}{300}}
\put(6312,1957){\ellipse{600}{300}}
\put(6312,157){\ellipse{600}{300}}
\path(612,157)(612,307)
\path(1512,157)(1512,307)
\path(6012,157)(6012,1957)
\path(6612,1957)(6612,157)
\path(3312,2107)(5112,2107)
\blacken\thicklines
\path(4992.000,2077.000)(5112.000,2107.000)(4992.000,2137.000)(4992.000,2077.000)
\thinlines
\path(612,1957)(612,1956)(612,1953)
        (612,1948)(613,1940)(613,1929)
        (614,1915)(616,1898)(618,1879)
        (621,1857)(625,1833)(631,1807)
        (638,1780)(646,1751)(657,1721)
        (669,1689)(685,1656)(704,1621)
        (726,1583)(753,1543)(785,1500)
        (822,1455)(864,1407)(912,1357)
        (950,1320)(990,1284)(1029,1250)
        (1067,1218)(1103,1190)(1137,1165)
        (1169,1143)(1198,1125)(1224,1109)
        (1249,1097)(1271,1086)(1291,1078)
        (1310,1071)(1328,1066)(1345,1061)
        (1362,1057)(1379,1053)(1396,1048)
        (1414,1043)(1433,1036)(1453,1028)
        (1475,1017)(1500,1005)(1526,989)
        (1555,971)(1587,949)(1621,924)
        (1657,896)(1695,864)(1734,830)
        (1774,794)(1812,757)(1860,707)
        (1902,659)(1939,614)(1971,571)
        (1998,531)(2020,493)(2039,458)
        (2055,425)(2067,393)(2078,363)
        (2086,334)(2093,307)(2099,281)
        (2103,257)(2106,235)(2108,216)
        (2110,199)(2111,185)(2111,174)
        (2112,166)(2112,161)(2112,158)(2112,157)
\path(12,157)(12,1957)
\path(612,307)(612,309)(613,313)
        (615,320)(617,330)(620,343)
        (624,359)(629,377)(635,396)
        (642,417)(650,438)(660,460)
        (671,483)(684,506)(699,530)
        (717,556)(738,581)(762,607)
        (785,629)(807,648)(827,664)
        (845,678)(860,689)(871,698)
        (881,705)(888,711)(894,716)
        (900,720)(905,723)(912,727)
        (920,730)(930,734)(943,739)
        (960,743)(980,748)(1005,752)
        (1032,756)(1062,757)(1092,756)
        (1119,752)(1144,748)(1164,743)
        (1181,739)(1194,734)(1204,730)
        (1212,727)(1219,723)(1225,719)
        (1230,716)(1236,711)(1243,705)
        (1253,698)(1264,689)(1279,678)
        (1297,664)(1317,648)(1339,629)
        (1362,607)(1386,581)(1407,556)
        (1425,530)(1440,506)(1453,483)
        (1464,460)(1474,438)(1482,417)
        (1489,396)(1495,377)(1500,359)
        (1504,343)(1507,330)(1509,320)
        (1511,313)(1512,309)(1512,307)
\path(2112,1957)(2112,1956)(2112,1953)
        (2112,1948)(2113,1940)(2113,1929)
        (2114,1915)(2116,1898)(2118,1879)
        (2121,1857)(2125,1833)(2131,1807)
        (2138,1780)(2146,1751)(2157,1721)
        (2169,1689)(2185,1656)(2204,1621)
        (2226,1583)(2253,1543)(2285,1500)
        (2322,1455)(2364,1407)(2412,1357)
        (2450,1320)(2490,1284)(2529,1250)
        (2567,1218)(2603,1190)(2637,1165)
        (2669,1143)(2698,1125)(2724,1109)
        (2749,1097)(2771,1086)(2791,1078)
        (2810,1071)(2828,1066)(2845,1061)
        (2862,1057)(2879,1053)(2896,1048)
        (2914,1043)(2933,1036)(2953,1028)
        (2975,1017)(3000,1005)(3026,989)
        (3055,971)(3087,949)(3121,924)
        (3157,896)(3195,864)(3234,830)
        (3274,794)(3312,757)(3360,707)
        (3402,659)(3439,614)(3471,571)
        (3498,531)(3520,493)(3539,458)
        (3555,425)(3567,393)(3578,363)
        (3586,334)(3593,307)(3599,281)
        (3603,257)(3606,235)(3608,216)
        (3610,199)(3611,185)(3611,174)
        (3612,166)(3612,161)(3612,158)(3612,157)
\path(1512,1957)(1512,1956)(1512,1953)
        (1512,1948)(1513,1940)(1513,1929)
        (1514,1915)(1516,1898)(1518,1879)
        (1521,1857)(1525,1833)(1531,1807)
        (1538,1780)(1546,1751)(1557,1721)
        (1569,1689)(1585,1656)(1604,1621)
        (1626,1583)(1653,1543)(1685,1500)
        (1722,1455)(1764,1407)(1812,1357)
        (1850,1320)(1890,1284)(1929,1250)
        (1967,1218)(2003,1190)(2037,1165)
        (2069,1143)(2098,1125)(2124,1109)
        (2149,1097)(2171,1086)(2191,1078)
        (2210,1071)(2228,1066)(2245,1061)
        (2262,1057)(2279,1053)(2296,1048)
        (2314,1043)(2333,1036)(2353,1028)
        (2375,1017)(2400,1005)(2426,989)
        (2455,971)(2487,949)(2521,924)
        (2557,896)(2595,864)(2634,830)
        (2674,794)(2712,757)(2760,707)
        (2802,659)(2839,614)(2871,571)
        (2898,531)(2920,493)(2939,458)
        (2955,425)(2967,393)(2978,363)
        (2986,334)(2993,307)(2999,281)
        (3003,257)(3006,235)(3008,216)
        (3010,199)(3011,185)(3011,174)
        (3012,166)(3012,161)(3012,158)(3012,157)
\end{picture}
}

\end{center}
\vspace{0.3cm}
which gives a natural isomorphism $a\colon (-*-)*- \nattrans -*(-*-)$
i.e. a collection of isomorphisms $a_{U,V,W}$ as required. The
isomorphism $r_U$ is obtained from a diffeomorphism

\vspace{0.3cm}
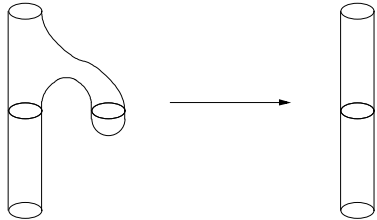
\begin{center}
\setlength{\unitlength}{0.00029167in}
\begingroup\makeatletter\ifx\SetFigFont\undefined%
\gdef\SetFigFont#1#2#3#4#5{%
  \reset@font\fontsize{#1}{#2pt}%
  \fontfamily{#3}\fontseries{#4}\fontshape{#5}%
  \selectfont}%
\fi\endgroup%
{\renewcommand{\dashlinestretch}{30}
\begin{picture}(6624,3929)(0,-10)
\put(1812.000,1807.000){\arc{600.000}{6.2832}{9.4248}}
\put(1812,1957){\ellipse{600}{300}}
\path(1512,1807)(1512,1957)
\path(2112,1957)(2112,1807)
\put(6312,3757){\ellipse{600}{300}}
\put(6312,1957){\ellipse{600}{300}}
\path(6012,1957)(6012,3757)
\path(6612,3757)(6612,1957)
\put(312,1957){\ellipse{600}{300}}
\put(312,157){\ellipse{600}{300}}
\path(12,157)(12,1957)
\path(612,1957)(612,157)
\put(6312,1957){\ellipse{600}{300}}
\put(6312,157){\ellipse{600}{300}}
\path(6012,157)(6012,1957)
\path(6612,1957)(6612,157)
\put(1812,1957){\ellipse{600}{300}}
\put(312,3757){\ellipse{600}{300}}
\put(312,1957){\ellipse{600}{300}}
\path(612,1957)(612,2107)
\path(1512,1957)(1512,2107)
\path(2937,2107)(5037,2107)
\blacken\thicklines
\path(4917.000,2077.000)(5037.000,2107.000)(4917.000,2137.000)(4917.000,2077.000)
\thinlines
\path(612,3757)(612,3756)(612,3753)
        (612,3748)(613,3740)(613,3729)
        (614,3715)(616,3698)(618,3679)
        (621,3657)(625,3633)(631,3607)
        (638,3580)(646,3551)(657,3521)
        (669,3489)(685,3456)(704,3421)
        (726,3383)(753,3343)(785,3300)
        (822,3255)(864,3207)(912,3157)
        (950,3120)(990,3084)(1029,3050)
        (1067,3018)(1103,2990)(1137,2965)
        (1169,2943)(1198,2925)(1224,2909)
        (1249,2897)(1271,2886)(1291,2878)
        (1310,2871)(1328,2866)(1345,2861)
        (1362,2857)(1379,2853)(1396,2848)
        (1414,2843)(1433,2836)(1453,2828)
        (1475,2817)(1500,2805)(1526,2789)
        (1555,2771)(1587,2749)(1621,2724)
        (1657,2696)(1695,2664)(1734,2630)
        (1774,2594)(1812,2557)(1860,2507)
        (1902,2459)(1939,2414)(1971,2371)
        (1998,2331)(2020,2293)(2039,2258)
        (2055,2225)(2067,2193)(2078,2163)
        (2086,2134)(2093,2107)(2099,2081)
        (2103,2057)(2106,2035)(2108,2016)
        (2110,1999)(2111,1985)(2111,1974)
        (2112,1966)(2112,1961)(2112,1958)(2112,1957)
\path(12,1957)(12,3757)
\path(612,2107)(612,2109)(613,2113)
        (615,2120)(617,2130)(620,2143)
        (624,2159)(629,2177)(635,2196)
        (642,2217)(650,2238)(660,2260)
        (671,2283)(684,2306)(699,2330)
        (717,2356)(738,2381)(762,2407)
        (785,2429)(807,2448)(827,2464)
        (845,2478)(860,2489)(871,2498)
        (881,2505)(888,2511)(894,2516)
        (900,2520)(905,2523)(912,2527)
        (920,2530)(930,2534)(943,2539)
        (960,2543)(980,2548)(1005,2552)
        (1032,2556)(1062,2557)(1092,2556)
        (1119,2552)(1144,2548)(1164,2543)
        (1181,2539)(1194,2534)(1204,2530)
        (1212,2527)(1219,2523)(1225,2519)
        (1230,2516)(1236,2511)(1243,2505)
        (1253,2498)(1264,2489)(1279,2478)
        (1297,2464)(1317,2448)(1339,2429)
        (1362,2407)(1386,2381)(1407,2356)
        (1425,2330)(1440,2306)(1453,2283)
        (1464,2260)(1474,2238)(1482,2217)
        (1489,2196)(1495,2177)(1500,2159)
        (1504,2143)(1507,2130)(1509,2120)
        (1511,2113)(1512,2109)(1512,2107)
\end{picture}
}

\end{center}
\vspace{0.3cm}
\noindent
and similarly for $l_U$.
 Finally, commutativity of the
 associativity pentagons and identity triangles holds since by
 definition the compositions differ by at most an isotopy.

Next we show how to obtain structure leading to a balanced
$\pi$-category.

\begin{bfseries}Crossing.\end{bfseries}
Let $\alpha, \gamma \in \pi$ and pick a cylinder $I$ and a map
$i_{\alpha, \gamma}\colon I \rightarrow X$ as shown below

\vspace{0.3cm}
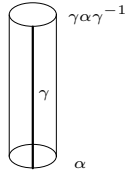
\begin{center}
\setlength{\unitlength}{0.00041667in}
\begingroup\makeatletter\ifx\SetFigFont\undefined%
\gdef\SetFigFont#1#2#3#4#5{%
  \reset@font\fontsize{#1}{#2pt}%
  \fontfamily{#3}\fontseries{#4}\fontshape{#5}%
  \selectfont}%
\fi\endgroup%
{\renewcommand{\dashlinestretch}{30}
\begin{picture}(2192,2176)(0,-10)
\put(312,2004){\ellipse{600}{300}}
\put(312,204){\ellipse{600}{300}}
\path(12,204)(12,2004)
\path(612,2004)(612,204)
\thicklines
\path(312,1854)(312,54)
\put(762,1929){\makebox(0,0)[lb]{\smash{{{\SetFigFont{6}{7.2}{\rmdefault}{\mddefault}{\updefault}$\gamma\alpha\gamma^{-1}$}}}}}
\put(387,954){\makebox(0,0)[lb]{\smash{{{\SetFigFont{6}{7.2}{\rmdefault}{\mddefault}{\updefault}$\gamma$}}}}}
\put(837,54){\makebox(0,0)[lb]{\smash{{{\SetFigFont{6}{7.2}{\rmdefault}{\mddefault}{\updefault}$\alpha$}}}}}
\end{picture}
}

\end{center}
\vspace{0.3cm}

Define 
$\varphi(\gamma)_\alpha :=F(i_{\alpha,\gamma})\colon \cA_\alpha \rightarrow
\cA_{\gamma\alpha\gamma^{-1}} $ and assemble these into a map $\varphi \colon \pi \rightarrow \oAut (\cA)$.
Note that $\varphi(\gamma)$ is invertible with inverse
$\varphi(\gamma^{-1})$ and also that 
$\varphi (\gamma)$ respects $*$  as required, by
the equality below.

\vspace{0.3cm}
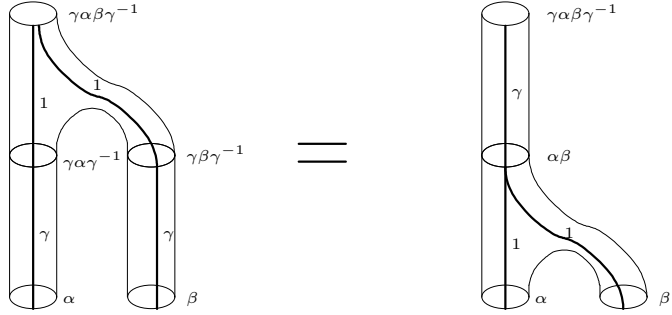
\begin{center}
\setlength{\unitlength}{0.00041667in}
\begingroup\makeatletter\ifx\SetFigFont\undefined%
\gdef\SetFigFont#1#2#3#4#5{%
  \reset@font\fontsize{#1}{#2pt}%
  \fontfamily{#3}\fontseries{#4}\fontshape{#5}%
  \selectfont}%
\fi\endgroup%
{\renewcommand{\dashlinestretch}{30}
\begin{picture}(8354,3944)(0,-10)
\put(6312,3772){\ellipse{600}{300}}
\put(6312,1972){\ellipse{600}{300}}
\path(6012,1972)(6012,3772)
\path(6612,3772)(6612,1972)
\put(312,1972){\ellipse{600}{300}}
\put(312,172){\ellipse{600}{300}}
\path(12,172)(12,1972)
\path(612,1972)(612,172)
\put(1812,1972){\ellipse{600}{300}}
\put(1812,172){\ellipse{600}{300}}
\path(1512,172)(1512,1972)
\path(2112,1972)(2112,172)
\put(7812,172){\ellipse{600}{300}}
\put(6312,1972){\ellipse{600}{300}}
\put(6312,172){\ellipse{600}{300}}
\path(6612,172)(6612,322)
\path(7512,172)(7512,322)
\path(6612,1972)(6612,1971)(6612,1968)
        (6612,1963)(6613,1955)(6613,1944)
        (6614,1930)(6616,1913)(6618,1894)
        (6621,1872)(6625,1848)(6631,1822)
        (6638,1795)(6646,1766)(6657,1736)
        (6669,1704)(6685,1671)(6704,1636)
        (6726,1598)(6753,1558)(6785,1515)
        (6822,1470)(6864,1422)(6912,1372)
        (6950,1335)(6990,1299)(7029,1265)
        (7067,1233)(7103,1205)(7137,1180)
        (7169,1158)(7198,1140)(7224,1124)
        (7249,1112)(7271,1101)(7291,1093)
        (7310,1086)(7328,1081)(7345,1076)
        (7362,1072)(7379,1068)(7396,1063)
        (7414,1058)(7433,1051)(7453,1043)
        (7475,1032)(7500,1020)(7526,1004)
        (7555,986)(7587,964)(7621,939)
        (7657,911)(7695,879)(7734,845)
        (7774,809)(7812,772)(7860,722)
        (7902,674)(7939,629)(7971,586)
        (7998,546)(8020,508)(8039,473)
        (8055,440)(8067,408)(8078,378)
        (8086,349)(8093,322)(8099,296)
        (8103,272)(8106,250)(8108,231)
        (8110,214)(8111,200)(8111,189)
        (8112,181)(8112,176)(8112,173)(8112,172)
\path(6012,172)(6012,1972)
\path(6612,322)(6612,324)(6613,328)
        (6615,335)(6617,345)(6620,358)
        (6624,374)(6629,392)(6635,411)
        (6642,432)(6650,453)(6660,475)
        (6671,498)(6684,521)(6699,545)
        (6717,571)(6738,596)(6762,622)
        (6785,644)(6807,663)(6827,679)
        (6845,693)(6860,704)(6871,713)
        (6881,720)(6888,726)(6894,731)
        (6900,735)(6905,738)(6912,742)
        (6920,745)(6930,749)(6943,754)
        (6960,758)(6980,763)(7005,767)
        (7032,771)(7062,772)(7092,771)
        (7119,767)(7144,763)(7164,758)
        (7181,754)(7194,749)(7204,745)
        (7212,742)(7219,738)(7225,734)
        (7230,731)(7236,726)(7243,720)
        (7253,713)(7264,704)(7279,693)
        (7297,679)(7317,663)(7339,644)
        (7362,622)(7386,596)(7407,571)
        (7425,545)(7440,521)(7453,498)
        (7464,475)(7474,453)(7482,432)
        (7489,411)(7495,392)(7500,374)
        (7504,358)(7507,345)(7509,335)
        (7511,328)(7512,324)(7512,322)
\thicklines
\path(3687,2122)(4287,2122)
\path(3687,1897)(4287,1897)
\thinlines
\put(1812,1972){\ellipse{600}{300}}
\put(312,3772){\ellipse{600}{300}}
\put(312,1972){\ellipse{600}{300}}
\path(612,1972)(612,2122)
\path(1512,1972)(1512,2122)
\thicklines
\path(312,1822)(312,22)
\path(1887,1822)(1887,22)
\path(312,3622)(312,1822)
\path(6312,1822)(6312,22)
\path(6312,3622)(6312,1822)
\thinlines
\path(612,3772)(612,3771)(612,3768)
        (612,3763)(613,3755)(613,3744)
        (614,3730)(616,3713)(618,3694)
        (621,3672)(625,3648)(631,3622)
        (638,3595)(646,3566)(657,3536)
        (669,3504)(685,3471)(704,3436)
        (726,3398)(753,3358)(785,3315)
        (822,3270)(864,3222)(912,3172)
        (950,3135)(990,3099)(1029,3065)
        (1067,3033)(1103,3005)(1137,2980)
        (1169,2958)(1198,2940)(1224,2924)
        (1249,2912)(1271,2901)(1291,2893)
        (1310,2886)(1328,2881)(1345,2876)
        (1362,2872)(1379,2868)(1396,2863)
        (1414,2858)(1433,2851)(1453,2843)
        (1475,2832)(1500,2820)(1526,2804)
        (1555,2786)(1587,2764)(1621,2739)
        (1657,2711)(1695,2679)(1734,2645)
        (1774,2609)(1812,2572)(1860,2522)
        (1902,2474)(1939,2429)(1971,2386)
        (1998,2346)(2020,2308)(2039,2273)
        (2055,2240)(2067,2208)(2078,2178)
        (2086,2149)(2093,2122)(2099,2096)
        (2103,2072)(2106,2050)(2108,2031)
        (2110,2014)(2111,2000)(2111,1989)
        (2112,1981)(2112,1976)(2112,1973)(2112,1972)
\path(12,1972)(12,3772)
\path(612,2122)(612,2124)(613,2128)
        (615,2135)(617,2145)(620,2158)
        (624,2174)(629,2192)(635,2211)
        (642,2232)(650,2253)(660,2275)
        (671,2298)(684,2321)(699,2345)
        (717,2371)(738,2396)(762,2422)
        (785,2444)(807,2463)(827,2479)
        (845,2493)(860,2504)(871,2513)
        (881,2520)(888,2526)(894,2531)
        (900,2535)(905,2538)(912,2542)
        (920,2545)(930,2549)(943,2554)
        (960,2558)(980,2563)(1005,2567)
        (1032,2571)(1062,2572)(1092,2571)
        (1119,2567)(1144,2563)(1164,2558)
        (1181,2554)(1194,2549)(1204,2545)
        (1212,2542)(1219,2538)(1225,2534)
        (1230,2531)(1236,2526)(1243,2520)
        (1253,2513)(1264,2504)(1279,2493)
        (1297,2479)(1317,2463)(1339,2444)
        (1362,2422)(1386,2396)(1407,2371)
        (1425,2345)(1440,2321)(1453,2298)
        (1464,2275)(1474,2253)(1482,2232)
        (1489,2211)(1495,2192)(1500,2174)
        (1504,2158)(1507,2145)(1509,2135)
        (1511,2128)(1512,2124)(1512,2122)
\thicklines
\path(387,3622)(387,3621)(387,3618)
        (387,3613)(388,3605)(388,3594)
        (389,3580)(391,3563)(393,3544)
        (396,3522)(400,3498)(406,3472)
        (413,3445)(421,3416)(432,3386)
        (444,3354)(460,3321)(479,3286)
        (501,3248)(528,3208)(560,3165)
        (597,3120)(639,3072)(687,3022)
        (725,2985)(765,2949)(804,2915)
        (842,2883)(878,2855)(912,2830)
        (944,2808)(973,2790)(999,2774)
        (1024,2762)(1046,2751)(1066,2743)
        (1085,2736)(1103,2731)(1120,2726)
        (1137,2722)(1154,2718)(1171,2713)
        (1189,2708)(1208,2701)(1228,2693)
        (1250,2682)(1275,2670)(1301,2654)
        (1330,2636)(1362,2614)(1396,2589)
        (1432,2561)(1470,2529)(1509,2495)
        (1549,2459)(1587,2422)(1635,2372)
        (1677,2324)(1714,2279)(1746,2236)
        (1773,2196)(1795,2158)(1814,2123)
        (1830,2090)(1842,2058)(1853,2028)
        (1861,1999)(1868,1972)(1874,1946)
        (1878,1922)(1881,1900)(1883,1881)
        (1885,1864)(1886,1850)(1886,1839)
        (1887,1831)(1887,1826)(1887,1823)(1887,1822)
\path(6312,1822)(6312,1821)(6312,1818)
        (6312,1813)(6313,1805)(6313,1794)
        (6314,1780)(6316,1763)(6318,1744)
        (6321,1722)(6325,1698)(6331,1672)
        (6338,1645)(6346,1616)(6357,1586)
        (6369,1554)(6385,1521)(6404,1486)
        (6426,1448)(6453,1408)(6485,1365)
        (6522,1320)(6564,1272)(6612,1222)
        (6650,1185)(6690,1149)(6729,1115)
        (6767,1083)(6803,1055)(6837,1030)
        (6869,1008)(6898,990)(6924,974)
        (6949,962)(6971,951)(6991,943)
        (7010,936)(7028,931)(7045,926)
        (7062,922)(7079,918)(7096,913)
        (7114,908)(7133,901)(7153,893)
        (7175,882)(7200,870)(7226,854)
        (7255,836)(7287,814)(7321,789)
        (7357,761)(7395,729)(7434,695)
        (7474,659)(7512,622)(7560,572)
        (7602,524)(7639,479)(7671,436)
        (7698,396)(7720,358)(7739,323)
        (7755,290)(7767,258)(7778,228)
        (7786,199)(7793,172)(7799,146)
        (7803,122)(7806,100)(7808,81)
        (7810,64)(7811,50)(7811,39)
        (7812,31)(7812,26)(7812,23)(7812,22)
\put(762,3697){\makebox(0,0)[lb]{\smash{{{\SetFigFont{6}{7.2}{\rmdefault}{\mddefault}{\updefault}$\gamma\alpha\beta\gamma^{-1}$}}}}}
\put(2262,1897){\makebox(0,0)[lb]{\smash{{{\SetFigFont{6}{7.2}{\rmdefault}{\mddefault}{\updefault}$\gamma\beta\gamma^{-1}$}}}}}
\put(687,1822){\makebox(0,0)[lb]{\smash{{{\SetFigFont{6}{7.2}{\rmdefault}{\mddefault}{\updefault}$\gamma\alpha\gamma^{-1}$}}}}}
\put(687,97){\makebox(0,0)[lb]{\smash{{{\SetFigFont{6}{7.2}{\rmdefault}{\mddefault}{\updefault}$\alpha$}}}}}
\put(2262,97){\makebox(0,0)[lb]{\smash{{{\SetFigFont{6}{7.2}{\rmdefault}{\mddefault}{\updefault}$\beta$}}}}}
\put(387,922){\makebox(0,0)[lb]{\smash{{{\SetFigFont{6}{7.2}{\rmdefault}{\mddefault}{\updefault}$\gamma$}}}}}
\put(1962,922){\makebox(0,0)[lb]{\smash{{{\SetFigFont{6}{7.2}{\rmdefault}{\mddefault}{\updefault}$\gamma$}}}}}
\put(387,2572){\makebox(0,0)[lb]{\smash{{{\SetFigFont{6}{7.2}{\rmdefault}{\mddefault}{\updefault}$1$}}}}}
\put(1062,2797){\makebox(0,0)[lb]{\smash{{{\SetFigFont{6}{7.2}{\rmdefault}{\mddefault}{\updefault}$1$}}}}}
\put(6387,2722){\makebox(0,0)[lb]{\smash{{{\SetFigFont{6}{7.2}{\rmdefault}{\mddefault}{\updefault}$\gamma$}}}}}
\put(6387,772){\makebox(0,0)[lb]{\smash{{{\SetFigFont{6}{7.2}{\rmdefault}{\mddefault}{\updefault}$1$}}}}}
\put(7062,922){\makebox(0,0)[lb]{\smash{{{\SetFigFont{6}{7.2}{\rmdefault}{\mddefault}{\updefault}$1$}}}}}
\put(6837,1897){\makebox(0,0)[lb]{\smash{{{\SetFigFont{6}{7.2}{\rmdefault}{\mddefault}{\updefault}$\alpha\beta$}}}}}
\put(6837,3697){\makebox(0,0)[lb]{\smash{{{\SetFigFont{6}{7.2}{\rmdefault}{\mddefault}{\updefault}$\gamma\alpha\beta\gamma^{-1}$}}}}}
\put(6687,97){\makebox(0,0)[lb]{\smash{{{\SetFigFont{6}{7.2}{\rmdefault}{\mddefault}{\updefault}$\alpha$}}}}}
\put(8262,97){\makebox(0,0)[lb]{\smash{{{\SetFigFont{6}{7.2}{\rmdefault}{\mddefault}{\updefault}$\beta$}}}}}
\end{picture}
}

\end{center}
\vspace{0.3cm}

Furthermore, $\varphi$ is a group homomorphism by construction, and the
monoidal unit and other structure is preserved.

\begin{bfseries}Braiding.\end{bfseries} Recall that the monoidal
structure comes from 
$*_{\alpha, \beta} = F({p_{\alpha,\beta}})$ for choices of
cobordisms ${p_{\alpha,\beta}}\colon P \rightarrow X$.  
Let $T$ be an untwisting diffeomorphism as pictured

\vspace{0.3cm}
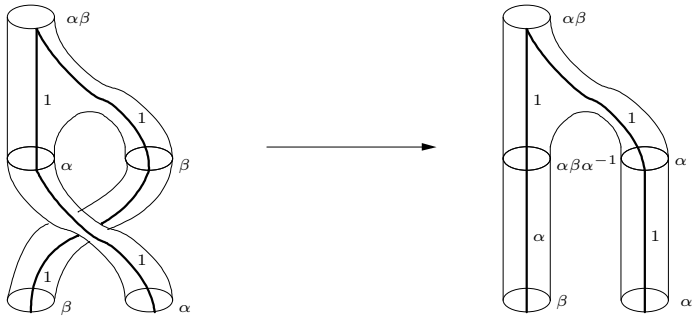
\begin{center}
\setlength{\unitlength}{0.00041667in}
\begingroup\makeatletter\ifx\SetFigFont\undefined%
\gdef\SetFigFont#1#2#3#4#5{%
  \reset@font\fontsize{#1}{#2pt}%
  \fontfamily{#3}\fontseries{#4}\fontshape{#5}%
  \selectfont}%
\fi\endgroup%
{\renewcommand{\dashlinestretch}{30}
\begin{picture}(8946,3976)(0,-10)
\put(6612,2004){\ellipse{600}{300}}
\put(6612,204){\ellipse{600}{300}}
\path(6312,204)(6312,2004)
\path(6912,2004)(6912,204)
\put(1812,2004){\ellipse{600}{300}}
\put(312,3804){\ellipse{600}{300}}
\put(312,2004){\ellipse{600}{300}}
\path(612,2004)(612,2154)
\path(1512,2004)(1512,2154)
\path(612,3804)(612,3803)(612,3800)
        (612,3795)(613,3787)(613,3776)
        (614,3762)(616,3745)(618,3726)
        (621,3704)(625,3680)(631,3654)
        (638,3627)(646,3598)(657,3568)
        (669,3536)(685,3503)(704,3468)
        (726,3430)(753,3390)(785,3347)
        (822,3302)(864,3254)(912,3204)
        (950,3167)(990,3131)(1029,3097)
        (1067,3065)(1103,3037)(1137,3012)
        (1169,2990)(1198,2972)(1224,2956)
        (1249,2944)(1271,2933)(1291,2925)
        (1310,2918)(1328,2913)(1345,2908)
        (1362,2904)(1379,2900)(1396,2895)
        (1414,2890)(1433,2883)(1453,2875)
        (1475,2864)(1500,2852)(1526,2836)
        (1555,2818)(1587,2796)(1621,2771)
        (1657,2743)(1695,2711)(1734,2677)
        (1774,2641)(1812,2604)(1860,2554)
        (1902,2506)(1939,2461)(1971,2418)
        (1998,2378)(2020,2340)(2039,2305)
        (2055,2272)(2067,2240)(2078,2210)
        (2086,2181)(2093,2154)(2099,2128)
        (2103,2104)(2106,2082)(2108,2063)
        (2110,2046)(2111,2032)(2111,2021)
        (2112,2013)(2112,2008)(2112,2005)(2112,2004)
\path(12,2004)(12,3804)
\path(612,2154)(612,2156)(613,2160)
        (615,2167)(617,2177)(620,2190)
        (624,2206)(629,2224)(635,2243)
        (642,2264)(650,2285)(660,2307)
        (671,2330)(684,2353)(699,2377)
        (717,2403)(738,2428)(762,2454)
        (785,2476)(807,2495)(827,2511)
        (845,2525)(860,2536)(871,2545)
        (881,2552)(888,2558)(894,2563)
        (900,2567)(905,2570)(912,2574)
        (920,2577)(930,2581)(943,2586)
        (960,2590)(980,2595)(1005,2599)
        (1032,2603)(1062,2604)(1092,2603)
        (1119,2599)(1144,2595)(1164,2590)
        (1181,2586)(1194,2581)(1204,2577)
        (1212,2574)(1219,2570)(1225,2566)
        (1230,2563)(1236,2558)(1243,2552)
        (1253,2545)(1264,2536)(1279,2525)
        (1297,2511)(1317,2495)(1339,2476)
        (1362,2454)(1386,2428)(1407,2403)
        (1425,2377)(1440,2353)(1453,2330)
        (1464,2307)(1474,2285)(1482,2264)
        (1489,2243)(1495,2224)(1500,2206)
        (1504,2190)(1507,2177)(1509,2167)
        (1511,2160)(1512,2156)(1512,2154)
\put(8112,2004){\ellipse{600}{300}}
\put(6612,3804){\ellipse{600}{300}}
\put(6612,2004){\ellipse{600}{300}}
\path(6912,2004)(6912,2154)
\path(7812,2004)(7812,2154)
\path(6912,3804)(6912,3803)(6912,3800)
        (6912,3795)(6913,3787)(6913,3776)
        (6914,3762)(6916,3745)(6918,3726)
        (6921,3704)(6925,3680)(6931,3654)
        (6938,3627)(6946,3598)(6957,3568)
        (6969,3536)(6985,3503)(7004,3468)
        (7026,3430)(7053,3390)(7085,3347)
        (7122,3302)(7164,3254)(7212,3204)
        (7250,3167)(7290,3131)(7329,3097)
        (7367,3065)(7403,3037)(7437,3012)
        (7469,2990)(7498,2972)(7524,2956)
        (7549,2944)(7571,2933)(7591,2925)
        (7610,2918)(7628,2913)(7645,2908)
        (7662,2904)(7679,2900)(7696,2895)
        (7714,2890)(7733,2883)(7753,2875)
        (7775,2864)(7800,2852)(7826,2836)
        (7855,2818)(7887,2796)(7921,2771)
        (7957,2743)(7995,2711)(8034,2677)
        (8074,2641)(8112,2604)(8160,2554)
        (8202,2506)(8239,2461)(8271,2418)
        (8298,2378)(8320,2340)(8339,2305)
        (8355,2272)(8367,2240)(8378,2210)
        (8386,2181)(8393,2154)(8399,2128)
        (8403,2104)(8406,2082)(8408,2063)
        (8410,2046)(8411,2032)(8411,2021)
        (8412,2013)(8412,2008)(8412,2005)(8412,2004)
\path(6312,2004)(6312,3804)
\path(6912,2154)(6912,2156)(6913,2160)
        (6915,2167)(6917,2177)(6920,2190)
        (6924,2206)(6929,2224)(6935,2243)
        (6942,2264)(6950,2285)(6960,2307)
        (6971,2330)(6984,2353)(6999,2377)
        (7017,2403)(7038,2428)(7062,2454)
        (7085,2476)(7107,2495)(7127,2511)
        (7145,2525)(7160,2536)(7171,2545)
        (7181,2552)(7188,2558)(7194,2563)
        (7200,2567)(7205,2570)(7212,2574)
        (7220,2577)(7230,2581)(7243,2586)
        (7260,2590)(7280,2595)(7305,2599)
        (7332,2603)(7362,2604)(7392,2603)
        (7419,2599)(7444,2595)(7464,2590)
        (7481,2586)(7494,2581)(7504,2577)
        (7512,2574)(7519,2570)(7525,2566)
        (7530,2563)(7536,2558)(7543,2552)
        (7553,2545)(7564,2536)(7579,2525)
        (7597,2511)(7617,2495)(7639,2476)
        (7662,2454)(7686,2428)(7707,2403)
        (7725,2377)(7740,2353)(7753,2330)
        (7764,2307)(7774,2285)(7782,2264)
        (7789,2243)(7795,2224)(7800,2206)
        (7804,2190)(7807,2177)(7809,2167)
        (7811,2160)(7812,2156)(7812,2154)
\put(8112,2004){\ellipse{600}{300}}
\put(8112,204){\ellipse{600}{300}}
\path(7812,204)(7812,2004)
\path(8412,2004)(8412,204)
\put(1812,2004){\ellipse{600}{300}}
\put(312,204){\ellipse{600}{300}}
\put(1812,204){\ellipse{600}{300}}
\put(312,2004){\ellipse{600}{300}}
\path(3312,2154)(5412,2154)
\blacken\thicklines
\path(5292.000,2124.000)(5412.000,2154.000)(5292.000,2184.000)(5292.000,2124.000)
\path(6612,54)(6612,3654)
\path(8112,54)(8112,1854)
\path(387,1854)(388,1853)(389,1851)
        (391,1847)(394,1840)(399,1831)
        (405,1820)(412,1806)(421,1790)
        (432,1771)(444,1750)(457,1728)
        (471,1704)(487,1679)(504,1652)
        (522,1625)(542,1596)(564,1566)
        (588,1534)(614,1501)(643,1465)
        (675,1427)(711,1387)(750,1345)
        (792,1300)(837,1254)(877,1214)
        (918,1176)(957,1140)(993,1107)
        (1027,1078)(1058,1053)(1086,1031)
        (1111,1013)(1133,999)(1152,987)
        (1169,977)(1185,970)(1199,964)
        (1212,959)(1225,954)(1237,949)
        (1250,944)(1264,938)(1279,931)
        (1295,921)(1314,909)(1335,895)
        (1359,877)(1385,855)(1414,830)
        (1446,801)(1480,768)(1515,732)
        (1551,694)(1587,654)(1629,604)
        (1667,556)(1700,511)(1729,468)
        (1755,428)(1776,390)(1795,355)
        (1811,322)(1824,290)(1836,260)
        (1846,231)(1854,204)(1861,178)
        (1868,154)(1873,132)(1877,113)
        (1880,96)(1883,82)(1884,71)
        (1886,63)(1886,58)(1887,55)(1887,54)
\path(312,54)(312,55)(312,59)
        (312,65)(312,74)(313,87)
        (313,103)(314,122)(315,145)
        (317,170)(319,197)(321,225)
        (325,255)(329,286)(334,318)
        (341,350)(348,384)(358,419)
        (369,455)(383,492)(398,531)
        (417,572)(438,613)(462,654)
        (488,694)(515,731)(542,764)
        (569,795)(595,823)(622,848)
        (647,870)(673,890)(698,909)
        (723,926)(747,942)(771,956)
        (794,970)(816,982)(837,992)
        (855,1002)(871,1010)(885,1016)
        (895,1021)(903,1025)(908,1027)
        (911,1028)(912,1029)
\path(1212,1179)(1214,1180)(1217,1182)
        (1223,1186)(1233,1191)(1245,1199)
        (1261,1209)(1278,1220)(1298,1233)
        (1320,1247)(1342,1262)(1366,1279)
        (1390,1297)(1415,1316)(1441,1338)
        (1468,1361)(1497,1386)(1526,1415)
        (1557,1446)(1587,1479)(1616,1513)
        (1643,1546)(1666,1577)(1687,1607)
        (1705,1634)(1721,1660)(1735,1685)
        (1748,1709)(1759,1732)(1770,1753)
        (1779,1774)(1787,1792)(1794,1809)
        (1800,1823)(1805,1835)(1808,1844)
        (1810,1849)(1811,1853)(1812,1854)
\path(387,1854)(387,3654)
\path(387,3654)(388,3653)(389,3651)
        (391,3647)(394,3640)(399,3631)
        (405,3620)(412,3606)(421,3590)
        (432,3571)(444,3550)(457,3528)
        (471,3504)(487,3479)(504,3452)
        (522,3425)(542,3396)(564,3366)
        (588,3334)(614,3301)(643,3265)
        (675,3227)(711,3187)(750,3145)
        (792,3100)(837,3054)(877,3014)
        (918,2976)(957,2940)(994,2907)
        (1028,2878)(1059,2853)(1087,2831)
        (1113,2813)(1135,2799)(1155,2787)
        (1173,2777)(1189,2770)(1204,2764)
        (1217,2759)(1231,2754)(1244,2749)
        (1258,2744)(1272,2738)(1287,2731)
        (1304,2721)(1323,2709)(1344,2695)
        (1368,2677)(1394,2655)(1422,2630)
        (1453,2601)(1485,2568)(1519,2532)
        (1554,2494)(1587,2454)(1626,2404)
        (1660,2356)(1689,2311)(1714,2268)
        (1735,2228)(1753,2190)(1767,2155)
        (1779,2122)(1788,2090)(1795,2060)
        (1801,2031)(1805,2004)(1808,1978)
        (1810,1954)(1812,1932)(1812,1913)
        (1813,1896)(1813,1882)(1813,1871)
        (1812,1863)(1812,1858)(1812,1855)(1812,1854)
\path(6612,3654)(6613,3653)(6614,3651)
        (6616,3647)(6619,3640)(6624,3631)
        (6630,3620)(6637,3606)(6646,3590)
        (6657,3571)(6669,3550)(6682,3528)
        (6696,3504)(6712,3479)(6729,3452)
        (6747,3425)(6767,3396)(6789,3366)
        (6813,3334)(6839,3301)(6868,3265)
        (6900,3227)(6936,3187)(6975,3145)
        (7017,3100)(7062,3054)(7102,3014)
        (7143,2976)(7182,2940)(7218,2907)
        (7252,2878)(7283,2853)(7311,2831)
        (7336,2813)(7358,2799)(7377,2787)
        (7394,2777)(7410,2770)(7424,2764)
        (7437,2759)(7450,2754)(7462,2749)
        (7475,2744)(7489,2738)(7504,2731)
        (7520,2721)(7539,2709)(7560,2695)
        (7584,2677)(7610,2655)(7639,2630)
        (7671,2601)(7705,2568)(7740,2532)
        (7776,2494)(7812,2454)(7854,2404)
        (7892,2356)(7925,2311)(7954,2268)
        (7980,2228)(8001,2190)(8020,2155)
        (8036,2122)(8049,2090)(8061,2060)
        (8071,2031)(8079,2004)(8086,1978)
        (8093,1954)(8098,1932)(8102,1913)
        (8105,1896)(8108,1882)(8109,1871)
        (8111,1863)(8111,1858)(8112,1855)(8112,1854)
\thinlines
\path(612,2004)(612,2003)(612,2000)
        (612,1995)(613,1987)(613,1976)
        (614,1962)(616,1945)(618,1926)
        (621,1904)(625,1880)(631,1854)
        (638,1827)(646,1798)(657,1768)
        (669,1736)(685,1703)(704,1668)
        (726,1630)(753,1590)(785,1547)
        (822,1502)(864,1454)(912,1404)
        (950,1367)(990,1331)(1029,1297)
        (1067,1265)(1103,1237)(1137,1212)
        (1169,1190)(1198,1172)(1224,1156)
        (1249,1144)(1271,1133)(1291,1125)
        (1310,1118)(1328,1113)(1345,1108)
        (1362,1104)(1379,1100)(1396,1095)
        (1414,1090)(1433,1083)(1453,1075)
        (1475,1064)(1500,1052)(1526,1036)
        (1555,1018)(1587,996)(1621,971)
        (1657,943)(1695,911)(1734,877)
        (1774,841)(1812,804)(1860,754)
        (1902,706)(1939,661)(1971,618)
        (1998,578)(2020,540)(2039,505)
        (2055,472)(2067,440)(2078,410)
        (2086,381)(2093,354)(2099,328)
        (2103,304)(2106,282)(2108,263)
        (2110,246)(2111,232)(2111,221)
        (2112,213)(2112,208)(2112,205)(2112,204)
\path(12,2004)(12,2003)(12,2000)
        (12,1995)(13,1987)(13,1976)
        (14,1962)(16,1945)(18,1926)
        (21,1904)(25,1880)(31,1854)
        (38,1827)(46,1798)(57,1768)
        (69,1736)(85,1703)(104,1668)
        (126,1630)(153,1590)(185,1547)
        (222,1502)(264,1454)(312,1404)
        (350,1367)(390,1331)(429,1297)
        (467,1265)(503,1237)(537,1212)
        (569,1190)(598,1172)(624,1156)
        (649,1144)(671,1133)(691,1125)
        (710,1118)(728,1113)(745,1108)
        (762,1104)(779,1100)(796,1095)
        (814,1090)(833,1083)(853,1075)
        (875,1064)(900,1052)(926,1036)
        (955,1018)(987,996)(1021,971)
        (1057,943)(1095,911)(1134,877)
        (1174,841)(1212,804)(1260,754)
        (1302,706)(1339,661)(1371,618)
        (1398,578)(1420,540)(1439,505)
        (1455,472)(1467,440)(1478,410)
        (1486,381)(1493,354)(1499,328)
        (1503,304)(1506,282)(1508,263)
        (1510,246)(1511,232)(1511,221)
        (1512,213)(1512,208)(1512,205)(1512,204)
\path(12,204)(12,205)(13,209)
        (14,214)(16,223)(19,236)
        (23,252)(27,271)(33,294)
        (39,319)(46,347)(53,377)
        (61,408)(70,440)(79,473)
        (89,507)(99,541)(110,576)
        (122,612)(135,648)(149,685)
        (164,723)(180,762)(198,801)
        (217,841)(237,879)(263,926)
        (289,967)(314,1003)(337,1033)
        (359,1059)(380,1080)(400,1099)
        (419,1114)(437,1128)(455,1139)
        (471,1149)(487,1157)(500,1164)
        (512,1169)(521,1173)(528,1176)
        (533,1178)(536,1179)(537,1179)
\path(912,1329)(914,1330)(917,1332)
        (923,1335)(933,1341)(946,1348)
        (963,1357)(983,1369)(1005,1382)
        (1030,1397)(1057,1413)(1085,1430)
        (1114,1448)(1144,1466)(1174,1486)
        (1204,1506)(1234,1527)(1264,1549)
        (1293,1572)(1323,1596)(1353,1621)
        (1382,1648)(1411,1676)(1437,1704)
        (1465,1737)(1487,1768)(1505,1797)
        (1518,1822)(1527,1845)(1533,1867)
        (1537,1886)(1538,1904)(1537,1921)
        (1535,1936)(1532,1950)(1529,1963)
        (1525,1974)(1521,1984)(1518,1992)
        (1515,1997)(1514,2001)(1513,2003)(1512,2004)
\path(612,204)(612,206)(613,210)
        (614,217)(616,228)(618,243)
        (621,261)(625,282)(630,305)
        (635,331)(642,357)(649,385)
        (657,414)(666,444)(677,475)
        (689,508)(704,542)(721,579)
        (740,616)(762,654)(785,690)
        (808,724)(831,754)(853,781)
        (875,805)(896,826)(916,845)
        (935,863)(954,879)(973,893)
        (991,906)(1007,918)(1022,928)
        (1035,936)(1045,943)(1053,948)
        (1058,951)(1061,953)(1062,954)
\path(1287,1104)(1288,1105)(1291,1106)
        (1297,1108)(1305,1111)(1316,1115)
        (1331,1121)(1348,1129)(1369,1137)
        (1392,1147)(1417,1158)(1444,1171)
        (1472,1184)(1501,1198)(1531,1213)
        (1561,1229)(1592,1246)(1623,1265)
        (1655,1285)(1687,1306)(1720,1329)
        (1753,1355)(1787,1383)(1821,1413)
        (1854,1445)(1887,1479)(1923,1520)
        (1954,1561)(1982,1600)(2005,1638)
        (2025,1674)(2042,1708)(2055,1740)
        (2067,1771)(2076,1802)(2084,1831)
        (2091,1859)(2096,1885)(2101,1910)
        (2104,1933)(2107,1953)(2109,1969)
        (2110,1983)(2111,1992)(2112,1999)
        (2112,2002)(2112,2004)
\put(2187,54){\makebox(0,0)[lb]{\smash{{{\SetFigFont{6}{7.2}{\rmdefault}{\mddefault}{\updefault}$\alpha$}}}}}
\put(687,54){\makebox(0,0)[lb]{\smash{{{\SetFigFont{6}{7.2}{\rmdefault}{\mddefault}{\updefault}$\beta$}}}}}
\put(2187,1854){\makebox(0,0)[lb]{\smash{{{\SetFigFont{6}{7.2}{\rmdefault}{\mddefault}{\updefault}$\beta$}}}}}
\put(762,3729){\makebox(0,0)[lb]{\smash{{{\SetFigFont{6}{7.2}{\rmdefault}{\mddefault}{\updefault}$\alpha\beta$}}}}}
\put(687,1854){\makebox(0,0)[lb]{\smash{{{\SetFigFont{6}{7.2}{\rmdefault}{\mddefault}{\updefault}$\alpha$}}}}}
\put(6987,129){\makebox(0,0)[lb]{\smash{{{\SetFigFont{6}{7.2}{\rmdefault}{\mddefault}{\updefault}$\beta$}}}}}
\put(8562,129){\makebox(0,0)[lb]{\smash{{{\SetFigFont{6}{7.2}{\rmdefault}{\mddefault}{\updefault}$\alpha$}}}}}
\put(6987,1854){\makebox(0,0)[lb]{\smash{{{\SetFigFont{6}{7.2}{\rmdefault}{\mddefault}{\updefault}$\alpha\beta\alpha^{-1}$}}}}}
\put(8487,1929){\makebox(0,0)[lb]{\smash{{{\SetFigFont{6}{7.2}{\rmdefault}{\mddefault}{\updefault}$\alpha$}}}}}
\put(7062,3729){\makebox(0,0)[lb]{\smash{{{\SetFigFont{6}{7.2}{\rmdefault}{\mddefault}{\updefault}$\alpha\beta$}}}}}
\put(462,2679){\makebox(0,0)[lb]{\smash{{{\SetFigFont{6}{7.2}{\rmdefault}{\mddefault}{\updefault}$1$}}}}}
\put(1662,2454){\makebox(0,0)[lb]{\smash{{{\SetFigFont{6}{7.2}{\rmdefault}{\mddefault}{\updefault}$1$}}}}}
\put(1662,654){\makebox(0,0)[lb]{\smash{{{\SetFigFont{6}{7.2}{\rmdefault}{\mddefault}{\updefault}$1$}}}}}
\put(462,429){\makebox(0,0)[lb]{\smash{{{\SetFigFont{6}{7.2}{\rmdefault}{\mddefault}{\updefault}$1$}}}}}
\put(6687,2679){\makebox(0,0)[lb]{\smash{{{\SetFigFont{6}{7.2}{\rmdefault}{\mddefault}{\updefault}$1$}}}}}
\put(7887,2454){\makebox(0,0)[lb]{\smash{{{\SetFigFont{6}{7.2}{\rmdefault}{\mddefault}{\updefault}$1$}}}}}
\put(6687,954){\makebox(0,0)[lb]{\smash{{{\SetFigFont{6}{7.2}{\rmdefault}{\mddefault}{\updefault}$\alpha$}}}}}
\put(8187,954){\makebox(0,0)[lb]{\smash{{{\SetFigFont{6}{7.2}{\rmdefault}{\mddefault}{\updefault}$1$}}}}}
\end{picture}
}

\end{center}
\vspace{0.3cm}

Letting $\tau_{\alpha,\beta}\colon J \rightarrow X$ be the crossed
cylinders in the picture above and
recalling that $F(\tau_{\alpha, \beta})$ is the functor changing
components it follows that for each $\alpha$ and $\beta$ there is a
2-morphism
\[
\xymatrix{P \circ J \ar[rr]^{T} \ar[dr]_{p_{\alpha,\beta}\circ
  \tau_{\beta,\alpha}} & & P\circ (I\sqcup I)
  \ar[dl]^{p_{\alpha\beta\alpha^{-1},\alpha}\circ (i_{\beta,\alpha}
  \sqcup i_{\alpha,1})}\\ & X}
\]
This induces a natural isomorphism (of functors $\cA_\beta \otimes
\cA_\alpha \rightarrow \cA_{\alpha\beta}$)
\[
s_{\alpha,\beta}\colon *_{\alpha,\beta} \circ
\mbox{twist}_{\beta,\alpha} \nattrans 
*_{\alpha\beta\alpha^{-1},\alpha} \circ (\varphi(\alpha)_\beta \otimes
id)
\]
Thus for each $U\in \cA_\alpha$ and $V\in \cA_\beta$ we get
isomorphisms 
\[
s_{U,V}\colon U*V \rightarrow \varphi(\alpha)V*U.
\]
Since $\varphi$ is a homomorphism into the group of $*$ preserving
automorphisms of $\cA$, the braiding result in \cite{Tillmann:S-Structures} now
implies condition (1) of Definition \ref{defn:braid}.  Condition (2)
follows immediately from the naturality of $s$.  

Condition (3) is
obtained as follows. By post-composing the two functors above with
$\varphi (\gamma)$ we get two functors $\cA_\beta \otimes \cA_\alpha
\rightarrow \cA_{\gamma\alpha\beta\gamma^{-1}}$. The natural
transformation $s$ induces a natural transformation of these two
functors which for objects $U\in\cA_\alpha$ and $V\in\cA_\beta$ is
given by $\varphi(\gamma)(s_{U,V})$. On the other hand, by
pre-composing by $\varphi(\gamma)$ we again get an induced natural
transformation which is given by
$s_{\varphi(\gamma)U,\varphi(\gamma)V}$ for objects $U\in\cA_\alpha$
and $V\in\cA_\beta$. However, the diffeomorphisms inducing the natural
transformation above are
the same as 2-morphisms in $\cS_X$ so the induced natural
transformations are the same i.e. $\varphi(\gamma)(s_{U,V}) =
s_{\varphi(\gamma)U,\varphi(\gamma)V}$, proving condition (3).

\begin{bfseries}Twist.\end{bfseries} The
Dehn twist of a cylinder  $D\colon I 
\rightarrow I$ gives for each $\alpha, \beta$ a 2-morphism
\[
\xymatrix{I \ar[rr]^{D} \ar[dr]_{i_{\alpha,\beta}} & &  I
  \ar[dl]^{i_{\alpha,\beta\alpha}}\\ & X}
\]
which induces a natural isomorphism $\varphi(\beta)_\alpha
 \nattrans 
\varphi(\beta\alpha)_\alpha $, which taking 
 $\beta = 1$ gives a natural isomorphism
\[
\theta_\alpha\colon \varphi(1)_\alpha
 \nattrans 
\varphi(\alpha)_\alpha
\]
that is, for each $U\in \cA_\alpha$ there is an isomorphism
$\theta_U\colon U \rightarrow \varphi(\alpha)U$. 
Conditions (1) and (2) of Definition \ref{defn:twist} are satisfied 
 by the equivalent mapping
class group 
identities of \cite{Tillmann:S-Structures} and the homotopy identity displayed 
below (where we represent a half twist by swapping the labels
of the incoming boundary components). 

\vspace{0.3cm}
\begin{center}
\setlength{\unitlength}{0.00041667in}
\begingroup\makeatletter\ifx\SetFigFont\undefined%
\gdef\SetFigFont#1#2#3#4#5{%
  \reset@font\fontsize{#1}{#2pt}%
  \fontfamily{#3}\fontseries{#4}\fontshape{#5}%
  \selectfont}%
\fi\endgroup%
{\renewcommand{\dashlinestretch}{30}
\begin{picture}(7959,6976)(0,-10)
\put(1950,204){\ellipse{600}{300}}
\put(450,2004){\ellipse{600}{300}}
\put(450,204){\ellipse{600}{300}}
\path(750,204)(750,354)
\path(1650,204)(1650,354)
\path(750,2004)(750,2003)(750,2000)
        (750,1995)(751,1987)(751,1976)
        (752,1962)(754,1945)(756,1926)
        (759,1904)(763,1880)(769,1854)
        (776,1827)(784,1798)(795,1768)
        (807,1736)(823,1703)(842,1668)
        (864,1630)(891,1590)(923,1547)
        (960,1502)(1002,1454)(1050,1404)
        (1088,1367)(1128,1331)(1167,1297)
        (1205,1265)(1241,1237)(1275,1212)
        (1307,1190)(1336,1172)(1362,1156)
        (1387,1144)(1409,1133)(1429,1125)
        (1448,1118)(1466,1113)(1483,1108)
        (1500,1104)(1517,1100)(1534,1095)
        (1552,1090)(1571,1083)(1591,1075)
        (1613,1064)(1638,1052)(1664,1036)
        (1693,1018)(1725,996)(1759,971)
        (1795,943)(1833,911)(1872,877)
        (1912,841)(1950,804)(1998,754)
        (2040,706)(2077,661)(2109,618)
        (2136,578)(2158,540)(2177,505)
        (2193,472)(2205,440)(2216,410)
        (2224,381)(2231,354)(2237,328)
        (2241,304)(2244,282)(2246,263)
        (2248,246)(2249,232)(2249,221)
        (2250,213)(2250,208)(2250,205)(2250,204)
\path(150,204)(150,2004)
\path(750,354)(750,356)(751,360)
        (753,367)(755,377)(758,390)
        (762,406)(767,424)(773,443)
        (780,464)(788,485)(798,507)
        (809,530)(822,553)(837,577)
        (855,603)(876,628)(900,654)
        (923,676)(945,695)(965,711)
        (983,725)(998,736)(1009,745)
        (1019,752)(1026,758)(1032,763)
        (1038,767)(1043,770)(1050,774)
        (1058,777)(1068,781)(1081,786)
        (1098,790)(1118,795)(1143,799)
        (1170,803)(1200,804)(1230,803)
        (1257,799)(1282,795)(1302,790)
        (1319,786)(1332,781)(1342,777)
        (1350,774)(1357,770)(1363,766)
        (1368,763)(1374,758)(1381,752)
        (1391,745)(1402,736)(1417,725)
        (1435,711)(1455,695)(1477,676)
        (1500,654)(1524,628)(1545,603)
        (1563,577)(1578,553)(1591,530)
        (1602,507)(1612,485)(1620,464)
        (1627,443)(1633,424)(1638,406)
        (1642,390)(1645,377)(1647,367)
        (1649,360)(1650,356)(1650,354)
\thicklines
\path(525,1779)(525,54)
\path(525,1854)(525,1853)(525,1850)
        (525,1845)(526,1837)(526,1826)
        (527,1812)(529,1795)(531,1776)
        (534,1754)(538,1730)(544,1704)
        (551,1677)(559,1648)(570,1618)
        (582,1586)(598,1553)(617,1518)
        (639,1480)(666,1440)(698,1397)
        (735,1352)(777,1304)(825,1254)
        (863,1217)(903,1181)(942,1147)
        (980,1115)(1016,1087)(1050,1062)
        (1082,1040)(1111,1022)(1137,1006)
        (1162,994)(1184,983)(1204,975)
        (1223,968)(1241,963)(1258,958)
        (1275,954)(1292,950)(1309,945)
        (1327,940)(1346,933)(1366,925)
        (1388,914)(1413,902)(1439,886)
        (1468,868)(1500,846)(1534,821)
        (1570,793)(1608,761)(1647,727)
        (1687,691)(1725,654)(1773,604)
        (1815,556)(1852,511)(1884,468)
        (1911,428)(1933,390)(1952,355)
        (1968,322)(1980,290)(1991,260)
        (1999,231)(2006,204)(2012,178)
        (2016,154)(2019,132)(2021,113)
        (2023,96)(2024,82)(2024,71)
        (2025,63)(2025,58)(2025,55)(2025,54)
\thinlines
\put(1950,5004){\ellipse{600}{300}}
\put(450,6804){\ellipse{600}{300}}
\put(450,5004){\ellipse{600}{300}}
\path(750,5004)(750,5154)
\path(1650,5004)(1650,5154)
\path(750,6804)(750,6803)(750,6800)
        (750,6795)(751,6787)(751,6776)
        (752,6762)(754,6745)(756,6726)
        (759,6704)(763,6680)(769,6654)
        (776,6627)(784,6598)(795,6568)
        (807,6536)(823,6503)(842,6468)
        (864,6430)(891,6390)(923,6347)
        (960,6302)(1002,6254)(1050,6204)
        (1088,6167)(1128,6131)(1167,6097)
        (1205,6065)(1241,6037)(1275,6012)
        (1307,5990)(1336,5972)(1362,5956)
        (1387,5944)(1409,5933)(1429,5925)
        (1448,5918)(1466,5913)(1483,5908)
        (1500,5904)(1517,5900)(1534,5895)
        (1552,5890)(1571,5883)(1591,5875)
        (1613,5864)(1638,5852)(1664,5836)
        (1693,5818)(1725,5796)(1759,5771)
        (1795,5743)(1833,5711)(1872,5677)
        (1912,5641)(1950,5604)(1998,5554)
        (2040,5506)(2077,5461)(2109,5418)
        (2136,5378)(2158,5340)(2177,5305)
        (2193,5272)(2205,5240)(2216,5210)
        (2224,5181)(2231,5154)(2237,5128)
        (2241,5104)(2244,5082)(2246,5063)
        (2248,5046)(2249,5032)(2249,5021)
        (2250,5013)(2250,5008)(2250,5005)(2250,5004)
\path(150,5004)(150,6804)
\path(750,5154)(750,5156)(751,5160)
        (753,5167)(755,5177)(758,5190)
        (762,5206)(767,5224)(773,5243)
        (780,5264)(788,5285)(798,5307)
        (809,5330)(822,5353)(837,5377)
        (855,5403)(876,5428)(900,5454)
        (923,5476)(945,5495)(965,5511)
        (983,5525)(998,5536)(1009,5545)
        (1019,5552)(1026,5558)(1032,5563)
        (1038,5567)(1043,5570)(1050,5574)
        (1058,5577)(1068,5581)(1081,5586)
        (1098,5590)(1118,5595)(1143,5599)
        (1170,5603)(1200,5604)(1230,5603)
        (1257,5599)(1282,5595)(1302,5590)
        (1319,5586)(1332,5581)(1342,5577)
        (1350,5574)(1357,5570)(1363,5566)
        (1368,5563)(1374,5558)(1381,5552)
        (1391,5545)(1402,5536)(1417,5525)
        (1435,5511)(1455,5495)(1477,5476)
        (1500,5454)(1524,5428)(1545,5403)
        (1563,5377)(1578,5353)(1591,5330)
        (1602,5307)(1612,5285)(1620,5264)
        (1627,5243)(1633,5224)(1638,5206)
        (1642,5190)(1645,5177)(1647,5167)
        (1649,5160)(1650,5156)(1650,5154)
\thicklines
\path(525,6579)(525,4854)
\path(525,6654)(525,6653)(525,6650)
        (525,6645)(526,6637)(526,6626)
        (527,6612)(529,6595)(531,6576)
        (534,6554)(538,6530)(544,6504)
        (551,6477)(559,6448)(570,6418)
        (582,6386)(598,6353)(617,6318)
        (639,6280)(666,6240)(698,6197)
        (735,6152)(777,6104)(825,6054)
        (863,6017)(903,5981)(942,5947)
        (980,5915)(1016,5887)(1050,5862)
        (1082,5840)(1111,5822)(1137,5806)
        (1162,5794)(1184,5783)(1204,5775)
        (1223,5768)(1241,5763)(1258,5758)
        (1275,5754)(1292,5750)(1309,5745)
        (1327,5740)(1346,5733)(1366,5725)
        (1388,5714)(1413,5702)(1439,5686)
        (1468,5668)(1500,5646)(1534,5621)
        (1570,5593)(1608,5561)(1647,5527)
        (1687,5491)(1725,5454)(1773,5404)
        (1815,5356)(1852,5311)(1884,5268)
        (1911,5228)(1933,5190)(1952,5155)
        (1968,5122)(1980,5090)(1991,5060)
        (1999,5031)(2006,5004)(2012,4978)
        (2016,4954)(2019,4932)(2021,4913)
        (2023,4896)(2024,4882)(2024,4871)
        (2025,4863)(2025,4858)(2025,4855)(2025,4854)
\thinlines
\put(7050,204){\ellipse{600}{300}}
\put(5550,2004){\ellipse{600}{300}}
\put(5550,204){\ellipse{600}{300}}
\path(5850,204)(5850,354)
\path(6750,204)(6750,354)
\path(5850,2004)(5850,2003)(5850,2000)
        (5850,1995)(5851,1987)(5851,1976)
        (5852,1962)(5854,1945)(5856,1926)
        (5859,1904)(5863,1880)(5869,1854)
        (5876,1827)(5884,1798)(5895,1768)
        (5907,1736)(5923,1703)(5942,1668)
        (5964,1630)(5991,1590)(6023,1547)
        (6060,1502)(6102,1454)(6150,1404)
        (6188,1367)(6228,1331)(6267,1297)
        (6305,1265)(6341,1237)(6375,1212)
        (6407,1190)(6436,1172)(6462,1156)
        (6487,1144)(6509,1133)(6529,1125)
        (6548,1118)(6566,1113)(6583,1108)
        (6600,1104)(6617,1100)(6634,1095)
        (6652,1090)(6671,1083)(6691,1075)
        (6713,1064)(6738,1052)(6764,1036)
        (6793,1018)(6825,996)(6859,971)
        (6895,943)(6933,911)(6972,877)
        (7012,841)(7050,804)(7098,754)
        (7140,706)(7177,661)(7209,618)
        (7236,578)(7258,540)(7277,505)
        (7293,472)(7305,440)(7316,410)
        (7324,381)(7331,354)(7337,328)
        (7341,304)(7344,282)(7346,263)
        (7348,246)(7349,232)(7349,221)
        (7350,213)(7350,208)(7350,205)(7350,204)
\path(5250,204)(5250,2004)
\path(5850,354)(5850,356)(5851,360)
        (5853,367)(5855,377)(5858,390)
        (5862,406)(5867,424)(5873,443)
        (5880,464)(5888,485)(5898,507)
        (5909,530)(5922,553)(5937,577)
        (5955,603)(5976,628)(6000,654)
        (6023,676)(6045,695)(6065,711)
        (6083,725)(6098,736)(6109,745)
        (6119,752)(6126,758)(6132,763)
        (6138,767)(6143,770)(6150,774)
        (6158,777)(6168,781)(6181,786)
        (6198,790)(6218,795)(6243,799)
        (6270,803)(6300,804)(6330,803)
        (6357,799)(6382,795)(6402,790)
        (6419,786)(6432,781)(6442,777)
        (6450,774)(6457,770)(6463,766)
        (6468,763)(6474,758)(6481,752)
        (6491,745)(6502,736)(6517,725)
        (6535,711)(6555,695)(6577,676)
        (6600,654)(6624,628)(6645,603)
        (6663,577)(6678,553)(6691,530)
        (6702,507)(6712,485)(6720,464)
        (6727,443)(6733,424)(6738,406)
        (6742,390)(6745,377)(6747,367)
        (6749,360)(6750,356)(6750,354)
\thicklines
\path(5625,1779)(5625,54)
\path(5625,1854)(5625,1853)(5625,1850)
        (5625,1845)(5626,1837)(5626,1826)
        (5627,1812)(5629,1795)(5631,1776)
        (5634,1754)(5638,1730)(5644,1704)
        (5651,1677)(5659,1648)(5670,1618)
        (5682,1586)(5698,1553)(5717,1518)
        (5739,1480)(5766,1440)(5798,1397)
        (5835,1352)(5877,1304)(5925,1254)
        (5963,1217)(6003,1181)(6042,1147)
        (6080,1115)(6116,1087)(6150,1062)
        (6182,1040)(6211,1022)(6237,1006)
        (6262,994)(6284,983)(6304,975)
        (6323,968)(6341,963)(6358,958)
        (6375,954)(6392,950)(6409,945)
        (6427,940)(6446,933)(6466,925)
        (6488,914)(6513,902)(6539,886)
        (6568,868)(6600,846)(6634,821)
        (6670,793)(6708,761)(6747,727)
        (6787,691)(6825,654)(6873,604)
        (6915,556)(6952,511)(6984,468)
        (7011,428)(7033,390)(7052,355)
        (7068,322)(7080,290)(7091,260)
        (7099,231)(7106,204)(7112,178)
        (7116,154)(7119,132)(7121,113)
        (7123,96)(7124,82)(7124,71)
        (7125,63)(7125,58)(7125,55)(7125,54)
\thinlines
\put(7050,5004){\ellipse{600}{300}}
\put(5550,6804){\ellipse{600}{300}}
\put(5550,5004){\ellipse{600}{300}}
\path(5850,5004)(5850,5154)
\path(6750,5004)(6750,5154)
\path(5850,6804)(5850,6803)(5850,6800)
        (5850,6795)(5851,6787)(5851,6776)
        (5852,6762)(5854,6745)(5856,6726)
        (5859,6704)(5863,6680)(5869,6654)
        (5876,6627)(5884,6598)(5895,6568)
        (5907,6536)(5923,6503)(5942,6468)
        (5964,6430)(5991,6390)(6023,6347)
        (6060,6302)(6102,6254)(6150,6204)
        (6188,6167)(6228,6131)(6267,6097)
        (6305,6065)(6341,6037)(6375,6012)
        (6407,5990)(6436,5972)(6462,5956)
        (6487,5944)(6509,5933)(6529,5925)
        (6548,5918)(6566,5913)(6583,5908)
        (6600,5904)(6617,5900)(6634,5895)
        (6652,5890)(6671,5883)(6691,5875)
        (6713,5864)(6738,5852)(6764,5836)
        (6793,5818)(6825,5796)(6859,5771)
        (6895,5743)(6933,5711)(6972,5677)
        (7012,5641)(7050,5604)(7098,5554)
        (7140,5506)(7177,5461)(7209,5418)
        (7236,5378)(7258,5340)(7277,5305)
        (7293,5272)(7305,5240)(7316,5210)
        (7324,5181)(7331,5154)(7337,5128)
        (7341,5104)(7344,5082)(7346,5063)
        (7348,5046)(7349,5032)(7349,5021)
        (7350,5013)(7350,5008)(7350,5005)(7350,5004)
\path(5250,5004)(5250,6804)
\path(5850,5154)(5850,5156)(5851,5160)
        (5853,5167)(5855,5177)(5858,5190)
        (5862,5206)(5867,5224)(5873,5243)
        (5880,5264)(5888,5285)(5898,5307)
        (5909,5330)(5922,5353)(5937,5377)
        (5955,5403)(5976,5428)(6000,5454)
        (6023,5476)(6045,5495)(6065,5511)
        (6083,5525)(6098,5536)(6109,5545)
        (6119,5552)(6126,5558)(6132,5563)
        (6138,5567)(6143,5570)(6150,5574)
        (6158,5577)(6168,5581)(6181,5586)
        (6198,5590)(6218,5595)(6243,5599)
        (6270,5603)(6300,5604)(6330,5603)
        (6357,5599)(6382,5595)(6402,5590)
        (6419,5586)(6432,5581)(6442,5577)
        (6450,5574)(6457,5570)(6463,5566)
        (6468,5563)(6474,5558)(6481,5552)
        (6491,5545)(6502,5536)(6517,5525)
        (6535,5511)(6555,5495)(6577,5476)
        (6600,5454)(6624,5428)(6645,5403)
        (6663,5377)(6678,5353)(6691,5330)
        (6702,5307)(6712,5285)(6720,5264)
        (6727,5243)(6733,5224)(6738,5206)
        (6742,5190)(6745,5177)(6747,5167)
        (6749,5160)(6750,5156)(6750,5154)
\thicklines
\path(5625,6579)(5625,4854)
\path(5625,6654)(5625,6653)(5625,6650)
        (5625,6645)(5626,6637)(5626,6626)
        (5627,6612)(5629,6595)(5631,6576)
        (5634,6554)(5638,6530)(5644,6504)
        (5651,6477)(5659,6448)(5670,6418)
        (5682,6386)(5698,6353)(5717,6318)
        (5739,6280)(5766,6240)(5798,6197)
        (5835,6152)(5877,6104)(5925,6054)
        (5963,6017)(6003,5981)(6042,5947)
        (6080,5915)(6116,5887)(6150,5862)
        (6182,5840)(6211,5822)(6237,5806)
        (6262,5794)(6284,5783)(6304,5775)
        (6323,5768)(6341,5763)(6358,5758)
        (6375,5754)(6392,5750)(6409,5745)
        (6427,5740)(6446,5733)(6466,5725)
        (6488,5714)(6513,5702)(6539,5686)
        (6568,5668)(6600,5646)(6634,5621)
        (6670,5593)(6708,5561)(6747,5527)
        (6787,5491)(6825,5454)(6873,5404)
        (6915,5356)(6952,5311)(6984,5268)
        (7011,5228)(7033,5190)(7052,5155)
        (7068,5122)(7080,5090)(7091,5060)
        (7099,5031)(7106,5004)(7112,4978)
        (7116,4954)(7119,4932)(7121,4913)
        (7123,4896)(7124,4882)(7124,4871)
        (7125,4863)(7125,4858)(7125,4855)(7125,4854)
\thinlines
\path(5550,4554)(5550,2454)
\blacken\thicklines
\path(5520.000,2574.000)(5550.000,2454.000)(5580.000,2574.000)(5520.000,2574.000)
\put(5925,3354){\makebox(0,0)[lb]{\smash{{{\SetFigFont{6}{7.2}{\rmdefault}{\mddefault}{\updefault}$\theta_V*\theta_U$}}}}}
\put(6000,1929){\makebox(0,0)[lb]{\smash{{{\SetFigFont{6}{7.2}{\rmdefault}{\mddefault}{\updefault}$\alpha\beta$}}}}}
\put(5925,129){\makebox(0,0)[lb]{\smash{{{\SetFigFont{6}{7.2}{\rmdefault}{\mddefault}{\updefault}$\beta$}}}}}
\put(7500,54){\makebox(0,0)[lb]{\smash{{{\SetFigFont{6}{7.2}{\rmdefault}{\mddefault}{\updefault}$\alpha$}}}}}
\put(5400,804){\makebox(0,0)[lb]{\smash{{{\SetFigFont{6}{7.2}{\rmdefault}{\mddefault}{\updefault}$\alpha\beta$}}}}}
\put(6825,729){\makebox(0,0)[lb]{\smash{{{\SetFigFont{6}{7.2}{\rmdefault}{\mddefault}{\updefault}$\alpha$}}}}}
\put(5400,5529){\makebox(0,0)[lb]{\smash{{{\SetFigFont{6}{7.2}{\rmdefault}{\mddefault}{\updefault}$\alpha$}}}}}
\put(6975,5379){\makebox(0,0)[lb]{\smash{{{\SetFigFont{6}{7.2}{\rmdefault}{\mddefault}{\updefault}$1$}}}}}
\put(6000,6729){\makebox(0,0)[lb]{\smash{{{\SetFigFont{6}{7.2}{\rmdefault}{\mddefault}{\updefault}$\alpha\beta$}}}}}
\put(6000,4854){\makebox(0,0)[lb]{\smash{{{\SetFigFont{6}{7.2}{\rmdefault}{\mddefault}{\updefault}$\beta$}}}}}
\put(7575,4854){\makebox(0,0)[lb]{\smash{{{\SetFigFont{6}{7.2}{\rmdefault}{\mddefault}{\updefault}$\alpha$}}}}}
\thinlines
\path(1050,4554)(1050,2454)
\blacken\thicklines
\path(1020.000,2574.000)(1050.000,2454.000)(1080.000,2574.000)(1020.000,2574.000)
\thinlines
\path(4650,1254)(2550,1254)
\blacken\thicklines
\path(2670.000,1284.000)(2550.000,1254.000)(2670.000,1224.000)(2670.000,1284.000)
\thinlines
\path(2550,5979)(4650,5979)
\blacken\thicklines
\path(4530.000,5949.000)(4650.000,5979.000)(4530.000,6009.000)(4530.000,5949.000)
\put(0,3354){\makebox(0,0)[lb]{\smash{{{\SetFigFont{6}{7.2}{\rmdefault}{\mddefault}{\updefault}$\theta_{U*V}$}}}}}
\put(900,1929){\makebox(0,0)[lb]{\smash{{{\SetFigFont{6}{7.2}{\rmdefault}{\mddefault}{\updefault}$\alpha\beta$}}}}}
\put(825,129){\makebox(0,0)[lb]{\smash{{{\SetFigFont{6}{7.2}{\rmdefault}{\mddefault}{\updefault}$\alpha$}}}}}
\put(2400,54){\makebox(0,0)[lb]{\smash{{{\SetFigFont{6}{7.2}{\rmdefault}{\mddefault}{\updefault}$\beta$}}}}}
\put(300,804){\makebox(0,0)[lb]{\smash{{{\SetFigFont{6}{7.2}{\rmdefault}{\mddefault}{\updefault}$\alpha\beta$}}}}}
\put(1800,654){\makebox(0,0)[lb]{\smash{{{\SetFigFont{6}{7.2}{\rmdefault}{\mddefault}{\updefault}$\alpha\beta$}}}}}
\put(3825,579){\makebox(0,0)[lb]{\smash{{{\SetFigFont{6}{7.2}{\rmdefault}{\mddefault}{\updefault}$s_{V,U}$}}}}}
\put(900,6729){\makebox(0,0)[lb]{\smash{{{\SetFigFont{6}{7.2}{\rmdefault}{\mddefault}{\updefault}$\alpha\beta$}}}}}
\put(900,4854){\makebox(0,0)[lb]{\smash{{{\SetFigFont{6}{7.2}{\rmdefault}{\mddefault}{\updefault}$\alpha$}}}}}
\put(2400,4929){\makebox(0,0)[lb]{\smash{{{\SetFigFont{6}{7.2}{\rmdefault}{\mddefault}{\updefault}$\beta$}}}}}
\put(300,5529){\makebox(0,0)[lb]{\smash{{{\SetFigFont{6}{7.2}{\rmdefault}{\mddefault}{\updefault}$1$}}}}}
\put(1800,5454){\makebox(0,0)[lb]{\smash{{{\SetFigFont{6}{7.2}{\rmdefault}{\mddefault}{\updefault}$1$}}}}}
\put(3525,6204){\makebox(0,0)[lb]{\smash{{{\SetFigFont{6}{7.2}{\rmdefault}{\mddefault}{\updefault}$s_{U,V}$}}}}}
\end{picture}
}

\end{center}
\vspace{0.3cm}

Condition (3) follows immediately from the naturality of $s$ and
condition (4) follows by pre and post composing $\theta$ with $\varphi
(\gamma) $, observing that the underlying geometry is the same 
in both cases, so the induced natural transformations agree. This
implies $\theta_{\varphi(\gamma)U} = \varphi(\gamma)(\theta_U)$.
 
This completes the proof of part (a) of Theorem \ref{prop:balanced}.

For part (b) let $X$ be any space. We claim that $\cA_1$ is a balanced
category with $\pi_2X$-action (see Appendix \ref{app:defns} for definitions).
Let $c\colon S_1 \rightarrow X$ be the collapse map and suppose this
is the choice made to represent $1\in \pi_1$ i.e. $\cA_1 = F(c)$. By
restricting to degenerate loops we obtain a monoidal product $*$,
monoidal unit $\monunit$ , braiding $s$ and twist $\theta$ in a
similar way to that in part (a) and $(\cA_1, *, \monunit, s, \theta)$
is a balanced category.  The only genuinely new structure in
part (b)  is the $\pi_2X$-action. Let $I$ be a cylinder and let $I
\rightarrow X$ be a map as indicated below

\vspace{0.3cm}
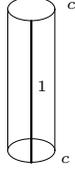
\begin{center}
\setlength{\unitlength}{0.00041667in}
\begingroup\makeatletter\ifx\SetFigFont\undefined%
\gdef\SetFigFont#1#2#3#4#5{%
  \reset@font\fontsize{#1}{#2pt}%
  \fontfamily{#3}\fontseries{#4}\fontshape{#5}%
  \selectfont}%
\fi\endgroup%
{\renewcommand{\dashlinestretch}{30}
\begin{picture}(854,2144)(0,-10)
\put(312,1972){\ellipse{600}{300}}
\put(312,172){\ellipse{600}{300}}
\path(12,172)(12,1972)
\path(612,1972)(612,172)
\thicklines
\path(312,1822)(312,22)
\put(687,22){\makebox(0,0)[lb]{\smash{{{\SetFigFont{6}{7.2}{\rmdefault}{\mddefault}{\updefault}$c$}}}}}
\put(387,922){\makebox(0,0)[lb]{\smash{{{\SetFigFont{6}{7.2}{\rmdefault}{\mddefault}{\updefault}1}}}}}
\put(762,1972){\makebox(0,0)[lb]{\smash{{{\SetFigFont{6}{7.2}{\rmdefault}{\mddefault}{\updefault}$c$}}}}}
\end{picture}
}

\end{center}
\vspace{0.3cm}

Such a map determines an element of $\pi_2(X)$ and any two maps giving
the same element are homotopic and hence 2-isomorphic. For each $g\in
\pi_2(X)$ choose $i_g\colon I\rightarrow X$ and define $\rho\colon
\pi_2(X) \rightarrow \oAut(\cA)$ by
\[
\rho(g) = F({i_g})
\]
The composition of morphisms ${i_g}$ corresponds to addition
in $\pi_2X$ and hence $\rho$ is a group homomorphism.
There is an RS-equivalence

\vspace{0.3cm}
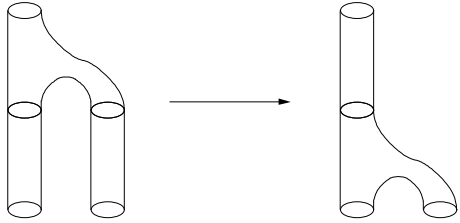
\begin{center}
\setlength{\unitlength}{0.00029167in}
\begingroup\makeatletter\ifx\SetFigFont\undefined%
\gdef\SetFigFont#1#2#3#4#5{%
  \reset@font\fontsize{#1}{#2pt}%
  \fontfamily{#3}\fontseries{#4}\fontshape{#5}%
  \selectfont}%
\fi\endgroup%
{\renewcommand{\dashlinestretch}{30}
\begin{picture}(8124,3929)(0,-10)
\put(6312,3757){\ellipse{600}{300}}
\put(6312,1957){\ellipse{600}{300}}
\path(6012,1957)(6012,3757)
\path(6612,3757)(6612,1957)
\put(312,1957){\ellipse{600}{300}}
\put(312,157){\ellipse{600}{300}}
\path(12,157)(12,1957)
\path(612,1957)(612,157)
\put(1812,1957){\ellipse{600}{300}}
\put(1812,157){\ellipse{600}{300}}
\path(1512,157)(1512,1957)
\path(2112,1957)(2112,157)
\put(7812,157){\ellipse{600}{300}}
\put(6312,1957){\ellipse{600}{300}}
\put(6312,157){\ellipse{600}{300}}
\path(6612,157)(6612,307)
\path(7512,157)(7512,307)
\path(6612,1957)(6612,1956)(6612,1953)
        (6612,1948)(6613,1940)(6613,1929)
        (6614,1915)(6616,1898)(6618,1879)
        (6621,1857)(6625,1833)(6631,1807)
        (6638,1780)(6646,1751)(6657,1721)
        (6669,1689)(6685,1656)(6704,1621)
        (6726,1583)(6753,1543)(6785,1500)
        (6822,1455)(6864,1407)(6912,1357)
        (6950,1320)(6990,1284)(7029,1250)
        (7067,1218)(7103,1190)(7137,1165)
        (7169,1143)(7198,1125)(7224,1109)
        (7249,1097)(7271,1086)(7291,1078)
        (7310,1071)(7328,1066)(7345,1061)
        (7362,1057)(7379,1053)(7396,1048)
        (7414,1043)(7433,1036)(7453,1028)
        (7475,1017)(7500,1005)(7526,989)
        (7555,971)(7587,949)(7621,924)
        (7657,896)(7695,864)(7734,830)
        (7774,794)(7812,757)(7860,707)
        (7902,659)(7939,614)(7971,571)
        (7998,531)(8020,493)(8039,458)
        (8055,425)(8067,393)(8078,363)
        (8086,334)(8093,307)(8099,281)
        (8103,257)(8106,235)(8108,216)
        (8110,199)(8111,185)(8111,174)
        (8112,166)(8112,161)(8112,158)(8112,157)
\path(6012,157)(6012,1957)
\path(6612,307)(6612,309)(6613,313)
        (6615,320)(6617,330)(6620,343)
        (6624,359)(6629,377)(6635,396)
        (6642,417)(6650,438)(6660,460)
        (6671,483)(6684,506)(6699,530)
        (6717,556)(6738,581)(6762,607)
        (6785,629)(6807,648)(6827,664)
        (6845,678)(6860,689)(6871,698)
        (6881,705)(6888,711)(6894,716)
        (6900,720)(6905,723)(6912,727)
        (6920,730)(6930,734)(6943,739)
        (6960,743)(6980,748)(7005,752)
        (7032,756)(7062,757)(7092,756)
        (7119,752)(7144,748)(7164,743)
        (7181,739)(7194,734)(7204,730)
        (7212,727)(7219,723)(7225,719)
        (7230,716)(7236,711)(7243,705)
        (7253,698)(7264,689)(7279,678)
        (7297,664)(7317,648)(7339,629)
        (7362,607)(7386,581)(7407,556)
        (7425,530)(7440,506)(7453,483)
        (7464,460)(7474,438)(7482,417)
        (7489,396)(7495,377)(7500,359)
        (7504,343)(7507,330)(7509,320)
        (7511,313)(7512,309)(7512,307)
\put(1812,1957){\ellipse{600}{300}}
\put(312,3757){\ellipse{600}{300}}
\put(312,1957){\ellipse{600}{300}}
\path(612,1957)(612,2107)
\path(1512,1957)(1512,2107)
\path(2937,2107)(5037,2107)
\blacken\thicklines
\path(4917.000,2077.000)(5037.000,2107.000)(4917.000,2137.000)(4917.000,2077.000)
\thinlines
\path(612,3757)(612,3756)(612,3753)
        (612,3748)(613,3740)(613,3729)
        (614,3715)(616,3698)(618,3679)
        (621,3657)(625,3633)(631,3607)
        (638,3580)(646,3551)(657,3521)
        (669,3489)(685,3456)(704,3421)
        (726,3383)(753,3343)(785,3300)
        (822,3255)(864,3207)(912,3157)
        (950,3120)(990,3084)(1029,3050)
        (1067,3018)(1103,2990)(1137,2965)
        (1169,2943)(1198,2925)(1224,2909)
        (1249,2897)(1271,2886)(1291,2878)
        (1310,2871)(1328,2866)(1345,2861)
        (1362,2857)(1379,2853)(1396,2848)
        (1414,2843)(1433,2836)(1453,2828)
        (1475,2817)(1500,2805)(1526,2789)
        (1555,2771)(1587,2749)(1621,2724)
        (1657,2696)(1695,2664)(1734,2630)
        (1774,2594)(1812,2557)(1860,2507)
        (1902,2459)(1939,2414)(1971,2371)
        (1998,2331)(2020,2293)(2039,2258)
        (2055,2225)(2067,2193)(2078,2163)
        (2086,2134)(2093,2107)(2099,2081)
        (2103,2057)(2106,2035)(2108,2016)
        (2110,1999)(2111,1985)(2111,1974)
        (2112,1966)(2112,1961)(2112,1958)(2112,1957)
\path(12,1957)(12,3757)
\path(612,2107)(612,2109)(613,2113)
        (615,2120)(617,2130)(620,2143)
        (624,2159)(629,2177)(635,2196)
        (642,2217)(650,2238)(660,2260)
        (671,2283)(684,2306)(699,2330)
        (717,2356)(738,2381)(762,2407)
        (785,2429)(807,2448)(827,2464)
        (845,2478)(860,2489)(871,2498)
        (881,2505)(888,2511)(894,2516)
        (900,2520)(905,2523)(912,2527)
        (920,2530)(930,2534)(943,2539)
        (960,2543)(980,2548)(1005,2552)
        (1032,2556)(1062,2557)(1092,2556)
        (1119,2552)(1144,2548)(1164,2543)
        (1181,2539)(1194,2534)(1204,2530)
        (1212,2527)(1219,2523)(1225,2519)
        (1230,2516)(1236,2511)(1243,2505)
        (1253,2498)(1264,2489)(1279,2478)
        (1297,2464)(1317,2448)(1339,2429)
        (1362,2407)(1386,2381)(1407,2356)
        (1425,2330)(1440,2306)(1453,2283)
        (1464,2260)(1474,2238)(1482,2217)
        (1489,2196)(1495,2177)(1500,2159)
        (1504,2143)(1507,2130)(1509,2120)
        (1511,2113)(1512,2109)(1512,2107)
\end{picture}
}

\end{center}
\vspace{0.3cm}

and so $F(p\circ (i_g\cup i_1)) = F(i_g\circ p)$ for $g\in H_2(X)$.
It follows that for $U,V\in\cA_1$
\[
(\rho(g)U)* V = \rho(g)(U*V)\;\;\;\mbox{ and } \rho (g)f*h = \rho (g) (f*h)
\]
and similarly 
\[
U*\rho (g)V = \rho (g)(U*V)\;\;\;\;\mbox{ and } 
 \rho (g) (f*h) = f*\rho (g)h.
\]

The proof of conditions (5) and (6) in definition \ref{defn:balancedG}
are similar to the proofs of \ref{defn:braid}(3) and
\ref{defn:twist}(4) in part (a) of Theorem \ref{prop:balanced}.
Conditions (3) and (4) can also be shown using similar arguments.

This completes the proof of part (b) of Theorem \ref{prop:balanced}.

To prove part (c) of Theorem \ref{prop:balanced} we prove $\cA$ is
semi-simple Artinian for the case $X=K(\pi,1)$. A similar argument
shows $\cA_1$ is semi-simple Artinian for a general space $X$. By the
work of Tillmann (see Appendix \ref{app:add}) we must produce a
non-degenerate form $\langle - , - \rangle \colon \cA_\alpha \otimes
\cA_{\alpha^{-1}} \rightarrow \hat{k}$ for each $\alpha \in\pi$.  Let
$C$ be a cylinder with two inputs and for $\alpha\in\pi$ choose a map
$c_\alpha\colon C\rightarrow X$ as indicated below.

\vspace{0.3cm}
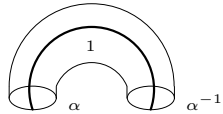
\begin{center}
\setlength{\unitlength}{0.00041667in}
\begingroup\makeatletter\ifx\SetFigFont\undefined%
\gdef\SetFigFont#1#2#3#4#5{%
  \reset@font\fontsize{#1}{#2pt}%
  \fontfamily{#3}\fontseries{#4}\fontshape{#5}%
  \selectfont}%
\fi\endgroup%
{\renewcommand{\dashlinestretch}{30}
\begin{picture}(2877,1428)(0,-10)
\put(1062.000,166.500){\arc{975.000}{3.5364}{5.8884}}
\put(1062.000,354.000){\arc{2100.000}{3.1416}{6.2832}}
\put(1812,204){\ellipse{600}{300}}
\put(312,204){\ellipse{600}{300}}
\path(612,204)(612,354)
\path(12,204)(12,354)
\path(1512,204)(1512,354)
\path(2112,204)(2112,354)
\thicklines
\put(1062.000,311.143){\arc{1585.714}{2.8113}{6.6135}}
\put(762,54){\makebox(0,0)[lb]{\smash{{{\SetFigFont{6}{7.2}{\rmdefault}{\mddefault}{\updefault}$\alpha$}}}}}
\put(2262,54){\makebox(0,0)[lb]{\smash{{{\SetFigFont{6}{7.2}{\rmdefault}{\mddefault}{\updefault}$\alpha^{-1}$}}}}}
\put(987,804){\makebox(0,0)[lb]{\smash{{{\SetFigFont{6}{7.2}{\rmdefault}{\mddefault}{\updefault}$1$}}}}}
\end{picture}
}

\end{center}
\vspace{0.3cm}

Now define
\begin{eqnarray*}
\langle -, - \rangle_\alpha:=F({c_\alpha}):\cA_\alpha \otimes
\cA_{\alpha^{-1}} \rightarrow \hat{k}&&
\end{eqnarray*}
and 
\[
\Delta_\alpha := F(\rotate {{c_\alpha}}) \colon
\hat{k}\rightarrow \cA_{\alpha} \otimes 
\cA_{\alpha^{-1}}
\]
and write the image of the canonical element $k\in \hat{k}$
as 
\[
E_{\alpha}: =\Delta_\alpha (k)=\sum_{i=1}^{n_\alpha} P_{\alpha}^i\otimes
Q_{\alpha^{-1}}^i
\]
where $P_{\alpha}^i \in \cA_{\alpha}$, $Q_{\alpha^{-1}}^i \in
\cA_{\alpha^{-1}}$. 

Standard topological quantum field theory manipulations show that $\langle -, - \rangle_\alpha$
and $E_{\alpha}$ provide a non-degenerate form and it follows from
Appendix \ref{app:add} that $\cA_\alpha$ is semi-simple Artinian proving
part (c) of  Theorem \ref{prop:balanced}.

\section{Tortile categories from self-dual $\cS_X$-structures} \label{section:tortile}

In this section
we consider under what conditions we can guarantee the category $\cA$
has (rigid) duals. In particular we
show that an $\cS_X$-structure that is lax self dual with respect to hom does have duals. Balanced
categories with duals are known as tortile categories (ribbon
categories) and we give appropriate variants of these in Appendix
\ref{app:defns}. Referring to these 
definitions  we prove the following theorem.

\begin{thm}\label{thm:main}
(a) Let $\pi$ be a discrete group and let $X$ be a based
  Eilenberg-Maclane space $K(\pi,1)$. The $k$-additive category $\cA$
  associated to an $\cS_X$-structure which is lax self dual with respect
  to hom,  is a semi-simple Artinian lax tortile $\pi$-category.  \\
(b) For any space $X$, the subcategory $\cA_1$ is a 
semi-simple Artinian tortile category with lax $\pi_2X$-action.
\end{thm}

Until further notice let $X=K(\pi,1)$ and let $F\colon \cS_X
\rightarrow \additive $ be an $\cS_X$-structure which is self dual with
respect to hom.  The key to getting duality in the category $\cA$ is
to use Tillmann's involutions $(-)^*\colon \cA_\alpha \rightarrow
\cA_{\alpha^{-1}}$ which arise from the non-degenerate forms $\langle
- , - \rangle_\alpha$ of the previous section. 
Non-degeneracy is courtesy of functors
\begin{eqnarray*}
&&\cI\colon \cA_{\alpha^{-1}} \rightarrow \cA_\alpha^\vee \;\;\;\;\;\;
Y\mapsto \langle Y,-\rangle\\ 
&&\cJ\colon\cA_\alpha^\vee\rightarrow \cA_{\alpha^{-1}}\;\;\;\;\;\;
H\mapsto \sum_{i=1}^nH(P^i_{\alpha})\otimes
Q^i_{\alpha^{-1}} 
\end{eqnarray*}
for which there are natural transformations $N\colon\mbox{id}\cong
\cI\cJ$ and $M\colon\mbox{id}\cong \cJ\cI$.

Notice that these
forms satisfy a ``Frobenius'' condition, namely, for $U\in \cA_\alpha, V\in
\cA_\beta$ and $W\in\cA_\gamma$ with $\alpha\beta\gamma = 1$ there are
natural isomorphisms
\[ 
\langle U, V*W \rangle_\alpha \simeq 
\langle U*V, W \rangle_{\gamma^{-1}}
\]
which arise from the diffeomorphism  indicated below.

\vspace{0.3cm}
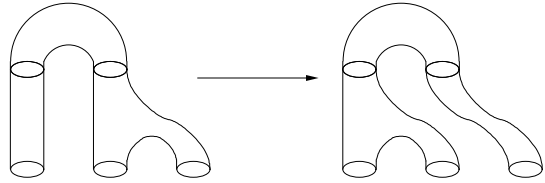
\begin{center}
\setlength{\unitlength}{0.00029167in}
\begingroup\makeatletter\ifx\SetFigFont\undefined%
\gdef\SetFigFont#1#2#3#4#5{%
  \reset@font\fontsize{#1}{#2pt}%
  \fontfamily{#3}\fontseries{#4}\fontshape{#5}%
  \selectfont}%
\fi\endgroup%
{\renewcommand{\dashlinestretch}{30}
\begin{picture}(9624,3181)(0,-10)
\put(3312,157){\ellipse{600}{300}}
\put(1812,1957){\ellipse{600}{300}}
\put(1812,157){\ellipse{600}{300}}
\path(2112,157)(2112,307)
\path(3012,157)(3012,307)
\path(2112,1957)(2112,1956)(2112,1953)
        (2112,1948)(2113,1940)(2113,1929)
        (2114,1915)(2116,1898)(2118,1879)
        (2121,1857)(2125,1833)(2131,1807)
        (2138,1780)(2146,1751)(2157,1721)
        (2169,1689)(2185,1656)(2204,1621)
        (2226,1583)(2253,1543)(2285,1500)
        (2322,1455)(2364,1407)(2412,1357)
        (2450,1320)(2490,1284)(2529,1250)
        (2567,1218)(2603,1190)(2637,1165)
        (2669,1143)(2698,1125)(2724,1109)
        (2749,1097)(2771,1086)(2791,1078)
        (2810,1071)(2828,1066)(2845,1061)
        (2862,1057)(2879,1053)(2896,1048)
        (2914,1043)(2933,1036)(2953,1028)
        (2975,1017)(3000,1005)(3026,989)
        (3055,971)(3087,949)(3121,924)
        (3157,896)(3195,864)(3234,830)
        (3274,794)(3312,757)(3360,707)
        (3402,659)(3439,614)(3471,571)
        (3498,531)(3520,493)(3539,458)
        (3555,425)(3567,393)(3578,363)
        (3586,334)(3593,307)(3599,281)
        (3603,257)(3606,235)(3608,216)
        (3610,199)(3611,185)(3611,174)
        (3612,166)(3612,161)(3612,158)(3612,157)
\path(1512,157)(1512,1957)
\path(2112,307)(2112,309)(2113,313)
        (2115,320)(2117,330)(2120,343)
        (2124,359)(2129,377)(2135,396)
        (2142,417)(2150,438)(2160,460)
        (2171,483)(2184,506)(2199,530)
        (2217,556)(2238,581)(2262,607)
        (2285,629)(2307,648)(2327,664)
        (2345,678)(2360,689)(2371,698)
        (2381,705)(2388,711)(2394,716)
        (2400,720)(2405,723)(2412,727)
        (2420,730)(2430,734)(2443,739)
        (2460,743)(2480,748)(2505,752)
        (2532,756)(2562,757)(2592,756)
        (2619,752)(2644,748)(2664,743)
        (2681,739)(2694,734)(2704,730)
        (2712,727)(2719,723)(2725,719)
        (2730,716)(2736,711)(2743,705)
        (2753,698)(2764,689)(2779,678)
        (2797,664)(2817,648)(2839,629)
        (2862,607)(2886,581)(2907,556)
        (2925,530)(2940,506)(2953,483)
        (2964,460)(2974,438)(2982,417)
        (2989,396)(2995,377)(3000,359)
        (3004,343)(3007,330)(3009,320)
        (3011,313)(3012,309)(3012,307)
\put(312,1957){\ellipse{600}{300}}
\put(312,157){\ellipse{600}{300}}
\path(12,157)(12,1957)
\path(612,1957)(612,157)
\put(1062.000,1919.500){\arc{975.000}{3.5364}{5.8884}}
\put(1062.000,2107.000){\arc{2100.000}{3.1416}{6.2832}}
\put(1812,1957){\ellipse{600}{300}}
\put(312,1957){\ellipse{600}{300}}
\path(612,1957)(612,2107)
\path(12,1957)(12,2107)
\path(1512,1957)(1512,2107)
\path(2112,1957)(2112,2107)
\put(7062.000,1919.500){\arc{975.000}{3.5364}{5.8884}}
\put(7062.000,2107.000){\arc{2100.000}{3.1416}{6.2832}}
\put(7812,1957){\ellipse{600}{300}}
\put(6312,1957){\ellipse{600}{300}}
\path(6612,1957)(6612,2107)
\path(6012,1957)(6012,2107)
\path(7512,1957)(7512,2107)
\path(8112,1957)(8112,2107)
\put(7812,157){\ellipse{600}{300}}
\put(6312,1957){\ellipse{600}{300}}
\put(6312,157){\ellipse{600}{300}}
\put(9312,157){\ellipse{600}{300}}
\put(7812,1957){\ellipse{600}{300}}
\path(6612,157)(6612,307)
\path(7512,157)(7512,307)
\path(6612,1957)(6612,1956)(6612,1953)
        (6612,1948)(6613,1940)(6613,1929)
        (6614,1915)(6616,1898)(6618,1879)
        (6621,1857)(6625,1833)(6631,1807)
        (6638,1780)(6646,1751)(6657,1721)
        (6669,1689)(6685,1656)(6704,1621)
        (6726,1583)(6753,1543)(6785,1500)
        (6822,1455)(6864,1407)(6912,1357)
        (6950,1320)(6990,1284)(7029,1250)
        (7067,1218)(7103,1190)(7137,1165)
        (7169,1143)(7198,1125)(7224,1109)
        (7249,1097)(7271,1086)(7291,1078)
        (7310,1071)(7328,1066)(7345,1061)
        (7362,1057)(7379,1053)(7396,1048)
        (7414,1043)(7433,1036)(7453,1028)
        (7475,1017)(7500,1005)(7526,989)
        (7555,971)(7587,949)(7621,924)
        (7657,896)(7695,864)(7734,830)
        (7774,794)(7812,757)(7860,707)
        (7902,659)(7939,614)(7971,571)
        (7998,531)(8020,493)(8039,458)
        (8055,425)(8067,393)(8078,363)
        (8086,334)(8093,307)(8099,281)
        (8103,257)(8106,235)(8108,216)
        (8110,199)(8111,185)(8111,174)
        (8112,166)(8112,161)(8112,158)(8112,157)
\path(6012,157)(6012,1957)
\path(6612,307)(6612,309)(6613,313)
        (6615,320)(6617,330)(6620,343)
        (6624,359)(6629,377)(6635,396)
        (6642,417)(6650,438)(6660,460)
        (6671,483)(6684,506)(6699,530)
        (6717,556)(6738,581)(6762,607)
        (6785,629)(6807,648)(6827,664)
        (6845,678)(6860,689)(6871,698)
        (6881,705)(6888,711)(6894,716)
        (6900,720)(6905,723)(6912,727)
        (6920,730)(6930,734)(6943,739)
        (6960,743)(6980,748)(7005,752)
        (7032,756)(7062,757)(7092,756)
        (7119,752)(7144,748)(7164,743)
        (7181,739)(7194,734)(7204,730)
        (7212,727)(7219,723)(7225,719)
        (7230,716)(7236,711)(7243,705)
        (7253,698)(7264,689)(7279,678)
        (7297,664)(7317,648)(7339,629)
        (7362,607)(7386,581)(7407,556)
        (7425,530)(7440,506)(7453,483)
        (7464,460)(7474,438)(7482,417)
        (7489,396)(7495,377)(7500,359)
        (7504,343)(7507,330)(7509,320)
        (7511,313)(7512,309)(7512,307)
\path(8112,1957)(8112,1956)(8112,1953)
        (8112,1948)(8113,1940)(8113,1929)
        (8114,1915)(8116,1898)(8118,1879)
        (8121,1857)(8125,1833)(8131,1807)
        (8138,1780)(8146,1751)(8157,1721)
        (8169,1689)(8185,1656)(8204,1621)
        (8226,1583)(8253,1543)(8285,1500)
        (8322,1455)(8364,1407)(8412,1357)
        (8450,1320)(8490,1284)(8529,1250)
        (8567,1218)(8603,1190)(8637,1165)
        (8669,1143)(8698,1125)(8724,1109)
        (8749,1097)(8771,1086)(8791,1078)
        (8810,1071)(8828,1066)(8845,1061)
        (8862,1057)(8879,1053)(8896,1048)
        (8914,1043)(8933,1036)(8953,1028)
        (8975,1017)(9000,1005)(9026,989)
        (9055,971)(9087,949)(9121,924)
        (9157,896)(9195,864)(9234,830)
        (9274,794)(9312,757)(9360,707)
        (9402,659)(9439,614)(9471,571)
        (9498,531)(9520,493)(9539,458)
        (9555,425)(9567,393)(9578,363)
        (9586,334)(9593,307)(9599,281)
        (9603,257)(9606,235)(9608,216)
        (9610,199)(9611,185)(9611,174)
        (9612,166)(9612,161)(9612,158)(9612,157)
\path(7512,1957)(7512,1956)(7512,1953)
        (7512,1948)(7513,1940)(7513,1929)
        (7514,1915)(7516,1898)(7518,1879)
        (7521,1857)(7525,1833)(7531,1807)
        (7538,1780)(7546,1751)(7557,1721)
        (7569,1689)(7585,1656)(7604,1621)
        (7626,1583)(7653,1543)(7685,1500)
        (7722,1455)(7764,1407)(7812,1357)
        (7850,1320)(7890,1284)(7929,1250)
        (7967,1218)(8003,1190)(8037,1165)
        (8069,1143)(8098,1125)(8124,1109)
        (8149,1097)(8171,1086)(8191,1078)
        (8210,1071)(8228,1066)(8245,1061)
        (8262,1057)(8279,1053)(8296,1048)
        (8314,1043)(8333,1036)(8353,1028)
        (8375,1017)(8400,1005)(8426,989)
        (8455,971)(8487,949)(8521,924)
        (8557,896)(8595,864)(8634,830)
        (8674,794)(8712,757)(8760,707)
        (8802,659)(8839,614)(8871,571)
        (8898,531)(8920,493)(8939,458)
        (8955,425)(8967,393)(8978,363)
        (8986,334)(8993,307)(8999,281)
        (9003,257)(9006,235)(9008,216)
        (9010,199)(9011,185)(9011,174)
        (9012,166)(9012,161)(9012,158)(9012,157)
\path(3387,1807)(5487,1807)
\blacken\thicklines
\path(5367.000,1777.000)(5487.000,1807.000)(5367.000,1837.000)(5367.000,1777.000)
\end{picture}
}

\end{center}
\vspace{0.3cm}

To define Tillmann's involutions let $\hom\colon \cA_\alpha \rightarrow
\cA_\alpha^\vee$ be  defined by $Y
\rightarrow \cA_\alpha(Y,-)$ and set $(-)^* = \cJ \circ \hom$.
 These functors satisfy $(-)^{**} \cong \mbox{id}$ and $\cA_\alpha(-,-) \simeq \langle -^*,-
\rangle_{\alpha^{-1}}$ and both the forms and
involutions extend to tensor products. In particular we have we have
non-degenerate forms and involutions
\[
\langle - , - \rangle_\aub \colon \cA_\aub \otimes \cA_\aubi
\rightarrow \hat{k}
\qquad \qquad
 (-)^*_\aub \colon \cA_\aub \rightarrow \cA_\aubi
\]
where $\aub = (\alpha_1, \ldots , \alpha_n)$, $\aubi = (\alpha^{-1}_n,
\ldots , \alpha_1^{-1})$ and $\cA_\aub = \cA_{\alpha_1} \otimes \cdots
\otimes \cA_{\alpha_n}$.  We will write $\cI$ and $\cJ$ for the
associated maps $\cA_\aub \rightarrow \cA_\aubi^\vee$ and $
\cA_\aubi^\vee \rightarrow \cA_\aub$.  See Appendix \ref{app:add} for
further discussion.

The functors $\cJ$ have the following property.

\begin{lem}\label{lem:jumpJ}
Let $g\colon \Sigma
\rightarrow X$ be a morphism in $\cS_X$. Then there is a natural isomorphism
\[
t_{F(g)}: \cJ\circ F(\rotate{g})^\vee
  \nattrans F(g)\circ \cJ
\]
\end{lem}
\begin{proof}
By using the diffeomorphism below we have natural isomorphisms 
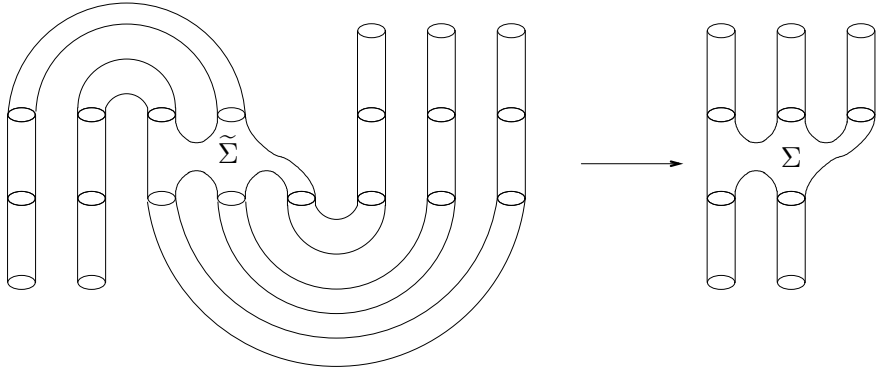
\[
k\colon \cJ\circ F(\rotate{g})^\vee \circ \cI \nattrans F(g).
\]

\vspace{0.3cm}
\begin{center}
\setlength{\unitlength}{0.00025000in}
\begingroup\makeatletter\ifx\SetFigFont\undefined%
\gdef\SetFigFont#1#2#3#4#5{%
  \reset@font\fontsize{#1}{#2pt}%
  \fontfamily{#3}\fontseries{#4}\fontshape{#5}%
  \selectfont}%
\fi\endgroup%
{\renewcommand{\dashlinestretch}{30}
\begin{picture}(18624,7834)(0,-10)
\put(4812,3609){\ellipse{600}{300}}
\put(3312,5409){\ellipse{600}{300}}
\put(3312,3609){\ellipse{600}{300}}
\put(6312,3609){\ellipse{600}{300}}
\put(4812,5409){\ellipse{600}{300}}
\path(3612,3609)(3612,3759)
\path(4512,3609)(4512,3759)
\path(5112,3609)(5112,3759)
\path(6012,3609)(6012,3759)
\path(3612,5409)(3612,5259)
\path(4512,5409)(4512,5259)
\path(3012,3609)(3012,5409)
\path(3612,3759)(3612,3761)(3613,3765)
        (3615,3772)(3617,3782)(3620,3795)
        (3624,3811)(3629,3829)(3635,3848)
        (3642,3869)(3650,3890)(3660,3912)
        (3671,3935)(3684,3958)(3699,3982)
        (3717,4008)(3738,4033)(3762,4059)
        (3785,4081)(3807,4100)(3827,4116)
        (3845,4130)(3860,4141)(3871,4150)
        (3881,4157)(3888,4163)(3894,4168)
        (3900,4172)(3905,4175)(3912,4179)
        (3920,4182)(3930,4186)(3943,4191)
        (3960,4195)(3980,4200)(4005,4204)
        (4032,4208)(4062,4209)(4092,4208)
        (4119,4204)(4144,4200)(4164,4195)
        (4181,4191)(4194,4186)(4204,4182)
        (4212,4179)(4219,4175)(4225,4171)
        (4230,4168)(4236,4163)(4243,4157)
        (4253,4150)(4264,4141)(4279,4130)
        (4297,4116)(4317,4100)(4339,4081)
        (4362,4059)(4386,4033)(4407,4008)
        (4425,3982)(4440,3958)(4453,3935)
        (4464,3912)(4474,3890)(4482,3869)
        (4489,3848)(4495,3829)(4500,3811)
        (4504,3795)(4507,3782)(4509,3772)
        (4511,3765)(4512,3761)(4512,3759)
\path(5112,5409)(5112,5408)(5112,5405)
        (5112,5400)(5113,5392)(5113,5381)
        (5114,5367)(5116,5350)(5118,5331)
        (5121,5309)(5125,5285)(5131,5259)
        (5138,5232)(5146,5203)(5157,5173)
        (5169,5141)(5185,5108)(5204,5073)
        (5226,5035)(5253,4995)(5285,4952)
        (5322,4907)(5364,4859)(5412,4809)
        (5450,4772)(5490,4736)(5529,4702)
        (5567,4670)(5603,4642)(5637,4617)
        (5669,4595)(5698,4577)(5724,4561)
        (5749,4549)(5771,4538)(5791,4530)
        (5810,4523)(5828,4518)(5845,4513)
        (5862,4509)(5879,4505)(5896,4500)
        (5914,4495)(5933,4488)(5953,4480)
        (5975,4469)(6000,4457)(6026,4441)
        (6055,4423)(6087,4401)(6121,4376)
        (6157,4348)(6195,4316)(6234,4282)
        (6274,4246)(6312,4209)(6360,4159)
        (6402,4111)(6439,4066)(6471,4023)
        (6498,3983)(6520,3945)(6539,3910)
        (6555,3877)(6567,3845)(6578,3815)
        (6586,3786)(6593,3759)(6599,3733)
        (6603,3709)(6606,3687)(6608,3668)
        (6610,3651)(6611,3637)(6611,3626)
        (6612,3618)(6612,3613)(6612,3610)(6612,3609)
\path(5112,3759)(5112,3761)(5113,3765)
        (5115,3772)(5117,3782)(5120,3795)
        (5124,3811)(5129,3829)(5135,3848)
        (5142,3869)(5150,3890)(5160,3912)
        (5171,3935)(5184,3958)(5199,3982)
        (5217,4008)(5238,4033)(5262,4059)
        (5285,4081)(5307,4100)(5327,4116)
        (5345,4130)(5360,4141)(5371,4150)
        (5381,4157)(5388,4163)(5394,4168)
        (5400,4172)(5405,4175)(5412,4179)
        (5420,4182)(5430,4186)(5443,4191)
        (5460,4195)(5480,4200)(5505,4204)
        (5532,4208)(5562,4209)(5592,4208)
        (5619,4204)(5644,4200)(5664,4195)
        (5681,4191)(5694,4186)(5704,4182)
        (5712,4179)(5719,4175)(5725,4171)
        (5730,4168)(5736,4163)(5743,4157)
        (5753,4150)(5764,4141)(5779,4130)
        (5797,4116)(5817,4100)(5839,4081)
        (5862,4059)(5886,4033)(5907,4008)
        (5925,3982)(5940,3958)(5953,3935)
        (5964,3912)(5974,3890)(5982,3869)
        (5989,3848)(5995,3829)(6000,3811)
        (6004,3795)(6007,3782)(6009,3772)
        (6011,3765)(6012,3761)(6012,3759)
\path(3612,5259)(3612,5257)(3613,5253)
        (3615,5246)(3617,5236)(3620,5223)
        (3624,5207)(3629,5189)(3635,5170)
        (3642,5149)(3650,5128)(3660,5106)
        (3671,5083)(3684,5060)(3699,5036)
        (3717,5010)(3738,4985)(3762,4959)
        (3785,4937)(3807,4918)(3827,4902)
        (3845,4888)(3860,4877)(3871,4868)
        (3881,4861)(3888,4855)(3894,4850)
        (3900,4846)(3905,4843)(3912,4839)
        (3920,4836)(3930,4832)(3943,4827)
        (3960,4823)(3980,4818)(4005,4814)
        (4032,4810)(4062,4809)(4092,4810)
        (4119,4814)(4144,4818)(4164,4823)
        (4181,4827)(4194,4832)(4204,4836)
        (4212,4839)(4219,4843)(4225,4847)
        (4230,4850)(4236,4855)(4243,4861)
        (4253,4868)(4264,4877)(4279,4888)
        (4297,4902)(4317,4918)(4339,4937)
        (4362,4959)(4386,4985)(4407,5010)
        (4425,5036)(4440,5060)(4453,5083)
        (4464,5106)(4474,5128)(4482,5149)
        (4489,5170)(4495,5189)(4500,5207)
        (4504,5223)(4507,5236)(4509,5246)
        (4511,5253)(4512,5257)(4512,5259)
\put(7062.000,3646.500){\arc{975.000}{0.3948}{2.7468}}
\put(7062.000,3459.000){\arc{2100.000}{6.2832}{9.4248}}
\put(7812,3609){\ellipse{600}{300}}
\put(6312,3609){\ellipse{600}{300}}
\path(6612,3609)(6612,3459)
\path(6012,3609)(6012,3459)
\path(7512,3609)(7512,3459)
\path(8112,3609)(8112,3459)
\put(7812,5409){\ellipse{600}{300}}
\put(7812,3609){\ellipse{600}{300}}
\path(7512,3609)(7512,5409)
\path(8112,5409)(8112,3609)
\put(2562.000,5371.500){\arc{975.000}{3.5364}{5.8884}}
\put(2562.000,5559.000){\arc{2100.000}{3.1416}{6.2832}}
\put(3312,5409){\ellipse{600}{300}}
\put(1812,5409){\ellipse{600}{300}}
\path(2112,5409)(2112,5559)
\path(1512,5409)(1512,5559)
\path(3012,5409)(3012,5559)
\path(3612,5409)(3612,5559)
\put(9312,5409){\ellipse{600}{300}}
\put(9312,3609){\ellipse{600}{300}}
\path(9012,3609)(9012,5409)
\path(9612,5409)(9612,3609)
\put(10812,5409){\ellipse{600}{300}}
\put(10812,3609){\ellipse{600}{300}}
\path(10512,3609)(10512,5409)
\path(11112,5409)(11112,3609)
\put(1812,5409){\ellipse{600}{300}}
\put(1812,3609){\ellipse{600}{300}}
\path(1512,3609)(1512,5409)
\path(2112,5409)(2112,3609)
\put(312,5409){\ellipse{600}{300}}
\put(312,3609){\ellipse{600}{300}}
\path(12,3609)(12,5409)
\path(612,5409)(612,3609)
\put(16812,5409){\ellipse{600}{300}}
\put(15312,3609){\ellipse{600}{300}}
\put(15312,5409){\ellipse{600}{300}}
\put(18312,5409){\ellipse{600}{300}}
\put(16812,3609){\ellipse{600}{300}}
\path(15612,5409)(15612,5259)
\path(16512,5409)(16512,5259)
\path(17112,5409)(17112,5259)
\path(18012,5409)(18012,5259)
\path(15612,3609)(15612,3759)
\path(16512,3609)(16512,3759)
\path(15012,5409)(15012,3609)
\path(15612,5259)(15612,5257)(15613,5253)
        (15615,5246)(15617,5236)(15620,5223)
        (15624,5207)(15629,5189)(15635,5170)
        (15642,5149)(15650,5128)(15660,5106)
        (15671,5083)(15684,5060)(15699,5036)
        (15717,5010)(15738,4985)(15762,4959)
        (15785,4937)(15807,4918)(15827,4902)
        (15845,4888)(15860,4877)(15871,4868)
        (15881,4861)(15888,4855)(15894,4850)
        (15900,4846)(15905,4843)(15912,4839)
        (15920,4836)(15930,4832)(15943,4827)
        (15960,4823)(15980,4818)(16005,4814)
        (16032,4810)(16062,4809)(16092,4810)
        (16119,4814)(16144,4818)(16164,4823)
        (16181,4827)(16194,4832)(16204,4836)
        (16212,4839)(16219,4843)(16225,4847)
        (16230,4850)(16236,4855)(16243,4861)
        (16253,4868)(16264,4877)(16279,4888)
        (16297,4902)(16317,4918)(16339,4937)
        (16362,4959)(16386,4985)(16407,5010)
        (16425,5036)(16440,5060)(16453,5083)
        (16464,5106)(16474,5128)(16482,5149)
        (16489,5170)(16495,5189)(16500,5207)
        (16504,5223)(16507,5236)(16509,5246)
        (16511,5253)(16512,5257)(16512,5259)
\path(17112,3609)(17112,3610)(17112,3613)
        (17112,3618)(17113,3626)(17113,3637)
        (17114,3651)(17116,3668)(17118,3687)
        (17121,3709)(17125,3733)(17131,3759)
        (17138,3786)(17146,3815)(17157,3845)
        (17169,3877)(17185,3910)(17204,3945)
        (17226,3983)(17253,4023)(17285,4066)
        (17322,4111)(17364,4159)(17412,4209)
        (17450,4246)(17490,4282)(17529,4316)
        (17567,4348)(17603,4376)(17637,4401)
        (17669,4423)(17698,4441)(17724,4457)
        (17749,4469)(17771,4480)(17791,4488)
        (17810,4495)(17828,4500)(17845,4505)
        (17862,4509)(17879,4513)(17896,4518)
        (17914,4523)(17933,4530)(17953,4538)
        (17975,4549)(18000,4561)(18026,4577)
        (18055,4595)(18087,4617)(18121,4642)
        (18157,4670)(18195,4702)(18234,4736)
        (18274,4772)(18312,4809)(18360,4859)
        (18402,4907)(18439,4952)(18471,4995)
        (18498,5035)(18520,5073)(18539,5108)
        (18555,5141)(18567,5173)(18578,5203)
        (18586,5232)(18593,5259)(18599,5285)
        (18603,5309)(18606,5331)(18608,5350)
        (18610,5367)(18611,5381)(18611,5392)
        (18612,5400)(18612,5405)(18612,5408)(18612,5409)
\path(17112,5259)(17112,5257)(17113,5253)
        (17115,5246)(17117,5236)(17120,5223)
        (17124,5207)(17129,5189)(17135,5170)
        (17142,5149)(17150,5128)(17160,5106)
        (17171,5083)(17184,5060)(17199,5036)
        (17217,5010)(17238,4985)(17262,4959)
        (17285,4937)(17307,4918)(17327,4902)
        (17345,4888)(17360,4877)(17371,4868)
        (17381,4861)(17388,4855)(17394,4850)
        (17400,4846)(17405,4843)(17412,4839)
        (17420,4836)(17430,4832)(17443,4827)
        (17460,4823)(17480,4818)(17505,4814)
        (17532,4810)(17562,4809)(17592,4810)
        (17619,4814)(17644,4818)(17664,4823)
        (17681,4827)(17694,4832)(17704,4836)
        (17712,4839)(17719,4843)(17725,4847)
        (17730,4850)(17736,4855)(17743,4861)
        (17753,4868)(17764,4877)(17779,4888)
        (17797,4902)(17817,4918)(17839,4937)
        (17862,4959)(17886,4985)(17907,5010)
        (17925,5036)(17940,5060)(17953,5083)
        (17964,5106)(17974,5128)(17982,5149)
        (17989,5170)(17995,5189)(18000,5207)
        (18004,5223)(18007,5236)(18009,5246)
        (18011,5253)(18012,5257)(18012,5259)
\path(15612,3759)(15612,3761)(15613,3765)
        (15615,3772)(15617,3782)(15620,3795)
        (15624,3811)(15629,3829)(15635,3848)
        (15642,3869)(15650,3890)(15660,3912)
        (15671,3935)(15684,3958)(15699,3982)
        (15717,4008)(15738,4033)(15762,4059)
        (15785,4081)(15807,4100)(15827,4116)
        (15845,4130)(15860,4141)(15871,4150)
        (15881,4157)(15888,4163)(15894,4168)
        (15900,4172)(15905,4175)(15912,4179)
        (15920,4182)(15930,4186)(15943,4191)
        (15960,4195)(15980,4200)(16005,4204)
        (16032,4208)(16062,4209)(16092,4208)
        (16119,4204)(16144,4200)(16164,4195)
        (16181,4191)(16194,4186)(16204,4182)
        (16212,4179)(16219,4175)(16225,4171)
        (16230,4168)(16236,4163)(16243,4157)
        (16253,4150)(16264,4141)(16279,4130)
        (16297,4116)(16317,4100)(16339,4081)
        (16362,4059)(16386,4033)(16407,4008)
        (16425,3982)(16440,3958)(16453,3935)
        (16464,3912)(16474,3890)(16482,3869)
        (16489,3848)(16495,3829)(16500,3811)
        (16504,3795)(16507,3782)(16509,3772)
        (16511,3765)(16512,3761)(16512,3759)
\put(7812,7209){\ellipse{600}{300}}
\put(7812,5409){\ellipse{600}{300}}
\path(7512,5409)(7512,7209)
\path(8112,7209)(8112,5409)
\put(9312,7209){\ellipse{600}{300}}
\put(9312,5409){\ellipse{600}{300}}
\path(9012,5409)(9012,7209)
\path(9612,7209)(9612,5409)
\put(10812,7209){\ellipse{600}{300}}
\put(10812,5409){\ellipse{600}{300}}
\path(10512,5409)(10512,7209)
\path(11112,7209)(11112,5409)
\put(1812,3609){\ellipse{600}{300}}
\put(1812,1809){\ellipse{600}{300}}
\path(1512,1809)(1512,3609)
\path(2112,3609)(2112,1809)
\put(312,3609){\ellipse{600}{300}}
\put(312,1809){\ellipse{600}{300}}
\path(12,1809)(12,3609)
\path(612,3609)(612,1809)
\put(15312,7209){\ellipse{600}{300}}
\put(15312,5409){\ellipse{600}{300}}
\path(15012,5409)(15012,7209)
\path(15612,7209)(15612,5409)
\put(16812,7209){\ellipse{600}{300}}
\put(16812,5409){\ellipse{600}{300}}
\path(16512,5409)(16512,7209)
\path(17112,7209)(17112,5409)
\put(18312,7209){\ellipse{600}{300}}
\put(18312,5409){\ellipse{600}{300}}
\path(18012,5409)(18012,7209)
\path(18612,7209)(18612,5409)
\put(15312,3609){\ellipse{600}{300}}
\put(15312,1809){\ellipse{600}{300}}
\path(15012,1809)(15012,3609)
\path(15612,3609)(15612,1809)
\put(16812,3609){\ellipse{600}{300}}
\put(16812,1809){\ellipse{600}{300}}
\path(16512,1809)(16512,3609)
\path(17112,3609)(17112,1809)
\put(7061.265,3609.735){\arc{3901.471}{0.0004}{3.1028}}
\put(7062.000,4092.750){\arc{6967.500}{0.1393}{3.0023}}
\put(7057.184,4091.662){\arc{8166.884}{0.1185}{3.0046}}
\put(7060.303,3686.885){\arc{5105.770}{0.0305}{3.0817}}
\put(2560.303,5256.115){\arc{5105.771}{3.2015}{6.2527}}
\put(2561.265,5408.265){\arc{3901.470}{3.1804}{6.2828}}
\put(9312,3609){\ellipse{600}{300}}
\put(10812,3609){\ellipse{600}{300}}
\put(312,5409){\ellipse{600}{300}}
\path(12312,4359)(14412,4359)
\thicklines
\path(14232.000,4314.000)(14412.000,4359.000)(14232.000,4404.000)
\put(16587,4284){\makebox(0,0)[lb]{\smash{{{\SetFigFont{11}{13.2}{\rmdefault}{\mddefault}{\updefault}$\Sigma$}}}}}
\put(4512,4359){\makebox(0,0)[lb]{\smash{{{\SetFigFont{11}{13.2}{\rmdefault}{\mddefault}{\updefault}$\widetilde{\Sigma}$}}}}}
\end{picture}
}

\end{center}
\vspace{0.3cm}

Pre composing with $\cJ$ and using the equivalence $N:id\simeq \cI\cJ$
gives the required natural isomorphisms.
\end{proof}

For the crossing $\varphi$ we have the following additional
property. 

\begin{lem} \label{lem:crossJ}
The following diagram commutes.
\[
\xymatrix{ \cJ \circ \varphi (\alpha^{-1})^\vee
  \ar[rr]^{t_{\varphi(\alpha)}} & & \varphi (\alpha)\circ \cJ
  \ar[dl]^{\theta_{\varphi(\alpha)\cJ}} \\ & \cJ \ar[ul]^{J\theta^\vee}
}
\]
\end{lem}

\begin{proof}
Consider the
following diagram, where $H\in \cA^\vee$.
\[
\xymatrix{
J\circ \vai ^\vee (H)  \ar[rr]^-{J\vai ^\vee (N_{
    H})} && \cJ\circ\vai^ \vee \circ \cI \circ \cJ  (H)
\ar[rr]^-{k_{\va,\cJ(H)}} && \va \circ \cJ  (H)
\ar[d]^{\theta_{\va\cJ(H)}}\\ 
\cJ (H) \ar[u]^{\cJ\theta^\vee_{H}}
\ar[rr]_{\cJ(N_{H})} && \cJ \circ \cI \circ \cJ  (H)
\ar[u]^{\cJ\theta^\vee_{\cI\cJ(H)}}\ar[rr]_{M_{\cJ(H)}} &&
\cJ(H)} 
\]  

The left square commutes by naturality of
$\theta^\vee$ and the right square by the mapping class group
identity given below.

\vspace{0.3cm}
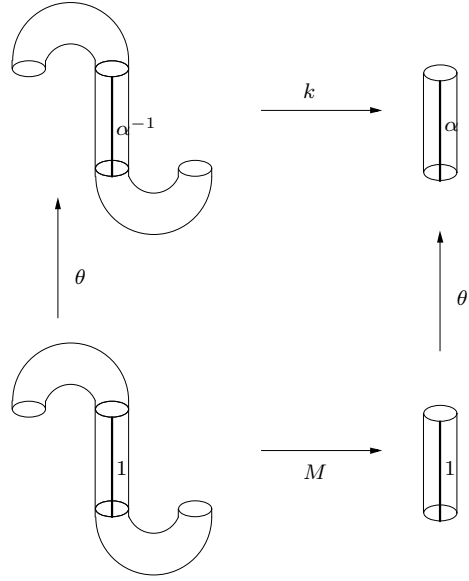
\begin{center}
\setlength{\unitlength}{0.00029167in}
\begingroup\makeatletter\ifx\SetFigFont\undefined%
\gdef\SetFigFont#1#2#3#4#5{%
  \reset@font\fontsize{#1}{#2pt}%
  \fontfamily{#3}\fontseries{#4}\fontshape{#5}%
  \selectfont}%
\fi\endgroup%
{\renewcommand{\dashlinestretch}{30}
\begin{picture}(8775,10383)(0,-10)
\put(7737,9084){\ellipse{600}{300}}
\put(7737,7284){\ellipse{600}{300}}
\path(7437,7284)(7437,9084)
\path(8037,9084)(8037,7284)
\put(1812,9159){\ellipse{600}{300}}
\put(1812,7359){\ellipse{600}{300}}
\path(1512,7359)(1512,9159)
\path(2112,9159)(2112,7359)
\put(2562.000,7396.500){\arc{975.000}{0.3948}{2.7468}}
\put(2562.000,7209.000){\arc{2100.000}{6.2832}{9.4248}}
\put(3312,7359){\ellipse{600}{300}}
\put(1812,7359){\ellipse{600}{300}}
\path(2112,7359)(2112,7209)
\path(1512,7359)(1512,7209)
\path(3012,7359)(3012,7209)
\path(3612,7359)(3612,7209)
\put(1062.000,9121.500){\arc{975.000}{3.5364}{5.8884}}
\put(1062.000,9309.000){\arc{2100.000}{3.1416}{6.2832}}
\put(1812,9159){\ellipse{600}{300}}
\put(312,9159){\ellipse{600}{300}}
\path(612,9159)(612,9309)
\path(12,9159)(12,9309)
\path(1512,9159)(1512,9309)
\path(2112,9159)(2112,9309)
\thicklines
\path(1812,7209)(1812,9009)
\thinlines
\put(7737,2934){\ellipse{600}{300}}
\put(7737,1134){\ellipse{600}{300}}
\path(7437,1134)(7437,2934)
\path(8037,2934)(8037,1134)
\put(1812,3009){\ellipse{600}{300}}
\put(1812,1209){\ellipse{600}{300}}
\path(1512,1209)(1512,3009)
\path(2112,3009)(2112,1209)
\put(2562.000,1246.500){\arc{975.000}{0.3948}{2.7468}}
\put(2562.000,1059.000){\arc{2100.000}{6.2832}{9.4248}}
\put(3312,1209){\ellipse{600}{300}}
\put(1812,1209){\ellipse{600}{300}}
\path(2112,1209)(2112,1059)
\path(1512,1209)(1512,1059)
\path(3012,1209)(3012,1059)
\path(3612,1209)(3612,1059)
\put(1062.000,2971.500){\arc{975.000}{3.5364}{5.8884}}
\put(1062.000,3159.000){\arc{2100.000}{3.1416}{6.2832}}
\put(1812,3009){\ellipse{600}{300}}
\put(312,3009){\ellipse{600}{300}}
\path(612,3009)(612,3159)
\path(12,3009)(12,3159)
\path(1512,3009)(1512,3159)
\path(2112,3009)(2112,3159)
\thicklines
\path(1812,1059)(1812,2859)
\thinlines
\path(4512,2259)(6612,2259)
\blacken\thicklines
\path(6492.000,2229.000)(6612.000,2259.000)(6492.000,2289.000)(6492.000,2229.000)
\path(7737,984)(7737,2784)
\thinlines
\path(4512,8409)(6612,8409)
\blacken\thicklines
\path(6492.000,8379.000)(6612.000,8409.000)(6492.000,8439.000)(6492.000,8379.000)
\path(7737,7134)(7737,8934)
\thinlines
\path(837,4659)(837,6759)
\blacken\thicklines
\path(867.000,6639.000)(837.000,6759.000)(807.000,6639.000)(867.000,6639.000)
\thinlines
\path(7737,4059)(7737,6159)
\blacken\thicklines
\path(7767.000,6039.000)(7737.000,6159.000)(7707.000,6039.000)(7767.000,6039.000)
\put(1887,7959){\makebox(0,0)[lb]{\smash{{{\SetFigFont{8}{9.6}{\rmdefault}{\mddefault}{\updefault}$\alpha^{-1}$}}}}}
\put(1887,1809){\makebox(0,0)[lb]{\smash{{{\SetFigFont{8}{9.6}{\rmdefault}{\mddefault}{\updefault}$1$}}}}}
\put(7812,1809){\makebox(0,0)[lb]{\smash{{{\SetFigFont{8}{9.6}{\rmdefault}{\mddefault}{\updefault}$1$}}}}}
\put(7812,8034){\makebox(0,0)[lb]{\smash{{{\SetFigFont{8}{9.6}{\rmdefault}{\mddefault}{\updefault}$\alpha$}}}}}
\put(5262,8634){\makebox(0,0)[lb]{\smash{{{\SetFigFont{8}{9.6}{\rmdefault}{\mddefault}{\updefault}$k$}}}}}
\put(5262,1734){\makebox(0,0)[lb]{\smash{{{\SetFigFont{8}{9.6}{\rmdefault}{\mddefault}{\updefault}$M$}}}}}
\put(8037,4884){\makebox(0,0)[lb]{\smash{{{\SetFigFont{8}{9.6}{\rmdefault}{\mddefault}{\updefault}$\theta$}}}}}
\put(1137,5259){\makebox(0,0)[lb]{\smash{{{\SetFigFont{8}{9.6}{\rmdefault}{\mddefault}{\updefault}$\theta$}}}}}
\end{picture}
}
\end{center}
\vspace{0.3cm}

The result now follows because the map along the top is
$t_{\varphi(\alpha)}$ and the bottom map is the
identity.
\end{proof}

Now we investigate how the involution interacts with the monoidal
structure and the crossing.

\begin{prop}\label{defn:p}
Let $U\in\cA_\alpha ,V\in\cA_\beta$ and $\gamma\in\pi$. There are natural isomorphisms
\[
\xymatrix{p_{U,V}: (U*V)^*
  \ar[r]^-{\simeq} & V^**U^* }
\]
and 
\[
\xymatrix{c_{\gamma,U}: (\varphi(\gamma)U)^*
  \ar[r]^-{\simeq} & \varphi(\gamma)U^* }
\]
and moreover the $p_{U,V}$ and $c_{\gamma,U}$ commute, that is, the
following diagram commutes:
\[
\xymatrix{ (\vg U*\vg V)^* \ar[d]_= \ar[r]^p & (\vg V)^**(\vg U)^*
  \ar[r]^{c*c} &\vg V^**\vg U^* \ar[d]^=\\ 
  (\vg(U*V))^*\ar[r]^c & \vg(U*V)^* \ar[r]^{\vg(p)} & \vg (V^**U^*)}
\] 

\end{prop}

\begin{proof} 
  Let $p\colon P\rightarrow X$ be the map giving the monoidal product
  $\cA_\alpha \otimes \cA_\beta \rightarrow \cA_{\alpha\beta}$.
  Notice that $\rotate{\xyreflect{p}}$ induces the monoidal product
  $\cA_{\alpha^{-1}} \otimes \cA_{\beta^{-1}} \rightarrow
  \cA_{\beta^{-1}\alpha^{-1}}$.  By property S-I of self-dual
  $\cS_X$-structures we have natural isomorphisms
\[
\cA(U*V,-)\nattrans \cA\otimes\cA(U\otimes V,\Delta -)
\]
i.e. isomorphisms
\[
q_{U,V}\colon \ohom \circ F(p)(U\otimes V) \rightarrow F(\xyreflect
p)^\vee \circ \ohom (U\otimes V)
\]
Applying the functor $\cJ$ to these and composing with
$t\colon \cJ\circ F(\xyreflect
p)^\vee \rightarrow F(\rotate{\xyreflect{p}})\circ \cJ$ we obtain
\[
\xymatrix{
\cJ \circ \ohom \circ *(U\otimes V)\ar[r]^{J(q_{U,V})} &
\cJ\circ F(\xyreflect p)^\vee \circ \ohom (U\otimes V)
\ar[r]^{t_{\ohom (U\otimes V)}} & {*}\circ \cJ \circ \ohom
(U\otimes V)  
} 
\]
which define a collection of natural transformations
\[
p_{U,V}\colon (U*V)^* \rightarrow V^* * U^*.
\] 

To construct $c_{\gamma,U}$ proceed as above but starting with the
natural transformations
\[
w_{\gamma,U}: \cA(\vg U,-) \nattrans \cA(U,\vgi -)
\]
again coming from property S-I of a
self-dual $\cS_X$-structure.

To show $c$ and $p$ commute we claim the following diagram commutes
where composition has been omitted 
and a superscript 2 refers to tensor product.

{\SMALL 
\[
\xymatrix{
&&{*}\cJ^2{\vgi^\vee}^2 \cI^2\cJ^2\ohom^2\ar[r]_{*(k^2)}&{*}\vg^2\cJ^2\ohom^2 \ar[ddddddd]^=\\
&{*}\cJ^2\ohom^2\vg^2 \ar@/^3pc/[urr]^{c*c} \ar[r]& {*}\cJ^2{\vgi^\vee}^2\ohom^2\ar[u]_{*\cJ^2{\vgi^\vee}^2(N)}&\\
&\cJ\Delta^\vee\cI^2\cJ^2\ohom^2\vg^2 \ar[u]_k&\cJ\Delta^\vee\cI^2\cJ^2{\vgi^\vee}^2\ohom^2 \ar[u]_k&\\
\cJ\ohom {*}\vg^2\ar[uur]^p \ar[r] \ar[d]_=&\cJ\Delta^\vee\ohom^2\vg^2 \ar[u]_{\cJ\Delta^\vee(N)} \ar[r] &\cJ\Delta^\vee{\vgi^\vee}^2\ohom^2 \ar[d]_=\ar[u]_{\cJ\Delta^\vee(N)}&\\
\cJ\ohom\vg{*}\ar[r] \ar[ddr]_c&\cJ\vgi^\vee\ohom {*}\ar[d]^{\cJ\vgi^\vee(N)}\ar[r]&\cJ\vgi^\vee\Delta^\vee\ohom^2\ar[d]^{\cJ\vgi^\vee(N)}&\\
&\cJ\vgi^\vee\cI\cJ\ohom {*}\ar[d]^k&\cJ\vgi^\vee\cI\cJ\Delta^\vee\ohom^2\ar[d]^k&\\
&\vg\cJ\ohom {*}\ar@/_3pc/[drr]_{\vg(p)} \ar[r]&\vg\cJ\Delta^\vee\ohom^2\ar[d]^{\vg\cJ\Delta^\vee(N)}&\\
&&\vg\cJ\Delta^\vee\cI^2\cJ^2\ohom^2\ar[r]^{\vg(k)}&\vg{*}\cJ^2\ohom^2\\
}
\]
}

Let $m\colon P \rightarrow X$ and $i\colon I\rightarrow X$ induce the
monoidal structure and crossing $\varphi(\gamma)$ on $\cA$. The surface diagram
for the crossing in section \ref{section:balanced} induces an equality $m\circ (i\sqcup i) =
i\circ m$, thus the natural transformations in S-I of the self-duality
has two decompositions, which result in the following commutative diagram.
\[
\xymatrix{
\cA(\vg U*\vg V,-)\ar[d]_= \ar[r] & \cA^2(\vg V\otimes \vg U,\Delta -)\ar[r] & \cA^2(V\otimes U,\vgi^2\Delta -)\ar[d]^=\\
\cA(\vg(U*V),-)\ar[r] & \cA(U*V,\vgi -) \ar[r] & \cA^2(V\otimes U,\Delta\vgi -)
}
\]
Applying $\cJ$ shows that the central horizontal slice of the main
diagram above commutes. The two large central rectangles
commute by naturality, and the triangles and sections with curved
arrow commute by the
definitions of $p$ and $c$. Finally the far right part of the diagram
follows from identities
in the mapping class group, similar to those in Lemma \ref{lem:jumpJ}.

\end{proof}

After these preliminaries we now construct the duality.  Notice that 
for $U\in\cA_\alpha, V\in\cA_\beta$ and $W\in\cA_\gamma$ with
$\alpha = \beta \gamma$ there are natural isomorphisms

\[
\xymatrix{
\cA_\alpha (U*W,V)  \simeq  \la (U*W)^* , V \ra
\ar[r]^-{\la p , id \ra}  &    
\langle W^* *U^*, V \rangle\simeq  \la W^*,U^* *V\ra
\simeq  \cA_\gamma (W,U^* * V)
}
\]

Use these to define $b_U$ and $d_U$ by
\[
\cA_\alpha (U,U) \simeq \cA_{1}(\monunit, U^* *U) \;\; id \mapsto b_U
\]
\[
\cA_{\alpha^{-1}}(U,U) \simeq \cA_{1}(U*U^*,\monunit) \;\; id \mapsto d_U
\]

\begin{prop}\label{prop:dual}
$b_U$ and $d_U$ provide $\cA$ with a right duality.
\end{prop}

\begin{proof}
  Since $\cA$ is a monoidal category it can be thought of as a
  2-category with one object and it follows from the definition of
  adjoints in 2-categories that $b_U$ and $d_U$ give a right duality
  iff $U$ is right adjoint to $U^*$ in this category. It follows from
  the theory of adjoints in 2-categories (see for example
  \cite{Gray:FormalCategoryTheory} page 158) that $U$ is left adjoint to $U^*$ iff there
  are natural equivalences
\[
\cA_\alpha (U,V*W) \simeq \cA_\gamma (V^* * U,W)
\]
and
\[
(U^**V)*U \simeq U^**(V*U)
\]

The result now follows by the natural isomorphisms above  and associativity.
\end{proof}

We now prove that the $c_{\gamma,U}$'s defined above
satisfy the definition of lax tortile $\pi$-category, but first we
need one result about the crossing.

\begin{lem}\label{lem:top}
The following diagram commutes
\[
\xymatrix{
  \cA(\vg U,\vg U) \ar[d]_\simeq  && \cA(U,U) \ar[ll]_\vg \ar[d]^\simeq\\
  \la (\vg U)^*,\vg U\ra \ar[r]^{\la c_{\gamma,U},id\ra} & \la \vg
  U^*,\vg U\ra \ar[r]^-{=} &\la,U^*,U\ra \\ 
  }
\]
\end{lem}

\begin{proof}
  Recalling the definition of $w_{\gamma,U}$ from Proposition \ref{defn:p},
  since $\vg$ is induced from a straight cylinder, $w_{\gamma,U}$ acts
  as the functor $\vgi$ by the definition of lax
  self-dual $\cS_X$-structure.  Now consider the diagram below.
\[
\xymatrix{\ohom \circ \vg (U)\ar[d]_{w_{\gamma,U}} \ar[r]^N & \cI
  \circ \cJ \circ \ohom \circ \vg (U) \ar[d]^{I\circ J
  (w_{\gamma,U})}\\ 
  \vgi^\vee \circ \ohom (U) \ar[r]^N \ar[d]_{\vgi^\vee(N)} &
  \cI\circ\cJ\circ \vgi^\vee\circ\ohom(U)
  \ar[d]^{\cI\circ\cJ\circ\vgi^\vee(N)}\\ 
  \vgi^\vee \circ \cI\circ\cJ\circ\ohom(U) \ar[d]_= \ar[r]^N &
  \cI\circ\cJ\circ\vgi^\vee\circ\cI\circ\cJ\circ\ohom(U)
  \ar[dl]^{\cI(k)} \\ 
  \cI \circ \vg \circ \cJ \circ \ohom (U) 
}
\] 

To see
that this diagram commutes: the two top sections follow immediately by naturality of
$N$; the bottom triangle is a consequence of the mapping class group
identity below.

\vspace{0.3cm}
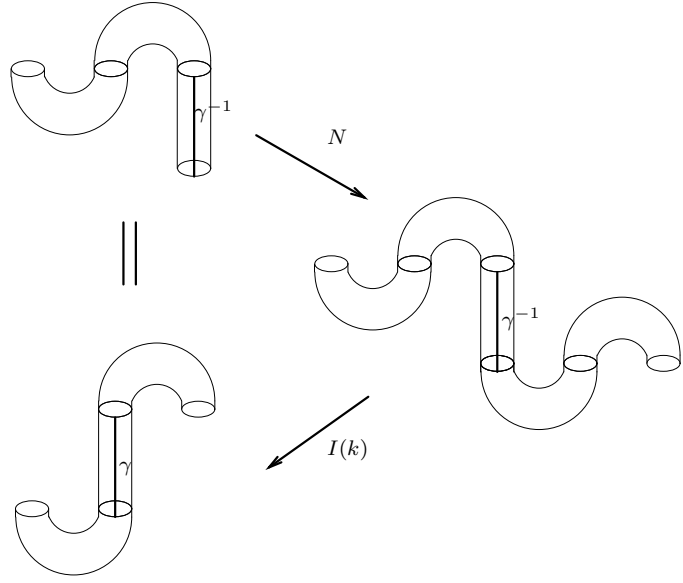
\begin{center}
\setlength{\unitlength}{0.00029167in}
\begingroup\makeatletter\ifx\SetFigFont\undefined%
\gdef\SetFigFont#1#2#3#4#5{%
  \reset@font\fontsize{#1}{#2pt}%
  \fontfamily{#3}\fontseries{#4}\fontshape{#5}%
  \selectfont}%
\fi\endgroup%
{\renewcommand{\dashlinestretch}{30}
\begin{picture}(12099,10383)(0,-10)
\put(1062.000,9196.500){\arc{975.000}{0.3948}{2.7468}}
\put(1062.000,9009.000){\arc{2100.000}{6.2832}{9.4248}}
\put(1812,9159){\ellipse{600}{300}}
\put(312,9159){\ellipse{600}{300}}
\path(612,9159)(612,9009)
\path(12,9159)(12,9009)
\path(1512,9159)(1512,9009)
\path(2112,9159)(2112,9009)
\put(2562.000,9121.500){\arc{975.000}{3.5364}{5.8884}}
\put(2562.000,9309.000){\arc{2100.000}{3.1416}{6.2832}}
\put(3312,9159){\ellipse{600}{300}}
\put(1812,9159){\ellipse{600}{300}}
\path(2112,9159)(2112,9309)
\path(1512,9159)(1512,9309)
\path(3012,9159)(3012,9309)
\path(3612,9159)(3612,9309)
\put(3312,9159){\ellipse{600}{300}}
\put(3312,7359){\ellipse{600}{300}}
\path(3012,7359)(3012,9159)
\path(3612,9159)(3612,7359)
\put(1137.000,1246.500){\arc{975.000}{0.3948}{2.7468}}
\put(1137.000,1059.000){\arc{2100.000}{6.2832}{9.4248}}
\put(1887,1209){\ellipse{600}{300}}
\put(387,1209){\ellipse{600}{300}}
\path(687,1209)(687,1059)
\path(87,1209)(87,1059)
\path(1587,1209)(1587,1059)
\path(2187,1209)(2187,1059)
\put(1887,3009){\ellipse{600}{300}}
\put(1887,1209){\ellipse{600}{300}}
\path(1587,1209)(1587,3009)
\path(2187,3009)(2187,1209)
\put(2637.000,2971.500){\arc{975.000}{3.5364}{5.8884}}
\put(2637.000,3159.000){\arc{2100.000}{3.1416}{6.2832}}
\put(3387,3009){\ellipse{600}{300}}
\put(1887,3009){\ellipse{600}{300}}
\path(2187,3009)(2187,3159)
\path(1587,3009)(1587,3159)
\path(3087,3009)(3087,3159)
\path(3687,3009)(3687,3159)
\put(9537.000,3871.500){\arc{975.000}{0.3948}{2.7468}}
\put(9537.000,3684.000){\arc{2100.000}{6.2832}{9.4248}}
\put(10287,3834){\ellipse{600}{300}}
\put(8787,3834){\ellipse{600}{300}}
\path(9087,3834)(9087,3684)
\path(8487,3834)(8487,3684)
\path(9987,3834)(9987,3684)
\path(10587,3834)(10587,3684)
\put(11037.000,3796.500){\arc{975.000}{3.5364}{5.8884}}
\put(11037.000,3984.000){\arc{2100.000}{3.1416}{6.2832}}
\put(11787,3834){\ellipse{600}{300}}
\put(10287,3834){\ellipse{600}{300}}
\path(10587,3834)(10587,3984)
\path(9987,3834)(9987,3984)
\path(11487,3834)(11487,3984)
\path(12087,3834)(12087,3984)
\put(8037.000,5596.500){\arc{975.000}{3.5364}{5.8884}}
\put(8037.000,5784.000){\arc{2100.000}{3.1416}{6.2832}}
\put(8787,5634){\ellipse{600}{300}}
\put(7287,5634){\ellipse{600}{300}}
\path(7587,5634)(7587,5784)
\path(6987,5634)(6987,5784)
\path(8487,5634)(8487,5784)
\path(9087,5634)(9087,5784)
\put(6537.000,5671.500){\arc{975.000}{0.3948}{2.7468}}
\put(6537.000,5484.000){\arc{2100.000}{6.2832}{9.4248}}
\put(7287,5634){\ellipse{600}{300}}
\put(5787,5634){\ellipse{600}{300}}
\path(6087,5634)(6087,5484)
\path(5487,5634)(5487,5484)
\path(6987,5634)(6987,5484)
\path(7587,5634)(7587,5484)
\thicklines
\path(4823.165,2146.686)(4662.000,1959.000)(4892.527,2048.763)
\path(4662,1959)(6462,3234)
\put(5712,2184){\makebox(0,0)[lb]{\smash{{{\SetFigFont{8}{9.6}{\rmdefault}{\mddefault}{\updefault}$I(k)$}}}}}
\thinlines
\put(8787,3834){\ellipse{600}{300}}
\put(8787,5634){\ellipse{600}{300}}
\put(8787,3834){\ellipse{600}{300}}
\thicklines
\path(2037,6384)(2037,5259)
\path(2262,6384)(2262,5259)
\thinlines
\path(8487,3834)(8487,3684)
\path(9087,3834)(9087,3684)
\thicklines
\path(8787,3684)(8787,5484)
\thinlines
\path(8487,3834)(8487,5634)
\path(9087,5634)(9087,3834)
\thicklines
\path(3312,7209)(3312,9009)
\path(1887,1059)(1887,2859)
\path(6149.132,6901.962)(6387.000,6834.000)(6209.099,7005.905)
\path(6387,6834)(4437,7959)
\put(8862,4509){\makebox(0,0)[lb]{\smash{{{\SetFigFont{8}{9.6}{\rmdefault}{\mddefault}{\updefault}$\gamma^{-1}$}}}}}
\put(3312,8259){\makebox(0,0)[lb]{\smash{{{\SetFigFont{8}{9.6}{\rmdefault}{\mddefault}{\updefault}$\gamma^{-1}$}}}}}
\put(1962,1884){\makebox(0,0)[lb]{\smash{{{\SetFigFont{8}{9.6}{\rmdefault}{\mddefault}{\updefault}$\gamma$}}}}}
\put(5712,7809){\makebox(0,0)[lb]{\smash{{{\SetFigFont{8}{9.6}{\rmdefault}{\mddefault}{\updefault}$N$}}}}}
\end{picture}
}
\end{center}
\vspace{0.3cm}

Now following
the diagram around the right hand path from the top left hand corner gives the map
$\cI(c_{\gamma,U})\circ N$ by definition of $c_{\gamma,U}$ and  the
left hand side  is the map $\vgi^\vee(N)\circ
\vgi$, by the remark made above. Evaluating the diagram of functors on
$\vg U$ gives the desired result.
\end{proof}

\begin{prop}\label{prop:tortilepicat}
$\cA$ is a lax tortile $\pi$-category.
\end{prop}

\begin{proof}
Diagram (1) of Definition
\ref{defn:tortile} is simply naturality of the $c_{\alpha,U}$.

To prove diagram (2) of 
\ref{defn:tortile} we claim the following diagram commutes.

\[
\xymatrix{
J\circ \ohom \circ \va (U) \ar[d]_{\cJ\ohom(\theta_U)}
\ar[rr]^{J(w_{\alpha,U})} && J\circ \vai^\vee \circ \ohom(U)
\ar[rr]^{t_{\va,\ohom (U)}} && \va \circ \cJ\circ\ohom(U)
\ar[dll]^{\theta_{\va\cJ\ohom(U)}}\\ 
J\circ \ohom (U) \ar[rr]^= && J\circ \ohom (U)
\ar[u]^{\cJ\theta^\vee_{\ohom(U)}} 
}
\]

The top line of this diagram is
$c_{\alpha,U}$ by definition, the left map $\theta_U^*$ and the right
map $\theta_{\va U^*}$.

By the consequence S-II of self-duality
  the natural transformations $w_{\alpha,U}$  satisfy 
\[
\xymatrix{\cA(\va U,-) \ar[d]_{\cA(\theta_U,-)} \ar[r]^{w_{\alpha,U}}
  & \cA(U,\vai -) \ar[d]^{\cA(U,\theta_{\vai -})}\\ 
  \cA(U,-) \ar[r]_{id} & \cA(U,-)}
\]
This is obtained from S-II by taking the diffeomorphism $T$ to be the
inverse twist, which is mapped to the twist by reflection. The bottom
map is the identity since the identity is sent to the identity in the
definition of pseudo 2-natural transformation.  Taking $\cJ$ of this
diagram proves the left square above commutes.

For the triangle on the right pre compose the diagram in Lemma \ref{lem:crossJ} with $\hom$.

Finally, diagram (3) for $b_U$ in definition \ref{defn:tortile} is
shown to commute by considering 
the following diagram.
{\SMALL 
\[
\xymatrix{\cA(\vg U,\vg U) \ar[d]_\simeq  &&& \cA(U,U) \ar[lll]_\vg \ar[d]^\simeq\\
\la (\vg U)^*,\vg U\ra \ar[d]_\simeq \ar[rr]^{\la c_{\gamma,U},id\ra} && \la \vg U^*,\vg U\ra \ar[r]^= &\la U^*,U\ra \ar[d]^\simeq\\
\la (\vg U^**\bf1)^*,\vg U\ra \ar[d]_{\la p, id\ra} \ar[r]^= & \la (\vg(U*\bf1))^*,\vg U\ra \ar[r]^{\la c_{\gamma,U*\bf1},id\ra}  & \la \vg(U*\bf1)^*,\vg U\ra \ar[r]^= \ar[d]^{\la \vg p,id\ra} & \la (U*\bf1)^*,U\ra \ar[d]^{\la p,id\ra}\\
\la (\vg \bf1)^* * (\vg U)^*,\vg U\ra \ar[d]_\simeq \ar[rr]^{\la c_{\gamma,\bf1}*c_{\gamma,U},id\ra} && \la \vg \bf1^* *\vg U^*,\vg U\ra \ar[r]^= & \la \bf1^* * U^*,U\ra \ar[d]^\simeq\\
\la (\vg \bf1)^*,(\vg U)^**\vg U\ra \ar[rr]^{\la c_{\gamma,\bf1},c_{\gamma,U}*id\ra} \ar[d]_\simeq && \la \vg \bf1^*, \vg U^* *\vg U\ra \ar[r]^= & \la \bf1^*, U^* *U\ra \ar[d]^\simeq\\
\cA(\bf1,(\vg U)^**\vg U) && \cA(\vg \bf1, \vg U^**\vg U) \ar[ll]_{\cA(id,c_{\gamma,U}*id)} & \cA(\bf1,U^**U) \ar[l]_\vg
}
\]
} 

Notice that the left and right sides of the diagram define $b_{\vg U}$
and $b_U$ respectively as the image of the identity morphism. Thus,
since $\vg$ is a functor, the commutativity of this diagram proves (3)
by mapping the identity
in the top right hand corner to the bottom left hand corner via the
two exterior paths. 

The top and bottom parts of
the diagram commute by lemma \ref{lem:top} and naturality. The left
part of the central section commutes by lemma \ref{defn:p}. The
remaining sections of the diagram commute by naturality.

To complete the proof it
remains to check diagram (3) for $d_U$. The arguments are essentially
the same as those for the $b_U$ case, so we simply give the
corresponding diagram without further comment.

\hspace*{-3cm}
{\SMALL
\[
\xymatrix{
\cA((\vg U)^*,(\vg U)^*) \ar[d]_\simeq && \cA(\vg U^*,\vg U^*) \ar[d]_\simeq \ar[ll]_{\cA(c_{\gamma,U},c_{\gamma,U}^{-1})} && \cA(U^*,U^*) \ar[d]^\simeq \ar[ll]_\vg\\
\la (\vg U)^**,(\vg U)^*\ra \ar[d]_\simeq \ar[rr]^{\la {c_{\gamma,U}^*}^{-1},c_{\gamma,U}\ra} && \la (\vg U^*)^*,\vg U^*\ra \ar[d]_\simeq \ar[r]^{\la c_{\gamma,U},id\ra} & \la \vg U^{**}, \vg U^*\ra \ar[r]^= & \la U^{**},U^*\ra \ar[d]^\simeq\\
\la (\vg U)^{**}, (\vg U)^* *\bf1\ra \ar[d]_\simeq \ar[rr]^{\la c_{\gamma,U}^{*^{-1}},c_{\gamma,U}*id\ra} && \la (\vg U^*)^*,\vg U^**\bf1\ra \ar[r]^{\la c_{\gamma,U},id\ra} & \la \vg U^{**}, \vg U^**\bf1\ra \ar[r]^= & \la U^{**},U^**\bf1\ra \ar[d]^\simeq\\
\la (\vg U)^{**}*(\vg U)^*,\bf1 \ra \ar[rr]^{\la c_{\gamma,U}^{*^{-1}}*id,id \ra} && \la (\vg U^*)^**(\vg U)^*,\bf1\ra \ar[r]^{\la c_{\gamma,U}*c_{\gamma,U},id\ra} & \la \vg U^{**}*\vg U^*,\bf1\ra \ar[r]^= & \la U^{**}*U^*,1\ra\\
\la (\vg U*(\vg U)^*)^*,\bf1\ra \ar[d]_\simeq \ar[u]^{\la p,id\ra} \ar[rr]^{\la (id*c_{\gamma,U}^{-1})^*,id\ra} && \la (\vg U*\vg U^*)^*,\bf1\ra \ar[u]^{\la p,id\ra} \ar[r]^{\la c_{\gamma,U*U^*},id\ra} \ar[d]_\simeq & \la \vg(U*U^*)^*,\bf1\ra \ar[u]^{\la \vg p,id\ra} \ar[r]^= & \la (U*U^*)^*,\bf1\ra \ar[u]_{\la p,id\ra} \ar[d]^\simeq\\
\cA(\vg U*(\vg U)^*,\bf1) && \cA(\vg U*\vg U^*,\bf1) \ar[ll]_{\cA(id*c_{\gamma,U},id)} && \cA(U*U^*,1) \ar[ll]_\vg
}
\]
}
This finishes the proof.
\end{proof}

Now let $X$ be an arbitrary space and claim that $\cA_1$ is a lax
tortile category with $\pi_2(X)$-action. Methods similar to those used for
$X=K(\pi,1)$ above, produce a duality on $\cA_1$. We need to relate
this to the
$\pi_2(X)$-action.

\begin{lem}
For all $U\in \cA_1$ and $g\in \pi_2X$ there are natural isomorphisms
\[
h_{g,U}:(\rho(g)U)^* \rightarrow \rho(g^{-1})U^*
\] 
\end{lem}
\begin{proof}
  The maps $i_g:I\rightarrow X$ giving $\pi_2(X)$-action
  satisfy $\xyreflect {i_g} = i_{g^{-1}}$ so, by property S-I, the
  self-duality induces natural isomorphisms
\[
v_{g,U}:\cA(\rho(g)U,-)\nattrans \cA(U,\rho(g^{-1})-)
\]
Applying the functor $\cJ$ and composing with the equivalence $t$ of \ref{lem:jumpJ} we get
\[
\xymatrix{
\cJ\circ\ohom\circ\rho(g)(U) \ar[r]^{J(v_{g,U})} &\cJ\circ\rho(g^{-1})^\vee \circ \ohom (U) \ar[r]^{t_{\rho(g),\ohom(U)}} &\rho(g^{-1})\circ \cJ \circ \ohom (U)
}
\]
giving the desired natural isomorphisms.
\end{proof}

\begin{prop}\label{prop:tortileGaction}
$\cA_1$ is a lax tortile category with  $\pi_2(X)$-action
\end{prop}

\begin{proof}
Diagram (1) of definition \ref{defn:ribbonG} is naturality of the
$h_{g,U}$'s.  To prove the diagrams in (2) for $b_U$ and $d_U$,
proceed as in the proof of Proposition \ref{prop:tortilepicat} replacing $c_{\gamma,
  U}$ with $h_{g,U}$. 
\end{proof}

Combining Propositions \ref{prop:dual}, \ref{prop:tortilepicat} and \ref{prop:tortileGaction}
with the results of section \ref{section:balanced} completes the proof of
Theorem \ref{thm:main}.

\appendix

\section{Balanced and tortile structures}\label{app:defns}

This appendix contains the variants of balanced and tortile categories
that we need. The similarity of our lax tortile $\pi$-category to
Turaev's ribbon $\pi$-category should be noted. This is not an
accident as both grew out of \cite{Turaev:HomotopyFieldTheory2D}, however because we  do not
have duals from the outset we must
define balanced categories first and add duals later. Also our structures
are lax versions of those defined by Turaev in \cite{Turaev:HomotopyFieldTheory3D}.

\begin{defn}
Let $\pi$ be a discrete group. A monoidal category $(\cA, *, \monunit)$
is said to be $\pi$-{\em graded} if it splits as a disjoint union of
full subcategories $\cA
= \bigsqcup_{\alpha\in \pi} \cA_\alpha$ such that
 objects belong to a unique grading,
 there are no morphisms between objects of different grading
and the monoidal structure adds gradings (in the group $\pi$).
\end{defn}

We
denote the set of invertible functors on $\cA$ by $\oAut (\cA)$.  A
functor is {\em monoidal} if the monoidal product, unit and all
structure morphisms are preserved.

\begin{defn}
A $\pi$-graded monoidal category is said to be {\em crossed} if it
comes equipped with a group homomorphism $\varphi\colon \pi \rightarrow
\oAut (\cA)$ with $\varphi(\gamma)\colon \cA_\beta \rightarrow
\cA_{\gamma\beta\gamma^{-1}}$ such that each $\varphi(\gamma)$ is a
monoidal functor.
\end{defn}

\begin{defn}\label{defn:braid}
A {\em braiding} on a crossed $\pi$-graded monoidal category $\cA$ is a
collection of isomorphism $s_{U,V}:U*V\rightarrow \varphi(\alpha)V*U$
for objects $U$ of $\cA_\alpha$, $V$ of $\cA_\beta$, such that 
\begin{enumerate}
\item For $U\in \cA_\alpha, V\in \cA_\beta$ and $W\in\cA_\gamma$ the following two diagrams commute: 
\[
\xymatrix{U*V*W  \ar[rr]^{s_{U*V,W}} \ar[dr]_{id*s_{V,W}} & & \varphi (\alpha\beta)W*U*V
  \\ & U*\varphi(\beta)W*V \ar[ur]_{s_{U,\varphi(\beta)W}*id}}
\]
\[
\xymatrix{U*V*W  \ar[rr]^{s_{U,V*W}} \ar[dr]_{s_{U,V}*id} & & \varphi (\alpha)V*\varphi (\alpha)W*U
  \\ & \varphi(\alpha)V*U*W \ar[ur]_{id*s_{U,W}}} 
\]
(associativity morphisms have been suppressed; insert them and one
gets the hexagons expected by the reader).
\item for $U,U^\prime \in \cA_\alpha$ and $V,V^\prime \in \cA_\beta$
  and morphisms $f\colon U \rightarrow U^\prime$ and $g\colon V
  \rightarrow V^\prime$ the following diagram commutes
\[
\xymatrix{U*V \ar[rr]^{f*g} \ar[dd]_{s_{U,V}} & & U^\prime * V^\prime
  \ar[dd]^{s_{U^\prime,V^\prime}}\\ 
& & \\
 \varphi(\alpha)V*U
  \ar[rr]^{\varphi(\alpha)g*f} & & \varphi(\alpha)V^\prime *U^\prime}
\]
\item for $U\in\cA_\alpha$ and $V\in \cA_\beta$ and $\gamma\in \pi$,
\[
s_{\varphi(\gamma)U, \varphi(\gamma)V} = \varphi(\gamma)(s_{U,V}).
\]
\end{enumerate}
\end{defn}

\begin{defn}\label{defn:twist}
A {\em twist} on a braided crossed $\pi$-graded monoidal category $\cA$ is a
 collection of
isomorphisms $\theta_U\colon U \rightarrow \varphi(\alpha)U$ for $U\in
\cA_\alpha$ such that
\begin{enumerate}
\item $\theta_\monunit=\oid_\monunit$
\item for $U\in \cA_\alpha$ and $V\in \cA_\beta$ 
\[
\theta_{U*V}=s_{\varphi(\alpha\beta)V, \varphi(\alpha)U}\circ 
(\theta_{\varphi(\alpha)V}*\theta_U)\circ s_{U,V}
\]
\item for $U,V\in \cA_\alpha$ and  morphism $f\colon U \rightarrow V$
  the following diagram commutes 
\[
\xymatrix{U \ar[rr]^{f} \ar[dd]_{\theta_{U}} & & V
  \ar[dd]^{\theta_{V}}\\ 
& & \\
 \varphi(\alpha)U
  \ar[rr]^{\varphi(\alpha)f} & & \varphi(\alpha)V}
\]
\item for $\gamma\in\pi$ and $V\in \cA$,
\[
\theta_{\varphi(\gamma)V}= \varphi(\gamma)(\theta_V)
\]
\end{enumerate}
\end{defn}

\begin{defn}
A {\em balanced $\pi$-category} is a braided crossed $\pi$-graded monoidal
category with twist.
\end{defn}

When $\pi=\{1\}$ the usual notion of a balanced category is obtained.
Now we introduce duality into balanced $\pi$-categories. A {\em right
  duality} assigns to an object $U\in\cA_\alpha$ an object $U^*\in
\cA_{\alpha^{-1}}$ and morphisms
\[
b_U\colon \monunit \rightarrow U^* * U \;\;\;\; d_U\colon U*U^*
\rightarrow \monunit
\]
such that the following compositions are the identity.
\[
\xymatrix{U^*\simeq \monunit * U^* \ar[r]^(0.3){b_U*id} & (U^**U)*U^*\simeq
  U^**(U*U^*) \ar[r]^(0.7){id*d_U} & U^**\monunit \simeq U^*
}
\]
\[
\xymatrix{U\simeq U*\monunit  \ar[r]^(0.3){id*b_U} & U*(U^**U)\simeq
  (U*U^*)*U  \ar[r]^(0.65){d_U*id} & \monunit *U \simeq U
}
\]
The above duality is also called rigid duality, and $U^*$ is a rigid dual, to distinguish it from weaker versions. By duality we shall always mean rigid duality.

\begin{defn}\label{defn:tortile}
A {\em lax tortile $\pi$-category} is a balanced $\pi$-category with right
duality and a collection of isomorphisms $c_{\gamma,U}:(\varphi(\gamma)U)^*\rightarrow
\varphi(\gamma)U^*$ for $U\in \cA$ and $\gamma \in \pi$ such that
\begin{enumerate}
\item For $U\in \cA_\alpha$ and $f:U\rightarrow V$ the following diagram commutes
\[
\xymatrix{(\varphi(\gamma)U)^* \ar[rr]^{c_{\gamma,U}}  & & \varphi(\gamma)U^*
  \\ 
& & \\
 (\varphi(\gamma)V)^*\ar[uu]^{(\varphi(\gamma)f)^*}
  \ar[rr]^{c_{\gamma,V}} & & \varphi(\gamma)V^*\ar[uu]_{\varphi(\gamma)f^*}}
\]
\item 
For $U\in \cA_\alpha$  the following diagram commutes 
\[
\xymatrix{(\varphi(\alpha)U)^* \ar[rr]^{c_{\alpha,U}} \ar[dr]_{(\theta_U)^*} & & \varphi(\alpha)U^*
\ar[dl]^{\theta_{\varphi(\alpha)U^*}}
\\
& U^*}
\]
\item
For $U\in \cA_\alpha$ the following diagram commutes
\[
\xymatrix{(\vg U)^**\vg U \ar[rr]^{c_{\gamma,U}*id} & &  \vg U^**\vg U\\
& \monunit \ar[ul]^{b_{\vg U}} \ar[ur]_{\vg (b_u)}&}
\]
and similarly for $d_u$.
\end{enumerate}
\end{defn}

Note that a \emph{strict} tortile $\pi$-category is one for which all
$c_{\gamma,U}$'s are identities, and is the same thing as a ribbon
crossed $\pi$-category (satisfying $\theta_{\monunit} = id_{\monunit}$).

Now we introduce an ungraded version with an action of a
group $G$.

\begin{defn}\label{defn:balancedG}
A {\em balanced  category with $G$-action}
is a balanced  category  $(\cA, * , \monunit, s,\theta)$  together
with a homomorphism 
$\rho:G\rightarrow \oAut(\cA)$  such that
\begin{enumerate}
\item $\rho (g)U*V = \rho (g) (U*V) = U*\rho (g)V$
\item for $f\colon U\rightarrow U^\prime$ and $h\colon V\rightarrow
  V^\prime$ we have $\rho (g)f*h = \rho (g) (f*h) = f*\rho (g)h$
\item $\rho (g)(a_{U,V,W})=a_{\rho (g)U,V,W}=a_{U,\rho (g)V,W}=a_{U,V,\rho (g)W}$
\item $\rho (g)(r_U)=r_{\rho (g)U}$ and $\rho (g)(l_U)=l_{\rho (g)U}$ 
\item $s_{\rho(g)U, V} = \rho(g)(s_{U,V}) = s_{U,\rho(g)V}$
\item $\theta_{\rho(g)U} = \rho(g)(\theta_U)$
\end{enumerate}
\end{defn}

Note that $\rho$ factors through the centre of $ \oAut(\cA)$ and on
objects $\rho (g)$ is determined by $\rho(g)\monunit$ since
$\rho(g)(U) = (\rho(g)\monunit)*U$. Incorporating duals we have the
following definition.

\begin{defn}\label{defn:ribbonG}
A {\em tortile  category with lax $G$-action} is a balanced  category with
$G$-action with right
duality and a collection of isomorphisms  $h_{g,U}: (\rho (g) U)^* \rightarrow \rho (g^{-1}) (U^*)$
such that the following diagrams commute
\begin{enumerate}
\item For $U,V\in \cA$ and $f:U\rightarrow V$,
\[
\xymatrix{(\rho(g)U)^* \ar[r]^{h_U}  & \rgi U^* \\
(\rg V)^* \ar[u]^{(\rho(g)f)^*} \ar[r]_{h_V} & \rgi V^* \ar[u]_{\rgi f^*}}
\]
\item
\[
\xymatrix{(\rho(g)U)^* *\rho(g)U \ar[rr]^{h_{g,U}*id}  & & U^* *U  \\
& \monunit \ar[ul]^{b_{\rho(g)U}} \ar[ur]_{b_U}&}
\]
and similarly for $d_U$.
\end{enumerate} 
\end{defn}

\section{$k$-additive categories}\label{app:add}
Let $k$ be an algebraically closed field. A {\em additive category over $k$}, or {\em $k$-additive category} is a category with the following properties:
\begin{enumerate}
\item
morphism spaces are complex vector spaces and composition is bilinear
\item
there is a finite direct sum on objects
\item
there is a zero object $0$ such that $\ohom(U,0)=\ohom(0,U)=0$ for all $U\in \cA$
\end{enumerate}

$k$-additive functors between $k$-additive categories are additive functors acting linearly on the morphism spaces. The tensor product of two $k$-additive categories is defined in the usual way. For more details on additive categories we refer to \cite{Maclane}.

The {\em idempotent completion} $\hat{\cA}$ of a category $\cA$ has objects pairs $(U,e)$ where $e:U\rightarrow U$ is an idempotent in $\cA$ and morphisms $f:(U,e)\rightarrow (U',e')$ where $f:U\rightarrow U'$ such that $fe=f=e'f$. The idempotent completion of an additive category over $k$ remains an additive category over $k$. Any functor from $\cA$ to an idempotent complete category $\cB$ can be completed to a functor on $\hat{\cA}$; this construction yields a natural equivalence between the categories of such functors.

\begin{defn}
Let $\additive$ be the 2-category with objects idempotent complete additive small categories over $k$, 1-morphisms $k$-additive functors between these and 2-morphisms natural transformations.
\end{defn}

The tensor product of k-additive categories induces a monoidal structure on $\additive$: since $k$ is algebraically closed the tensor product of idempotent complete $k$-additive categories remains idempotent complete. The category $\hat{k}$ of finite dimensional vector spaces over $k$ is the monoidal unit.

The set of $k$-additive functors between $\cA$ and
$\cB$ is denoted $\additive (\cA, \cB)$ and this is again an
idempotent $k$-additive category if both $\cA$ and $\cB$ are small, and $\cB$ is idempotent complete. If a $k$-additive category is small its {\em dual} can be defined as
$\cA^\vee = \additive(\cA, \hat k)$.

An object $X$ of an additive category is {\em simple} if $\cA (X,X)$ is one dimensional and 
$\cA$ is {\em semi-simple} if every object is isomorphic to a finite
sum of simple objects. $\cA$ is {\em Artinian} if there are finitely many isomorphism
classes of simple objects.

We now discuss Tillmann's work on duality in idempotent complete $k$-additive categories. 
A {\em non-degenerate form} consists of a $k$-additive functor $\langle
-,- \rangle \colon \cA \otimes \cB \rightarrow \hat{k}$ and an object
$\sum_{i=1}^n P_i\otimes Q_i$ in $\cB\otimes \cA$ such that the
functors
\begin{eqnarray*}
&&\cI\colon \cA \rightarrow \cB^\vee \text{ by }
Y\mapsto \langle Y,-\rangle\\ 
&&\cJ\colon\cB^\vee\rightarrow \cA\text{ by }
H\mapsto \sum_{i=1}^nH(P_i)\otimes
Q_i 
\end{eqnarray*}
provide an equivalence of categories.
Her methods prove that given a non-degenerate form and letting $X=\sum Q_i$ and
$A=\cA(X,X)$ then the following hold
\begin{enumerate}
\item $\cA$ is equivalent to the category of finite dimensional
  projective $A$-modules
\item $\cA(Y,Z)$ is a finitely generated vector space
\end{enumerate}

She also  defines a  contravariant functor
$(-)^*\colon \cB \rightarrow \cA$ and similarly a contravariant functor
$(-)^*\colon \cA \rightarrow \cB$ with the property that $(-)^*\circ (-)^*$ is
naturally isomorphic to $id_{\cA}$. To define these involutions consider
the contravariant functor
\[
\ohom\colon \cB \rightarrow \cB^\vee \;\;\;\; Y \rightarrow \cB(Y,-)
\]
which is defined since by the above $\cB(Y,Z)$ is a finitely generated
vector space  and hence an object of $\hat k$. Set $(-)^*=\cJ \circ
\ohom$. Tillmann proves that there is a natural equivalence $\cB(-,-) \simeq
\langle -^*,- \rangle$ and her methods also show that 
$\cA$ and $\cB$ are  semi-simple Artinian categories.

If we are given a collection of non-degenerate forms $\langle - , -
\rangle_{\alpha_i}\colon \cA_{\alpha_i} \otimes \cB_{\alpha_i}
\rightarrow \hat{k}$ for $i=1\ldots n$ we can construct a
non-degenerate form
\[
\langle - , -\rangle_\aub \colon \cA_{\alpha_1}\otimes \cdots \otimes
\cA_{\alpha_n}\otimes 
\cB_{\alpha_n}\otimes \cdots \otimes \cB_{\alpha_1} \rightarrow
\hat{k}
\]
and similarly an involution
\[
(-)^*\colon \cB_{\alpha_n}\otimes \cdots \otimes \cB_{\alpha_1}
\rightarrow \cA_{\alpha_1}\otimes \cdots \otimes \cA_{\alpha_n}
\]
which satisfies
$ \cB_{\alpha_n}\otimes \cdots \otimes \cB_{\alpha_1} (-,-) \simeq
\langle (-)^* , -\rangle_\aub .
$

\section*{Acknowledgements}
This work was partially supported by an EPSRC grant. The first author held a postdoctoral fellowship at the Centre de Recherches Math\'{e}matiques in Montr\'{e}al, Canada when it was carried out.

\end{document}